\documentclass[australian,3p,sort&compress,review]{elsarticle}
\usepackage[T1]{fontenc}
\usepackage[utf8]{inputenc}
\usepackage{units}
\usepackage{amsmath}
\usepackage{graphicx}
\usepackage{esint}

\makeatletter

\providecommand{\tabularnewline}{\\}

\journal{}
\usepackage{multirow}
\usepackage{upgreek}
\usepackage{paralist}
\usepackage[section]{placeins}
\usepackage{scalerel}

\newcommand{\acco}{\overset{\circ\circ}{\mathbf{u}}}
\newcommand{\circdot}[1]{\overset{\circ}{#1}}

\usepackage{xfrac}
\usepackage[bottom]{footmisc}

\usepackage{listings}

\usepackage{xcolor}
\definecolor{codegreen}{rgb}{0,0.6,0}
\definecolor{codegray}{rgb}{0.5,0.5,0.5}
\definecolor{codepurple}{rgb}{0.58,0,0.82}
\definecolor{backcolour}{rgb}{0.95,0.95,0.92}

\lstdefinestyle{mystyle}{
    commentstyle=\color{codegreen},
    keywordstyle=\color{blue},
    numberstyle=\tiny\color{gray},
    stringstyle=\color{purple},
    basicstyle=\ttfamily\footnotesize,
    breakatwhitespace=false,         
    breaklines=true,                 
    captionpos=b,                    
    keepspaces=true,                 
    numbers=left,                    
    numbersep=5pt,                  
    showspaces=false,                
    showstringspaces=false,
    showtabs=false,                  
    tabsize=2,
   lineskip=0.1\baselineskip}

\ifdefined\showcaptionsetup
 \PassOptionsToPackage{caption=false}{subfig}
\fi
\usepackage{subfig}
\makeatother

\usepackage{babel}
\begin{document}
\begin{frontmatter}
\title{A high-order implicit time integration method for linear and nonlinear
dynamics with efficient computation of accelerations}
\author[unsw]{Daniel O'Shea~}
\author[unsw]{Xiaoran Zhang~}
\author[unsw]{Shayan Mohammadian~}
\author[unsw]{Chongmin Song~\corref{cor1}}
\cortext[cor1]{Corresponding author. E-mail address: c.song@unsw.edu.au (C. Song)}
\address[unsw]{Centre for Infrastructure Engineering and Safety, School of Civil
and Environmental Engineering, University of New South Wales, Sydney,
NSW 2052, Australia.}
\begin{abstract}
An algorithm for a family of self-starting high-order implicit time
integration schemes with controllable numerical dissipation is proposed
for both linear and nonlinear transient problems. This work builds
on the previous works of the authors on elastodynamics by presenting
a new algorithm that eliminates the need for factorization of the
mass matrix providing benefit for the solution of nonlinear problems.
The improved algorithm directly obtains the acceleration at the same
order of accuracy of the displacement and velocity using vector operations
(without additional equation solutions). The nonlinearity is handled
by numerical integration within a time step to achieve the desired
order of accuracy. The new algorithm fully retains the desirable features
of the previous works: 1. The order of accuracy is not affected by
the presence of external forces and physical damping; 2. The amount
of numerical dissipation in the algorithm is controlled by a user-specified
parameter $\rho_{\infty}$, leading to schemes ranging from perfectly
nondissipative $A$-stable to $L$-stable; 3. The effective stiffness
matrix is a linear combination of the mass, damping, and stiffness
matrices as in the trapezoidal rule, leading to high efficiency for
large-scale problems. The proposed algorithm, with its elegance and
computational advantages, is shown to replicate the numerical results
demonstrated on linear problems in previous works. Additional numerical
examples of linear and nonlinear vibration and wave propagation are
presented herein. Notably, the proposed algorithms show the same convergence
rates for nonlinear problems as linear problems, and very high accuracy.
It was found that second-order time integration methods commonly used
in commercial software produce significantly polluted acceleration
responses for a common class of wave propagation problems. The high-order
time integration schemes presented here perform noticably better at
suppressing spurious high-frequency oscillations and producing reliable
and useable acceleration responses. The source code written in \texttt{MATLAB}
is available for download at: $\rm{https://github.com/ChongminSong/HighOrderTimeIngt\_PartialFraction.git}$
\end{abstract}
\begin{keyword}
High-order method\sep Implicit time integration\sep Numerical dissipation
\sep Partial fraction\sep Structural dynamics \sep Wave propagation

\end{keyword}
\end{frontmatter}

\section{Introduction}\label{sec:Introduction}

Direct time integration schemes are employed in the numerical simulation
of time-dependent problems in many disciplines of engineering and
science. The requirements for an effective time integration scheme
may vary according to the applications. For some classes of applications,
such as celestial motions~\citep{Fox1984}, the equation of motion
is directly formulated as an initial value problem of ordinary differential
equations. A sufficiently accurate approximation of the exact solution
of the ordinary differential equations under given initial conditions
is required. Several well-known classes of numerical methods, such
as single-step Euler methods, linear multi-step methods and Runge--Kutta
methods using intermediate sub-steps, have been widely studied and
implemented in computer software libraries, see e.g., \citep{Shampine1997,Givoli2023}.

There are many engineering applications where the governing equations
are formulated as partial differential equations in time and space.
One type of such applications addressed in this paper is the problem
of wave propagation in the continuum. Typically, the problem domain
is first discretized spatially by the finite element method \citep{BookBathe2014,Zienkiewvicz:2005}
or other numerical methods. The resulting semi-discrete equation of
motion is expressed as a system of ordinary differential equations
in time. The spatial discretization techniques are accurate for low-frequency
modes, but introduce large dispersion errors for high-frequency modes
\citep{Babuska1995,Kwon2020,Kim2021,Song}. When high-frequency modes
are excited, the solution to the ordinary differential equations will
be polluted by spurious (i.e., non-physical) oscillations. Therefore,
time integration schemes aiming to obtain the exact solution of the
ordinary differential equations, such as the Runge--Kutta methods,
are not used in finite element software packages for structural analysis.
A large number of direct time integration schemes aiming to achieve
high-frequency dissipation while minimizing low-frequency dissipation
have been developed in the area of structural dynamics. Only implicit
schemes are addressed in this paper.

Many implicit time integration methods that have the capacity of introducing
numerical dissipation of high-frequency modes have been developed,
e.g.,~\citep{Houbolt1950,Newmark1959,Wilson1972,Hilber1977,Reusch1988,Chung1993,Bathe2007,Kuo2012,Soares2015,Kim2017a,Noh2018,Kim2020,Behnoudfar2021,Malakiyeh2021,Soares2023,Wang2023}
to name a few. The integration methods that are most commonly found
in commercial finite element packages include the Houbolt method~\citep{Houbolt1950},
the Wilson-$\theta$ method~\citep{Wilson1972}, the Newmark method~\citep{Newmark1959},
the HHT-$\alpha$ method~\citep{Hilber1977}, and the generalized-$\alpha$
method~\citep{Chung1993}. All of these methods are second-order
accurate and share a common characteristic: the effective stiffness
matrix of the implicit equation has the same sparsity and size as
those of the static stiffness matrix, ensuring high computational
efficiency in time-stepping. In practical applications, two parameters
are typically selected: one is the time step size, and the other is
a parameter controlling the amount of numerical dissipation. For the
same numerical dissipation parameter, reducing the time step size
will also reduce the amount of numerical dissipation. Therefore, a
smaller time step size does not necessarily lead to a smaller solution
error in the sense of engineering application~\citep{Malakiyeh2021}.
The optimal choice of the time step size for a given time integration
scheme is conventionally recommended as a Courant-Friedrichs-Lewy
(CFL) number.

Various high-order time integration methods that have an order of
accuracy higher than two have been proposed~ \citep{Reusch1988,Fung1998,Kuo2012,Kim2017a,Kim2017b,Barucq2018,Kim2018,Kim2020,Soares2020,Behnoudfar2020,Behnoudfar2021,Kim2021,Choi2022}.
In comparison with standard second-order time integration schemes,
the potential benefits are twofold: 1) A much larger time step size
can be used to obtain a solution of similar accuracy (although it
is more expensive to advance one time step, the overall performance
can be more efficient), and 2) The numerical dissipation can be more
effectively controlled, i.e., stronger high-frequency dissipation
and less low-frequency dissipation. Numerical methods such as collocation,
differential quadrature, weighted residuals, and matrix exponentials
are applied to derive the time-stepping formulations. Reviews on the
recent progress have been reported by the authors in~\citep{Song2022,Song2024}.
Most of the high-order time integration methods either increase the
size or destroy the sparsity of the implicit equation, and, as a result,
the computational cost increases rapidly with the problem size. These
methods are then not capable of handling problems of as large a size
as the conventional second-order methods can.

Composite time integration approaches, which are constructed by splitting
a single time step into two or more sub-steps, have been adopted in
many recent works \citep{Bathe2007,Bathe2005,Bathe2012,Noh2013,Kim2017,Wen2017,Kim2018,Noh2018,Noh2019,Noh2019b,Malakiyeh2019,Kim2020b,Malakiyeh2021,Kwon2021,Choi2022,Li2022,Noh2023}.
By judiciously choosing different time integration methods or different
values of the parameters for the individual sub-steps, improved accuracy
and numerical properties are achieved. A comprehensive review and
an attempt to unify the different formulations have been reported
by Kim and Reddy in~\citep{Kim2020a}. The Bathe method has been
extended to third-order and fourth-order accuracy in~\citep{Kwon2021,Choi2022}.
A high-order composite time integration method is proposed in~\citep{Li2022},
where the accuracy order equals the number of sub-steps.

Self-starting with the initial displacement and velocity is a desirable
property of a time integration method~ \citep{Hilber1978}. Self-starting
implicit methods have been reported in \citep{Kim2020b,Kim2021a,Li2022,Wang2023},
to name a few. As the initial acceleration is not required, the factorization
of the mass matrix becomes unnecessary. This reduces the computational
time and memory costs when a consistent mass matrix is used~\citep{Kim2020b}.
In situations where acceleration responses are needed, for example,
in earthquake engineering and in structural health monitoring employing
accelerometers, the acceleration can be obtained from the interpolation
built in the schemes when the accuracy order is not higher than three
\citep{Li2022,Wang2023}. To maintain the same order of accuracy with
the displacement and velocity in a higher-order scheme, the acceleration
needs to be solved from the equation of motion, thus requiring factorization
of the mass matrix and compromising the computational advantages
of high-order self-starting methods.

A family of computationally effective high-order time integration
schemes has been developed in~\citep{Song,Song2024} by designing
rational approximations of the matrix exponential in the analytical
solution of the equation of motion. A salient feature of these schemes
is that the effective stiffness matrix is a linear combination of
the mass, damping, and stiffness matrices, similar to the conventional
second-order methods such as the trapezoidal rule and Newmark method.
These schemes conform to the concept of the composite time integration
method, where the solution at the intermediate time point is only
used to calculate the solution at the full step~\citep{Bathe2005,Kwon2021,Kim2018}.
Two types of rational approximations are adopted. In~\citep{Song},
the numerical dissipation is controlled by mixing the Padé expansions
of orders $(M,M)$ and $(M-1,M)$. When the rational function is a
Padé expansion of order $(M,M)$, i.e., the diagonal Padé expansion,
the scheme is strictly non-dissipative, $2M$-th order accurate and
\emph{A-}stable. Otherwise, the scheme is dissipative and $(2M-1)$-th
order accurate. At the other extreme, the Padé expansion of order
$(M-1,M)$ leads to an $L$-stable scheme. Each solution of an implicit
equation increases the order of accuracy by one when the effective
stiffness matrix is real and by two when the effective stiffness matrix
is complex. In~~\citep{Song2024}, the denominator of the rational
approximation is selected to have one real root, i.e., a single multiple
root. Unconditionally stable schemes up to the sixth-order accuracy
are identified. Each solution of an implicit equation increases the
order of accuracy by one. In this approach, the effective stiffness
matrices of all the implicit equations are the same, reducing the
computational cost of matrix factorization when a direct solver is
used. Although the algorithm in~\citep{Song,Song2024} is self-starting
without using the acceleration for time-stepping, the algorithm requires
the factorization of the mass matrix to compute $\mathbf{M}^{-1}\mathbf{f}$,
where $\mathbf{M}$ is the mass matrix and $\mathbf{f}$ is the external
force vector, and therefore the full benefits of a self-starting algorithm
are not realized. In a nonlinear analysis, this operation incurs additional
computational cost since the force vector $\mathbf{f}$ within a time
step is unknown beforehand. .

This paper introduces a new time-stepping algorithm for the family
of high-order schemes in~\citep{Song,Song2024}, and demonstrates
their performance on a range of linear and nonlinear transient problems.
Using partial fractions of the rational approximations, this algorithm
eliminates the need to compute $\mathbf{M}^{-1}\mathbf{f}$, and instead
it is sufficient to sample the force vector $\mathbf{f}$ at specified
time instances, leading to a family of efficient high-order schemes
for nonlinear dynamics. Furthermore, a new formulation to compute
the acceleration is developed, which involves only vector operations
and no additional equation solutions. The order of accuracy remains
the same as that of the displacement and velocity, i.e., equal to
the order of the time integration scheme. It is worthwhile to reemphasize
that acceleration is not needed for the time-stepping and the algorithm
itself is self-starting. The factorization of the mass matrix becomes
unnecessary even if high-order accuracy of acceleration is required.
The full benefits of a self-starting algorithm are realized for both
nonlinear and linear cases. In numerical examples, the proposed algorithms
show the same convergence rates for nonlinear problems as linear problems,
and very high accuracy.

The remainder of the paper is organized as follows: Section 2 presents
the time-stepping solution scheme using matrix exponentials, including
analytic treatment of the non-homogeneous term using a series expansion
of the excitation force vector and a recursive scheme to evaluate
the integrals. Section 3 presents the rational approximation of the
matrix exponential using Padé expansions. Section 4 presents the algorithmic
novelty of the work, deriving time-stepping algorithms that eliminate
the need to compute $\mathbf{M}^{-1}\mathbf{f}$ terms through expansions
using partial fractions for cases where the divisor polynomial possesses
either 1) distinct roots, or 2) a single multiple root. Solutions
to the time-stepping equation for both cases have an identical representation
which is presented in Section 5. In Section 6, the algorithm to obtain
acceleration to the same order of accuracy as displacement without
additional equation solution and using only vector operations is developed.
Discussion on treatment of the non-homogeneous term is offered in
Section 7, with sample algorithms for the time-stepping procedures
presented in Section 8. Section 9 provides a range of numerical examples
which demonstrate the efficacy of the proposed approach on accurately
modelling both linear and nonlinear dynamic problems, followed by
conclusions drawn. 

Additional examples of using the algorithms herein are presented in
the Appendix, alongside sample MATLAB code for many of the important
functions used.

\section{Time-stepping scheme using rational approximation of matrix exponential}

For a structural dynamics problem, the governing equations of motion
are expressed as a system of second-order ODEs as follows:
\begin{equation}
\mathbf{M\ddot{u}}(t)+\mathbf{f}_{I}(\mathbf{u}(t),\mathbf{\dot{u}}(t))=\mathbf{f}_{E}(t)\,,\thinspace\thinspace\thinspace\mathbf{u}_{0}=\mathbf{u}(0),\thinspace\dot{\mathbf{u}}_{0}=\dot{\mathbf{u}}(0)\thinspace,\label{eq:EqofMotionNonlinear}
\end{equation}
where $\mathbf{M}$ denotes the mass matrix, $\mathbf{f}_{I}$ gives
the internal force vector, and $\mathbf{f}_{E}$ represents the external
excitation force. The vectors $\mathbf{u}$, $\dot{\mathbf{u}}$ and
$\mathbf{\ddot{u}}$ denote the displacement, velocity, and acceleration,
respectively, each a function of time $t$, and $\mathbf{u}_{0},\dot{\mathbf{u}}_{0}$
provide the initial conditions. Using a series expansion about time
$t=t_{n-1}$, the linear part of the internal force vector can be
separated from the nonlinear remainder such that
\begin{equation}
\mathbf{M}\mathbf{\ddot{u}}(t)+\mathbf{C}_{n-1}\mathbf{\dot{u}}(t)+\mathbf{K}_{n-1}\mathbf{u}(t)=\mathbf{f}(\mathbf{u}(t),\dot{\mathbf{u}}(t))\,,\label{eq:linear=000020eq=000020of=000020motion}
\end{equation}
with
\begin{equation}
\mathbf{C}(t)=\frac{\partial\mathbf{f}_{I}}{\partial\dot{\mathbf{u}}}\,,\thinspace\thinspace\thinspace\mathbf{K}(t)=\frac{\partial\mathbf{f}_{I}}{\partial\mathbf{u}}\thinspace,\label{eq:damping-stiffness}
\end{equation}
denoting tangent damping and stiffness matrices at time $t$, respectively,
and the right-hand side of Eq.~(\ref{eq:linear=000020eq=000020of=000020motion})
being described by the nonlinear force vector
\begin{equation}
\mathbf{f}(\mathbf{u}(t),\dot{\mathbf{u}}(t))=\mathbf{f}_{E}(t)-\mathbf{f}_{I}(\mathbf{u}(t),\mathbf{\dot{u}}(t))+\mathbf{C}_{n-1}\mathbf{\dot{u}}(t)+\mathbf{K}_{n-1}\mathbf{u}(t)\ .\label{eq:nonlinear=000020force}
\end{equation}

The above equations of motion (Eq.~(\ref{eq:linear=000020eq=000020of=000020motion}))
are now considered in the time step $n$, denoting the time interval
$t_{n-1}\leq t\leq t_{n}$ (for $n=1,2,3,\ldots$). The time step
size is thus defined such that $\Delta t=t_{n}-t_{n-1}$. By introducing
a dimensionless time variable $s,$ any point in time within a given
time step $n$ can be determined by 
\begin{equation}
t(s)=t_{n-1}+s\Delta t,\quad\;0\leq s\leq1\,,\label{eq:dimensionless=000020time}
\end{equation}
with $t(0)=t_{n-1}$ and $t(1)=t_{n}$ returning the beginning and
end of the time step, respectively. In this work, the derivative of
a quantity with respect to the dimensionless time $s$ is denoted
through the use of a circle ($\circ$) above the symbol. The incremental
velocity and acceleration vectors within a given time step $t_{n-1}\leq t\leq t_{n}$
are therefore represented as 
\begin{equation}
\mathbf{\dot{u}}=\frac{1}{\Delta t}\frac{\mathrm{d\mathbf{u}}}{\mathrm{d}s}=\frac{1}{\Delta t}\overset{\circ}{\mathbf{u}}\,,\label{eq:velocity=000020in=000020dimensionless=000020time}
\end{equation}
\begin{equation}
\mathbf{\ddot{u}}=\frac{1}{\Delta t^{2}}\frac{\mathrm{d^{2}\mathbf{u}}}{\mathrm{d}s^{2}}=\frac{1}{\Delta t^{2}}\overset{\circ\circ}{\mathbf{u}}\,.\label{eq:acceleration=000020in=000020dimensinoless=000020time}
\end{equation}

Therefore, the equations of motion (Eq.~(\ref{eq:linear=000020eq=000020of=000020motion}))
can now be expressed with respect to the dimensionless time as:
\begin{equation}
\mathbf{M}\overset{\circ\circ}{\mathbf{u}}(t)+\Delta t\,\mathbf{C}_{n-1}\overset{\circ}{\mathbf{u}}(t)+\Delta t^{2}\,\mathbf{K}_{n-1}\mathbf{u}(t)=\Delta t^{2}\,\mathbf{f}(\mathbf{z}(t))\,.\label{eq:eqn=000020of=000020motion=000020in=000020dimensionless}
\end{equation}

Equation~(\ref{eq:eqn=000020of=000020motion=000020in=000020dimensionless})
can be transformed into a system of first-order ODEs by introducing
a state-space vector $\mathbf{z}$ defined as
\begin{equation}
\mathbf{z}(t)=\begin{Bmatrix}\overset{\circ}{\mathbf{u}}(t)\\
\mathbf{u}(t)
\end{Bmatrix}\,,\label{eq:state=000020space=000020vector}
\end{equation}
and subsequently taking the derivative of $\mathbf{z}$ with respect
to $s$:
\begin{equation}
\overset{\circ}{\mathbf{z}}(t)\equiv\frac{\mathrm{d\mathbf{z}}}{\mathrm{d}s}=\mathbf{A}\mathbf{z}(t)+\overline{\mathbf{f}}(\mathbf{z}(t))\,.\label{eq:1st=000020order=000020ODE}
\end{equation}
Here, $\mathbf{A}$ is the tangent coefficient matrix at time $t_{n-1}$
defined as
\begin{equation}
\mathbf{\mathbf{A}}=\begin{bmatrix}-\Delta t\mathbf{M}^{-1}\mathbf{C}_{n-1} & -\Delta t^{2}\mathbf{M}^{-1}\mathbf{K}_{n-1}\\
\mathbf{I} & \mathbf{0}
\end{bmatrix}\,,\label{eq:matrix=000020A}
\end{equation}
and $\overline{\mathbf{f}}$ is the nonlinear force vector (at time
$t$) given by:
\begin{equation}
\overline{\mathbf{f}}(\mathbf{z}(t))=\begin{Bmatrix}\Delta t^{2}\ \mathbf{M}^{-1}\mathbf{f}(\mathbf{z}(t))\\
\mathbf{0}
\end{Bmatrix}.
\end{equation}
The analytical solution to the system of first-order ODEs given by
Eq.~(\ref{eq:1st=000020order=000020ODE}) is expressed using the
matrix exponential function as 
\begin{equation}
\mathbf{z}(t_{n-1}+s\Delta t)=e^{\mathbf{A}s}\mathbf{z}_{n-1}+e^{\mathbf{A}s}\intop_{0}^{s}e^{-\mathbf{A}\tau}\overline{\mathbf{f}}(\mathbf{z}(t_{n-1}+\tau\thinspace\Delta t))\mathrm{d}\tau\,.\label{eq:general=000020solution}
\end{equation}

\subsection{Analytical approximation of non-homogeneous term}\label{subsec:Analytical-approximation-of}

Following \citep{Song2022}, the integration in Eq.~(\ref{eq:general=000020solution})
is performed analytically by approximating the vector $\mathbf{f}$
with Taylor expansion of order $p_{\mathrm{f}}$ at the mid-point
of the time step $n$
\begin{equation}
\mathbf{f}(s)=\sum\limits_{k=0}^{p_{\mathrm{f}}}\tilde{\mathbf{f}}^{(k)}(s-0.5)^{k}=\tilde{\mathbf{f}}^{(0)}+\tilde{\mathbf{f}}^{(1)}(s-0.5)+\tilde{\mathbf{f}}^{(2)}(s-0.5)^{2}+\ldots+\tilde{\mathbf{f}}^{(p_{\mathrm{f}})}(s-0.5)^{p_{\mathrm{f}}}\,,\label{eq:ForceExpansion}
\end{equation}
where $\tilde{\mathbf{f}}^{(k)}$ ($k=0,\,1,\,2,\,\ldots,\,p_{\mathrm{f}}$)
denotes the $k$-th order derivatives of the vector $\mathbf{f}$.
This series expansion is compactly written as the matrix-vector product
\begin{equation}
\mathbf{f}(s)=\tilde{\mathbf{F}}\cdot\mathbf{s}(s)\,,\label{eq:ForceFitting1}
\end{equation}
with the columns of matrix $\tilde{\mathbf{F}}$ storing force-derivatives
vectors and the vector $\mathbf{s}$ storing powers of $(s-0.5)$
according to
\begin{equation}
\tilde{\mathbf{F}}=\left[\begin{array}{ccccc}
\tilde{\mathbf{f}}^{(0)} & \tilde{\mathbf{f}}^{(1)} & \tilde{\mathbf{f}}^{(2)} & \ldots & \tilde{\mathbf{f}}^{(p_{\mathrm{f}})}\end{array}\right]\,,\label{eq:ForceFitting2}
\end{equation}
and
\begin{equation}
\mathbf{s}(s)=\left[\begin{array}{ccccc}
1 & (s-0.5) & (s-0.5)^{2} & \ldots & (s-0.5)^{p_{\mathrm{f}}}\end{array}\right]^{\mathrm{T}}\,.\label{eq:ForceFitting3}
\end{equation}
The evaluation of matrix $\tilde{\mathbf{F}}$ from the force $\mathbf{f}$
at specified time instances is discussed in Section~\ref{sec:Calculation-of-force-deriv}.

Substituting Eq.~(\ref{eq:ForceExpansion}) into Eq.~(\ref{eq:general=000020solution}),
the solution at the end of time step $t_{n}$ (i.e., $s=1$) is obtained
as
\begin{equation}
\mathbf{z}_{n}=e^{\mathbf{A}}\mathbf{z}_{n-1}+\sum\limits_{k=0}^{p_{\mathrm{f}}}\mathbf{B}_{k}\left\{ \begin{array}{c}
\Delta t^{2}\thinspace\mathbf{M^{-1}}\tilde{\mathbf{f}}^{(k)}\\
\mathbf{0}
\end{array}\right\} \,,\label{eq:sol=000020s=00003D00003D1}
\end{equation}
where the matrices $\mathbf{B}_{k}$ are obtained through integration
by parts and can be determined recursively:
\begin{equation}
\begin{split}\mathbf{B}_{k} & =\mathbf{A}^{-1}\left(k\mathbf{B}_{k-1}+\left(-\tfrac{1}{2}\right)^{k}\left(e^{\mathbf{A}}-(-1)^{k}\mathbf{I}\right)\right)\,,\quad\;\forall k=0,1,2,\ldots,p_{\mathrm{f}}\,,\end{split}
\label{eq:intergrationExpmI}
\end{equation}
with the starting value at $k=0$ given by
\begin{equation}
\mathbf{B}_{0}=\mathbf{A}^{-1}\left(e^{\mathbf{A}}-\mathbf{I}\right)\,.\label{eq:integrationExpm0}
\end{equation}

\section{Rational approximation of matrix exponential}\label{sec:Rational-approximation-of}

To derive an efficient time-stepping scheme, Padé expansions are often
employed (see \citep{Song2022,Song}) to approximate the matrix exponential
$e^{\mathbf{A}}$. Here, the roots of the denominator polynomial are
all distinct and can be either real or complex. In \citep{Song2024},
a rational approximation is constructed in such a way that the denominator
possesses only a single real root. The proposed time-stepping scheme
based on partial fractions is applicable to both cases of real or
complex roots. For generality, the following rational approximation
$\mathbf{R}=R(\mathbf{A})$ of the matrix exponential $e^{\mathbf{A}}$
is considered
\begin{equation}
e^{\mathbf{A}}\approx\mathbf{R}=\mathbf{\dfrac{\mathbf{P}}{Q}}=\dfrac{p_{0}\mathbf{I}+p_{1}\mathbf{A}+\ldots+p_{M}\mathbf{A}^{M}}{q_{0}\mathbf{I}+q_{1}\mathbf{A}+\cdots+q_{M}\mathbf{A}^{M}}\,,\label{eq:ExpmRational}
\end{equation}
where $\mathbf{P}=P(\mathbf{A})$ and $\mathbf{Q}=Q(\mathbf{A})$
are polynomials of matrix $\mathbf{A}$ of order $M$. Note that the
matrix product $\mathbf{Q}^{-1}\mathbf{P}$ is commutative (i.e.,
$\mathbf{P}\mathbf{Q}^{-1}=\mathbf{Q}^{-1}\mathbf{\mathbf{P}}$),
and the scalar coefficients $p_{i}$ and $q_{i}$ are real numbers.
To ensure that $e^{\mathbf{0}}=\mathbf{I}$ holds, it is necessary
for $p_{0}=q_{0}$. For a stable time-stepping scheme, it is necessary
to have $q_{M}\neq0$ and $|p_{M}/q_{M}|=\rho_{\infty}\leq1$, where
$\rho_{\infty}$ denotes the spectral radius at the high-frequency
limit. When $p_{M}=0$, it follows that $\rho_{\infty}=0$ and the
scheme is `$L$-stable'.

Replacing the matrix exponential $e^{\mathbf{A}}$ in Eqs.~(\ref{eq:sol=000020s=00003D00003D1}),
(\ref{eq:intergrationExpmI}) and (\ref{eq:integrationExpm0}) with
the rational approximation in Eq.~(\ref{eq:ExpmRational}), the time-stepping
equation becomes
\begin{equation}
\mathbf{z}_{n}=\frac{\mathbf{P}}{\mathbf{Q}}\mathbf{z}_{n-1}+\frac{1}{\mathbf{Q}}\sum\limits_{k=0}^{p_{\mathrm{f}}}\mathbf{C}_{k}\left\{ \begin{array}{c}
\Delta t^{2}\mathbf{M^{-1}}\tilde{\mathbf{f}}^{(k)}\\
\mathbf{0}
\end{array}\right\} \,,\label{eq:Stepping}
\end{equation}
where the matrices $\mathbf{C}_{k}$ (not to be confused with the
tangent damping matrix in Eq.~(\ref{eq:damping-stiffness})) are
obtained from Eqs.~(\ref{eq:intergrationExpmI}) and (\ref{eq:ExpmRational})
through the following recursive relation 
\begin{equation}
\mathbf{C}_{k}=\mathbf{Q}\mathbf{B}_{k}=\mathbf{A}^{-1}\left(k\mathbf{C}_{k-1}+\left(-0.5\right)^{k}(\mathbf{P}-(-1)^{k}\mathbf{Q})\right);\qquad(k=1,2,\ldots,p_{\mathrm{f}})\,,\label{eq:Stepping_Ck}
\end{equation}
which are evaluated starting from 
\begin{equation}
\mathbf{C}_{0}=\mathbf{Q}\mathbf{B}_{0}=\mathbf{A}^{-1}\left(\mathbf{P}-\mathbf{Q}\right)\,,\label{eq:Stepping_C0}
\end{equation}
using $\mathbf{B}_{0}$ in Eq.~(\ref{eq:integrationExpm0}). Since
both $\mathbf{P}$ and $\mathbf{Q}$ are matrix polynomials of order
$M$ , the matrices $\mathbf{C}_{k}$ ($k=0,1,\ldots,p_{\mathrm{f}}$)
are matrix polynomials with degrees at most $M-1$. They can hence
be expanded by the power series 
\begin{equation}
\mathbf{C}_{k}=C_{k}(\mathbf{A})=c_{k0}\mathbf{I}+c_{k1}\mathbf{A}+\ldots+c_{k(M-1)}\mathbf{A}^{M-1}=\sum_{i=0}^{M-1}c_{ki}\mathbf{A}^{i};\qquad(k=0,1,\ldots p_{\mathrm{f}})\,,\label{eq:Stepping_MatrixPolyCk}
\end{equation}
where $c_{ki}$ are scalar coefficients; stored for later use in a
coefficient matrix $\mathbf{c}$:
\begin{equation}
\mathbf{c}=[c_{ki}];\qquad(k=0,1,\ldots p_{\mathrm{f}};\quad i=0,1,\ldots,M-1)\,.\label{eq:ScalarCoeff_ckj}
\end{equation}

\section{Partial fraction expansion of rational approximations}\label{sec:Partial-fraction-expansion}

The time-stepping solution in Eq.~(\ref{eq:Stepping}) has been used
to produce self-starting high-order time integration schemes with
controllable numerical dissipation in \citep{Song2022,Song,Song2024},
and performance tested against a range of numerical examples. In those
works, the algorithms used to execute the time-stepping schemes required
calculation of the $\mathbf{M}^{-1}\tilde{\mathbf{f}}^{(k)}$ terms
in Eq.~(\ref{eq:Stepping}), plus additional equations solutions
to determine the acceleration. Therefore, whilst the algorithms were
self-starting, the full benefit of a self-starting algorithm (that
is, avoiding factorization of the mass matrix, or extraneous equation
solutions) were not realized. The computational benefit of avoiding
calculation of $\mathbf{M}^{-1}\mathbf{f}$ terms in time-stepping
algorithms is most prominent in nonlinear problems where $\mathbf{f}$
is unknown within the current time step, and iterative solution procedures
require additional calculations to be made.

In this section, a new approach to deriving an efficient time-stepping
scheme is presented by applying partial fraction expansion of the
rational functions $\mathbf{P}/\mathbf{Q}$ and $\mathbf{C}_{k}/\mathbf{Q}$
($k=1,\,2,\,\ldots,\,p_{\mathrm{f}}$) that are present in Eq.~(\ref{eq:Stepping}).
The general case of the divisor polynomial $\mathbf{Q}$ having distinct
roots is discussed in Section~\ref{subsec:Distinct-roots}, followed
by the specialized case of $\mathbf{Q}$ possessing a single multiple
real root presented in Section~\ref{subsec:Single-multiple-roots}.
The algorithmic implementation of both is then discussed. The algorithms
avoid any calculation of terms of the form $\mathbf{M}^{-1}\cdot\mathbf{f}$;
and also allow acceleration to be calculated at the same order of
accuracy as displacment and velocity using only vector operations.
Sample MATLAB code with illustrative examples are provided in \ref{subsec:Appendix-Single-multiple-root}.

\subsection{Case 1: Distinct roots }\label{subsec:Distinct-roots}

\subsubsection{Evaluation of matrix exponential}

Consider that the polynomial $Q(\mathbf{A})$ possesses distinct roots,
so the rational approximation Eq.~(\ref{eq:ExpmRational}) is evaluated
according to the Padé expansion as outlined in \citep{Song}. Both
the numerator $\mathbf{P}$ and denominator $\mathbf{Q}$ in Eq.~(\ref{eq:ExpmRational})
are selected as polynomials of order $M$. To begin the derivation,
the rational expansion in Eq.~(\ref{eq:ExpmRational}) is first decomposed
as 
\begin{equation}
\mathbf{\dfrac{\mathbf{P}}{Q}}=\dfrac{p_{M}}{q_{M}}\mathbf{I}+\dfrac{p_{l0}\mathbf{I}+p_{l1}\mathbf{A}+\ldots+p_{l(M-1)}\mathbf{A}^{M-1}}{q_{0}\mathbf{I}+q_{1}\mathbf{A}+\cdots+q_{M}\mathbf{A}^{M}}=\dfrac{p_{M}}{q_{M}}\mathbf{I}+\dfrac{\mathbf{P}_{\mathrm{L}}}{\mathbf{Q}}\,,\label{eq:Pade}
\end{equation}
where ${p_{M}}/{q_{M}}$ is the ratio of coefficients of the highest-order
term. The degree of the resulting numerator polynomial $\mathbf{P}_{\mathrm{L}}$
is equal to $M-1$, one degree lesser than that of the denominator
$\mathbf{Q}$. The order of accuracy of the Padé expansion of the
matrix exponential is then equal to $2M-1$ when numerical dissipation
is introduced, and $2M$ in the absence of numerical dissipation.
The coefficients of $\mathbf{P}_{\mathrm{L}}$ are simply 
\begin{equation}
p_{\mathrm{L}i}=p_{i}-q_{i}\dfrac{p_{M}}{q_{M}},\qquad(i=0,1,\ldots,M-1)\,.\label{eq:Coefficient-pli}
\end{equation}

The polynomial $\mathbf{Q}$ of degree $M$ is then factorized using
its roots $r_{i}$ \citep{Song} such that 
\begin{equation}
\mathbf{Q}=\prod_{i=1}^{M}\left(r_{i}\mathbf{I}-\mathbf{A}\right)=\left(r_{1}\mathbf{I}-\mathbf{A}\right)\left(r_{2}\mathbf{I}-\mathbf{A}\right)\ldots\left(r_{M}\mathbf{I}-\mathbf{A}\right)\,,\label{eq:Qfactor}
\end{equation}
with the roots being either real or pairs of complex conjugates \citep{Song}.
Since the roots are distinct, all the irreducible polynomials are
linear. In the partial fraction decomposition procedure below, it
is sufficient to consider the power $\mathbf{A}^{j}$ ($j<M$) as
the numerator since any polynomial can then be constructed by its
linear combination. The partial fraction decomposition is expressed
as 
\begin{equation}
\dfrac{\mathbf{A}^{j}}{\mathbf{Q}}=\dfrac{\mathbf{A}^{j}}{\left(r_{1}\mathbf{I}-\mathbf{A}\right)\left(r_{2}\mathbf{I}-\mathbf{A}\right)\ldots\left(r_{M}\mathbf{I}-\mathbf{A}\right)}=\dfrac{b_{1}}{r_{1}\mathbf{I}-\mathbf{A}}+\dfrac{b_{2}}{r_{2}\mathbf{I}-\mathbf{A}}+\ldots+\dfrac{b_{M}}{r_{M}\mathbf{I}-\mathbf{A}}\,,\label{eq:PF}
\end{equation}
where $b_{i}$ are constants to be determined. This can be written
compactly as: 
\begin{equation}
\sum_{i=1}^{M}\dfrac{b_{i}}{r_{i}\mathbf{I}-\mathbf{A}}=\dfrac{\mathbf{A}^{j}}{\prod_{i=1}^{M}\left(r_{i}\mathbf{I}-\mathbf{A}\right)}\,,\label{eq:PF_sum}
\end{equation}
and, multiplying both sides of this expression with $\prod_{i=1}^{M}\left(r_{i}\mathbf{I}-\mathbf{A}\right)$
leads to 
\begin{equation}
\sum_{i=1}^{M}\left(b_{i}\prod_{j=1,j\neq i}^{M}(r_{j}\mathbf{I}-\mathbf{A})\right)=\mathbf{A}^{j}\,,\label{eq:PF-expansionCoef1}
\end{equation}
noting that the denominator $(r_{i}\mathbf{I}-\mathbf{A})$ on the
left-hand side of Eq.~(\ref{eq:PF_sum}) cancels with the equivalent
term within the product. To determine the coefficients $b_{i}$, the
particular case of $\mathbf{A}=r_{i}\mathbf{I}$ can be considered.
Thus, all terms on the left-hand side of Eq.~(\ref{eq:PF-expansionCoef1})
which contain $(r_{i}\mathbf{I}-\mathbf{A})$ vanish, except for the
term associated with $b_{i}$, leading to: 
\begin{equation}
b_{i}\prod_{j=1,j\neq i}^{M}(r_{j}-r_{i})=r_{i}^{j}\,.\label{eq:PF-expansionCoef2}
\end{equation}
For conciseness, the scalar $a_{i}$ is introduced such that:
\begin{equation}
b_{i}=a_{i}r_{i}^{j}\qquad\mathrm{where}\quad a_{i}=\dfrac{1}{\prod_{j=1,j\neq i}^{M}(r_{j}-r_{i})}\,.\label{eq:PF-expansionSln}
\end{equation}

Using this definition for $b_{i}$ in the partial fraction decomposition
of the matrix polynomials, Eq.~(\ref{eq:PF_sum}) is now rewritten
as 
\begin{equation}
\dfrac{\mathbf{A}^{j}}{\mathbf{Q}}=\sum_{i=1}^{M}\dfrac{a_{i}r_{i}^{j}}{r_{i}\mathbf{I}-\mathbf{A}}\,.\label{eq:PF-noci}
\end{equation}

The partial fraction decomposition of the Padé expansion in Eq.~(\ref{eq:Pade})
can therefore be written as 
\begin{equation}
\dfrac{\mathbf{P}}{\mathbf{Q}}=\dfrac{p_{M}}{q_{M}}\mathbf{I}+\sum_{i=1}^{M}\dfrac{a_{i}}{r_{i}\mathbf{I}-\mathbf{A}}P_{\mathrm{L}}(r_{i})\,,\label{eq:PF_P/Q}
\end{equation}
with the following polynomials of each distinct root $r_{i}$ defined
\begin{equation}
P_{\mathrm{L}}(r_{i})=p_{\mathrm{L}0}+p_{\mathrm{L}1}r_{i}+\ldots+p_{\mathrm{L}(M-1)}r_{i}^{M-1}\,.\label{eq:PF_Pl(ri)}
\end{equation}

The coefficients of this polynomial are found using Eq.~(\ref{eq:Coefficient-pli}).
In a similar way to Eq.~(\ref{eq:PF_P/Q}), the partial fraction
decomposition of the matrix polynomials $\mathbf{C}_{k}$ in Eq.~(\ref{eq:Stepping_MatrixPolyCk})
is expressed as the summation
\begin{equation}
\dfrac{\mathbf{C}_{k}}{\mathbf{Q}}=\sum_{i=1}^{M}\dfrac{a_{i}}{r_{i}\mathbf{I}-\mathbf{A}}C_{k}(r_{i})\,,\label{eq:PF_Ck/Q}
\end{equation}
with polynomials
\begin{equation}
C_{k}(r_{i})=c_{k0}+c_{k1}r_{i}+\ldots+c_{k(M-1)}r_{i}^{M-1}=\sum_{j=0}^{M-1}c_{kj}r_{i}^{j}\,.\label{eq:PFCoeff_Ck(ri)}
\end{equation}

Therefore, the non-homogeneous term in the time-stepping equation
(Eq.~(\ref{eq:Stepping})) is able to be formulated as 
\begin{equation}
\dfrac{1}{\mathbf{Q}}\sum\limits_{k=0}^{p_{\mathrm{f}}}\mathbf{C}_{k}\left\{ \begin{array}{c}
\Delta t^{2}\mathbf{M^{-1}}\tilde{\mathbf{f}}^{(k)}\\
\mathbf{0}
\end{array}\right\} =\sum\limits_{k=0}^{p_{\mathrm{f}}}\sum_{i=1}^{M}\dfrac{a_{i}}{r_{i}\mathbf{I}-\mathbf{A}}C_{k}(r_{i})\left\{ \begin{array}{c}
\Delta t^{2}\mathbf{M^{-1}}\tilde{\mathbf{f}}^{(k)}\\
\mathbf{0}
\end{array}\right\} \,.\label{eq:PF-nonhomo-2sum}
\end{equation}

Given that the force vectors $\tilde{\mathbf{f}}^{(k)}$ only depend
on the polynomial order $p_{\mathrm{f}}$ and are independent of the
matrix $\mathbf{A}_{\mathrm{r}}$ (as indicated in Eq.~(\ref{eq:ForceExpansion})),
the double-summation is rearranged to consider the summation of the
force vectors first. Thus, for each sub-step $i$ the force term is
summed over prior to the matrix-vector multiplication, and the following
expression holds 
\begin{equation}
\dfrac{1}{\mathbf{Q}}\sum\limits_{k=0}^{p_{\mathrm{f}}}\mathbf{C}_{k}\left\{ \begin{array}{c}
\Delta t^{2}\mathbf{M^{-1}}\tilde{\mathbf{f}}^{(k)}\\
\mathbf{0}
\end{array}\right\} =\sum_{i=1}^{M}\dfrac{a_{i}}{r_{i}\mathbf{I}-\mathbf{A}}\left\{ \begin{array}{c}
\Delta t^{2}\mathbf{M^{-1}}\mathbf{f}_{\mathrm{r}i}\\
\mathbf{0}
\end{array}\right\} \,,\label{eq:PF-nonhomo-suminside}
\end{equation}
where, for simplicity in notation, vectors $\mathbf{f}_{\mathrm{r}i}$
are introduced 
\begin{equation}
\mathbf{f}_{\mathrm{r}i}=\sum\limits_{k=0}^{p_{\mathrm{f}}}\tilde{\mathbf{f}}^{(k)}C_{k}(r_{i})\,,\label{eq:PF-nonhomo-fri}
\end{equation}
which, using the matrix $\tilde{\mathbf{F}}$ and vector $\mathbf{c}$
previously defined in Eq.~(\ref{eq:ForceFitting2}) and Eq.~(\ref{eq:ScalarCoeff_ckj}),
respectively, can be condensed into matrix-vector form as 
\begin{equation}
\mathbf{f}_{\mathrm{r}i}=\tilde{\mathbf{F}}\cdot\mathbf{c}\mathbf{r}_{i}\,,\label{eq:PF-nonhomo-fri-mtx}
\end{equation}
using a vector $\mathbf{r}_{i}$ defined for each distinct root as
\begin{equation}
\mathbf{r}_{i}=\left[\begin{array}{ccccc}
1 & r_{i}^{1} & r_{i}^{2} & \ldots & r_{i}^{M-1}\end{array}\right]\,.\label{eq:PF-ri}
\end{equation}

Substituting the partial fraction expansions of Eq.~(\ref{eq:PF_P/Q})
and Eq.~(\ref{eq:PF-nonhomo-suminside}) into Eq.~(\ref{eq:Stepping}),
the time-stepping equation can now be reformulated as 
\begin{equation}
\mathbf{z}_{n}=\rho\thinspace\mathbf{z}_{n-1}+\sum_{i=1}^{M}a_{i}\mathbf{y}_{i}\,,\label{eq:PF-sol}
\end{equation}
where $\rho={p_{M}}/{q_{M}}=(-1)^{M}\rho_{\infty}$ is related to
the spectral radius at the high-frequency limit, and the following
auxiliary variables are introduced (related to each root $r_{i}$):
\begin{equation}
\mathbf{y}_{i}=\dfrac{1}{r_{i}\mathbf{I}-\mathbf{A}}\left(P_{\mathrm{L}}(r_{i})\mathbf{z}_{n-1}+\left\{ \begin{array}{c}
\Delta t^{2}\mathbf{M^{-1}}\mathbf{f}_{\mathrm{r}i}\\
\mathbf{0}
\end{array}\right\} \right)\,.\label{eq:DistinctRoot-yi}
\end{equation}

Note that for nonlinear problems, the vectors $\mathbf{\mathbf{f}}_{\mathrm{r}i}$
are functions of the unknown target state-space $\mathbf{z}_{n}$
and therefore an iterative solution procedure is required. The time-stepping
solution Eq.~(\ref{eq:PF-sol}) is therefore solved using standard
fixed-point iteration techniques. The evaluation of the state-space
$\mathbf{z}_{n}$ is presented in Section~\ref{subsec:Interpolation-of-displacements}.

To execute the time-stepping equation Eq.~(\ref{eq:PF-sol}), for
a chosen polynomial order $M$ and user-defined $\rho_{\infty}$,
the following quantities need to be determined prior to commencing
the time-stepping procedure:
\begin{enumerate}
\item the roots $r_{i}$ of $Q(x)$ 
\item the term $p_{M}/q_{M}$ and the coefficients and $p_{\mathrm{L}i}$
($i=0,\,1,\,2,\,\ldots,\,M-1$) 
\item the coefficients of the partial fractions $a_{i}$ ($i=0,\,1,\,2,\,\ldots,\,M$),
see Eq.~(\ref{eq:PF-expansionSln}) 
\item the coefficient matrix for external excitation vectors $\mathbf{f}_{\mathrm{r}i}$,
given by: $\mathbf{c}=[c_{ki}];\qquad(k=0,1,\ldots p_{f};\quad i=0,1,\ldots,M-1)$,
see Eq.~(\ref{eq:ScalarCoeff_ckj}) 
\end{enumerate}
The MATLAB code for determining the above parameters and an indicative
example at $M=3$ and $\rho_{\infty}=0.125$ are given in \ref{subsec:Appendix-Distinct-roots}. 

\subsubsection{Time-stepping scheme - distinct roots case}

For a chosen order $M$, the time-stepping equation Eq.~(\ref{eq:PF-sol})
requires evaluation of $M-1$ vectors $\mathbf{y}_{i}$. Multiplying
both sides of Eq.~(\ref{eq:DistinctRoot-yi}) with the matrix $r_{i}\mathbf{I}-\mathbf{A}$
and rearranging gives
\begin{equation}
(r_{i}\mathbf{I}-\mathbf{A})\mathbf{y}_{i}=\mathbf{g}^{(i)}+\left\{ \begin{array}{c}
\Delta t^{2}\mathbf{M^{-1}}\mathbf{f}_{\mathrm{r}i}\\
\mathbf{0}
\end{array}\right\} \,,\label{eq:DistinctRoot-step}
\end{equation}
which can be solved for each $\mathbf{y}_{i}$ (as demonstrated in
Section \ref{sec:Solution-to-the}). The following vectors $\mathbf{g}^{(i)}$
were introduced: 
\begin{equation}
\mathbf{g}^{(i)}=P_{\mathrm{L}}(r_{i})\mathbf{z}_{n-1}\,.\label{eq:DistinctRoot-stepg}
\end{equation}

In the case of a pair of complex conjugate roots, say $r_{i+1}=\bar{r}_{i},$
the solutions $\mathbf{y}_{i+1}=\bar{\mathbf{y}}_{i}$ will also be
complex conjugate. In this case it is sufficient to solve for $\mathbf{y}_{i}$
and use 
\begin{equation}
a_{i}\mathbf{y}_{i}+a_{i+1}\mathbf{y}_{i+1}=a_{i}\mathbf{y}_{i}+\bar{a}_{i}\bar{\mathbf{y}}_{i}=2\mathrm{Re}\left(a_{i}\mathbf{y}_{i}\right)\,,\label{eq:DistinctRoot-yconj}
\end{equation}
to obtain the solution for the conjugate root. 

\subsection{Case 2: Single multiple roots}\label{subsec:Single-multiple-roots}

\subsubsection{Evaluation of matrix exponential}

Alternatively, from Eq.~(\ref{eq:ExpmRational}), the matrix polynomial
$\mathbf{Q}$ in the denominator of the rational approximation can
be chosen to possess a single multiple root $r$ as in \citep{Song2024}.
Then, the rational approximation may be presented as:
\begin{equation}
e^{\mathbf{A}}\approx\dfrac{P(\mathbf{A})}{Q(\mathbf{A})}=\frac{p_{0}\mathbf{I}+p_{1}\mathbf{A}+\ldots+p_{M}\mathbf{A}^{M}}{(r\mathbf{I}-\mathbf{A})^{M}}\,.\label{eq:R-def}
\end{equation}

Here, the order of accuracy is equal to $M$, though the root $r$
is real and therefore complex solvers are not required. The value
of the spectral radius at the highest-frequency limit $\rho_{\infty}$
is defined by the user to match the coefficient of the highest-order
term in the numerator polynomial $P(\mathbf{A})$, i.e., $\rho_{\infty}=|p_{M}|$.
Note that the highest degree of the denominator is $M$ and cannot
be further reduced. The root $r$ is determined from generating a
polynomial equation (see \citep{Song2024} for details). The root
of this polynomial that leads to the spectral radius $\rho\leq1$
and the least relative period error in the low-frequency range is
then selected \citep{Song2024}. The remaining coefficients $p_{i}$
can be determined as functions of the root $r$ -- as demonstrated
in \ref{subsec:Appendix-Single-multiple-root}. In numerical algorithms
used to treat polynomials, it is often convenient to utilize Horner's
method. The polynomial $P(\mathbf{A})$ in Eq.~(\ref{eq:R-def})
is then expressed as 
\begin{equation}
P(\mathbf{A})=p_{0}\mathbf{I}+\mathbf{A}\left(p_{1}\mathbf{I}+\mathbf{A}\left(p_{2}\mathbf{I}\ldots+\mathbf{A}(p_{M-1}\mathbf{I}+p_{M}\mathbf{A})\ldots\right)\right)\,.\label{eq:HornerMethod}
\end{equation}

The polynomial can be manipulated recursively starting from the term
in the innermost parentheses. For conciseness, in the subsequent derivation
of the time-stepping algorithm in the form of Eq.~(\ref{eq:Stepping}),
the following matrix is introduced 
\begin{equation}
\mathbf{A}_{\mathrm{r}}=r\mathbf{I}-\mathbf{A}\,,\label{eq:Ar-def}
\end{equation}
which corresponds to the shift of the root of the polynomial.

Using Eq.~(\ref{eq:Ar-def}), the polynomials $\mathbf{P}=P(\mathbf{A})$
and $\mathbf{Q}=Q(\textbf{\ensuremath{\mathbf{A}}})$ in the rational
approximation form of Eq.~(\ref{eq:R-def}) are transformed to polynomials
in terms of the matrix $\mathbf{A}_{\mathrm{r}}$ as follows, by \emph{shifting}
the coefficients:
\begin{align}
P_{\mathrm{r}}(\mathbf{A}_{\mathrm{r}})\equiv P(\mathbf{A}) & =p_{0}\mathbf{I}+p_{1}(r\mathbf{I}-\textbf{A}_{r})+\ldots+p_{M}(r\mathbf{I}-\mathbf{A}_{\mathrm{r}})^{M}=\sum_{i=0}^{M}p_{\mathrm{r}i}\mathbf{A}_{\mathrm{r}}^{i}\,,\label{eq:Pr}\\
Q_{\mathrm{r}}(\mathbf{A}_{\mathrm{r}})\equiv Q(\textbf{\ensuremath{\mathbf{A}}}) & =\mathbf{A}_{\mathrm{r}}^{M}\,,\label{eq:Qr}
\end{align}
with shifted coefficients $p_{\mathrm{r}i}$ ($i=1,2,\ldots,M$) to
be determined accordingly. In the same way, the matrices $\mathbf{C}_{k}=\mathbf{C}_{k}(\mathbf{A})$
given in Eq.~(\ref{eq:Stepping_MatrixPolyCk}) are shifted and rewritten
as 
\begin{align}
\mathbf{C}_{rk}(\mathbf{A}_{r})=\sum_{i=0}^{M-1}c_{\mathrm{r}ki}\mathbf{A}_{\mathrm{r}}^{i}\equiv\mathbf{C}_{k}(\mathbf{A}) & =c_{k0}\mathbf{I}+c_{k1}(r\mathbf{I}-\textbf{A}_{r})+\ldots+c_{k(M-1)}(r\mathbf{I}-\textbf{A}_{r})^{M-1}\,.\label{eq:MatrixPoly_Crk}
\end{align}
The shifted scalar coefficients $c_{\mathrm{r}ki}$ are grouped and
stored as 
\begin{equation}
\mathbf{c}_{\mathrm{r}}=[c_{\mathrm{r}ki}]\,.\qquad(k=0,1,\ldots p_{\mathrm{f}};\quad i=0,1,\ldots,M-1)\,,\label{eq:ScalarCoeff_crki}
\end{equation}
analogously to those in Eq.~(\ref{eq:ScalarCoeff_ckj}) that operate
with $\mathbf{A}$.

The MATLAB function in Fig.~\ref{fig:MATLAB-function-shiftingroots}
performs the operation of shifting the polynomial roots in Eqs.~(\ref{eq:Pr})
and (\ref{eq:MatrixPoly_Crk}), returning $p_{\mathrm{r}i}$ and $c_{\mathrm{r}ki}$.
The algorithm is recursive, utilizing Horner's method in Eq.~(\ref{eq:HornerMethod}).

\begin{figure}

\texttt{function {[}prcoe{]} = shiftPolyCoefficients(pcoe,r)}

\texttt{M = size(pcoe,2) - 1; \%order of polynomial}

\texttt{zc = zeros(size(pcoe,1),1); \%vector of zeros}

\texttt{prcoe = pcoe;}

\texttt{for ii = M:-1:1}

\texttt{prcoe(:,ii:end) = {[}prcoe(:,ii)+r{*}prcoe(:,ii+1) r{*}prcoe(:,ii+2:end)
zc{]} ...}

\texttt{- {[}zc prcoe(:,ii+1:end){]};}

\texttt{end}\caption{MATLAB function to determine the coefficients of the polynomial $p_{\mathrm{r}}(x_{\mathrm{r}})$
resulting from shifting the roots of polynomial $p(x)$ with $x_{\mathrm{r}}=r-x$.
The input includes a row vector of size $M+1$ whose components are
the coefficients of $p(x)$.}\label{fig:MATLAB-function-shiftingroots}
\end{figure}

Substituting the polynomials $\mathbf{P}$, $\mathbf{Q}$ and matrices
$\mathbf{C}_{k}$ to equivalents in terms of the matrix $\mathbf{A}_{\mathrm{r}}$
(i.e., making use of Eqs.~(\ref{eq:Pr})-(\ref{eq:MatrixPoly_Crk})),
the time-stepping equation Eq.~(\ref{eq:Stepping}) is rewritten
as 
\begin{equation}
\mathbf{z}_{n}=\dfrac{\mathbf{P}_{\mathrm{r}}}{\mathbf{A}_{\mathrm{r}}^{M}}\mathbf{z}_{n-1}+\dfrac{1}{\mathbf{A}_{\mathrm{r}}^{M}}\sum\limits_{k=0}^{p_{\mathrm{f}}}\mathbf{C}_{\mathrm{r}k}\left\{ \begin{array}{c}
\Delta t^{2}\mathbf{M^{-1}}\tilde{\mathbf{f}}^{(k)}\\
\mathbf{0}
\end{array}\right\} \,.\label{eq:SteppingEq-Ar}
\end{equation}

Using Eq.~(\ref{eq:MatrixPoly_Crk}), the summation in the non-homogeneous
term may be written as a double-summation, and, given that the force
vectors $\tilde{\mathbf{f}}^{(k)}$ only depend on the polynomial
order $p_{\mathrm{f}}$ and are independent of the matrix $\mathbf{A}_{\mathrm{r}}$
(as indicated in Eq.~(\ref{eq:ForceExpansion})), the double-summation
is rearranged to consider the summation of the force vectors first:
\begin{equation}
\sum\limits_{k=0}^{p_{\mathrm{f}}}\mathbf{C}_{\mathrm{r}k}\left\{ \begin{array}{c}
\Delta t^{2}\mathbf{M^{-1}}\tilde{\mathbf{f}}^{(k)}\\
\mathbf{0}
\end{array}\right\} =\sum_{i=0}^{M-1}\mathbf{A}_{\mathrm{r}}^{i}\left\{ \begin{array}{c}
\Delta t^{2}\mathbf{M^{-1}}\mathbf{f}_{\mathrm{r}i}\\
\mathbf{0}
\end{array}\right\} \,,\label{eq:SteppingEq-Ar-SumRearrange}
\end{equation}
where, for simplicity in notation, vectors $\mathbf{f}_{\mathrm{r}i}$
representing the combination of the vectors $\tilde{\mathbf{f}}^{(k)}$
according to the coefficients $c_{\mathrm{r}ki}$ are introduced 
\begin{equation}
\mathbf{f}_{\mathrm{r}i}=\sum\limits_{k=0}^{p_{\mathrm{f}}}c_{\mathrm{r}ki}\tilde{\mathbf{f}}^{(k)}.\label{eq:fri}
\end{equation}

It is worthwhile to recall that, to ensure the order of their accuracy,
the coefficients $c_{\mathrm{r}ki}$ are determined using Eqs.~(\ref{eq:Stepping_C0})
and (\ref{eq:Stepping_MatrixPolyCk}). Each vector $\mathbf{f}_{\mathrm{r}i}$
is written in matrix form as 
\begin{equation}
\mathbf{f}_{\mathrm{r}i}=\tilde{\mathbf{F}}\cdot\mathbf{c}_{\mathrm{r}i}\,,\label{eq:SteppingEq-fri-mtx}
\end{equation}
where the matrix $\tilde{\mathbf{F}}$ was previously defined in Eq.~(\ref{eq:ForceFitting2}),
and the column vector $\mathbf{c}_{\mathrm{r}i}$ is defined correspondingly
as 
\begin{equation}
\mathbf{c}_{\mathrm{r}i}=\left[\begin{array}{ccccc}
c_{\mathrm{r}0i} & c_{\mathrm{r}1i} & \ldots & c_{\mathrm{r}2i} & c_{\mathrm{r}p_{\mathrm{f}}i}\end{array}\right]^{\mathrm{T}}\,.\label{eq:ScalarCoeff_cri}
\end{equation}

Now, substituting Eq.~(\ref{eq:SteppingEq-Ar-SumRearrange}) into
the time-stepping equation Eq.~(\ref{eq:SteppingEq-Ar}), provides
the following: 
\begin{equation}
\mathbf{z}_{n}=\dfrac{\mathbf{P}_{\mathrm{r}}}{\mathbf{A}_{\mathrm{r}}^{M}}\mathbf{z}_{n-1}+\dfrac{1}{\mathbf{A}_{\mathrm{r}}^{M}}\sum_{i=0}^{M-1}\mathbf{A}_{\mathrm{r}}^{i}\left\{ \begin{array}{c}
\Delta t^{2}\mathbf{M^{-1}}\mathbf{f}_{\mathrm{r}i}\\
\mathbf{0}
\end{array}\right\} \,,\label{eq:SteppingEq}
\end{equation}
and, after using the polynomial $\mathbf{P}_{\mathrm{r}}$ as in Eq.~(\ref{eq:Pr}),
the time-stepping solution using the partial fraction decomposition
is finally expressed as 
\begin{equation}
\mathbf{z}_{n}=p_{\mathrm{r}M}\mathbf{z}_{n-1}+\sum_{i=0}^{M-1}\dfrac{1}{\mathbf{A}_{\mathrm{r}}^{M-i}}\mathbf{b}_{i}\,,\label{eq:SteppingEq-final}
\end{equation}
with vectors $\mathbf{b}_{i}$ introduced as: 
\begin{equation}
\mathbf{b}_{i}=p_{\mathrm{r}i}\mathbf{z}_{n-1}+\left\{ \begin{array}{c}
\Delta t^{2}\mathbf{M^{-1}}\mathbf{f}_{\mathrm{r}i}\\
\mathbf{0}
\end{array}\right\} \,.\label{eq:SteppingEq-bi}
\end{equation}

To use the time-stepping scheme in Eq.\,(\ref{eq:SteppingEq-final}),
the user needs to specify just two parameters: $M$ (number of sub-steps)
and $\rho_{\infty}$ (spectral radius at the high-frequency limit).
The following quantities then need to be determined to initialize
the scheme:
\begin{enumerate}
\item the single multiple root $r$
\item coefficients of shifted polynomial $p_{\mathrm{r}i}$ ($i=0,\,1,\,2,\,\ldots,\,M$)
\item coefficient matrix used in determination of external excitation term
$\mathbf{F}_{\mathrm{r}}$, given by: $\mathbf{c}_{\mathrm{r}}=[c_{\mathrm{r}ki}];\quad(k=0,1,\ldots p_{\mathrm{f}};\quad i=0,1,\ldots,M-1)$
\end{enumerate}
Sample MATLAB code taking $M$ and $\rho_{\infty}$ as inputs to determine
the above parameters observed in Eq.~(\ref{eq:SteppingEq-final})
are given in Fig.~\ref{fig:MATLAB-function-CoeffForce}, \ref{subsec:Appendix-Single-multiple-root}.
The results for an indicative example at $M=3$ and $\rho_{\infty}=0.125$
are also provided.

\subsubsection{Time-stepping scheme - single multiple roots case}

To present an efficient algorithm for executing the derived time-stepping
scheme, the summation term in the time-stepping equation Eq.~(\ref{eq:SteppingEq-final})
may be expanded and by applying Horner's method (Eq.\,(\ref{eq:HornerMethod})),
this series of terms can then be grouped with the inverse of matrix
$\mathbf{A}_{\mathrm{r}}$ factored out to yield: 
\begin{equation}
\mathbf{z}_{n}=p_{\mathrm{r}M}\mathbf{z}_{n-1}+\dfrac{1}{\mathbf{A}_{\mathrm{r}}}\left(\mathbf{b}_{M-1}+\dfrac{1}{\mathbf{A}_{\mathrm{r}}}\left(\mathbf{b}_{M-2}+\ldots+\dfrac{1}{\mathbf{A}_{\mathrm{r}}}\left(\mathbf{b}_{2}+\dfrac{1}{\mathbf{A}_{\mathrm{r}}}\left(\mathbf{b}_{1}+\dfrac{1}{\mathbf{A}_{\mathrm{r}}}\mathbf{b}_{0}\right)\right)\right)\right)\,.\label{eq:SteppingEq-Horner}
\end{equation}

Introducing auxiliary variables $\mathbf{z}_{n}^{(i)}$, where $i=1,2,\ldots,M$,
the expression in Eq.~(\ref{eq:SteppingEq-Horner}) can be solved
recursively, according to
\begin{equation}
\mathbf{A}_{\mathrm{r}}\mathbf{z}_{n}^{(i+1)}=\mathbf{b}_{i}+\mathbf{z}_{n}^{(i)}\qquad(i=1,2,\ldots,M)\,.\label{eq:SingleRoot_recursiveEq}
\end{equation}

See \citep{Song2024} for details. This system of recursive equations
can be solved starting from $\mathbf{z}_{n}^{(0)}=\mathbf{0}$ to
the final term $\mathbf{z}_{n}^{(M)}$. Each step in this recursion
corresponds to solving an implicit equation for a given time step
in the time integration method. Following from Eq.~(\ref{eq:SteppingEq-Horner}),
the final solution for $\mathbf{z}_{n}$ is 
\begin{equation}
\mathbf{z}_{n}=p_{\mathrm{r}M}\mathbf{z}_{n-1}+\mathbf{z}_{n}^{(M)}\,.\label{eq:SteppingFinalSol}
\end{equation}
where $|p_{\mathrm{r}M}|=\rho_{\infty}$ \citep{Song2024}. Substituting
the definitions for $\mathbf{A}_{\mathrm{r}}$ in Eq.~(\ref{eq:Ar-def})
and $\mathbf{b}_{i}$ in Eq.~(\ref{eq:SteppingEq-bi}) into Eq.~(\ref{eq:SingleRoot_recursiveEq}),
the time-stepping equation is rewritten as:
\begin{equation}
(r\mathbf{I}-\mathbf{A})\mathbf{z}_{n}^{(i+1)}=\mathbf{g}^{(i)}+\left\{ \begin{array}{c}
\Delta t^{2}\mathbf{M^{-1}}\mathbf{f}_{\mathrm{r}i}\\
\mathbf{0}
\end{array}\right\} \,,\label{eq:SteppingWithA}
\end{equation}
with the following vectors $\mathbf{g}^{(i)}$ introduced 
\begin{equation}
\mathbf{g}^{(i)}=\mathbf{z}_{n}^{(i)}+p_{\mathrm{r}i}\mathbf{z}_{n-1}\,.\label{eq:SteppingSol_g}
\end{equation}

\section{Solution to the time-stepping equations}\label{sec:Solution-to-the}

The solution procedures for both the single multiple root case (Eq.~(\ref{eq:SteppingWithA}))
and distinct roots case (Eq.~(\ref{eq:DistinctRoot-step})) have
now each been presented in an identical form, which is generalized
as:
\begin{equation}
(r\mathbf{I}-\mathbf{A})\mathbf{x}=\mathbf{g}^{(i)}+\left\{ \begin{array}{c}
\Delta t^{2}\mathbf{M^{-1}}\mathbf{f}_{\mathrm{r}i}\\
\mathbf{0}
\end{array}\right\} \,.\label{eq:solution}
\end{equation}

For the case of distinct roots, $\mathbf{x}\equiv\mathbf{y}_{i}$,
which is calculated for $i=1,...,M$, and $\mathbf{g}^{(i)}$ is defined
by Eq.~(\ref{eq:DistinctRoot-stepg}). For the case of single multiple
roots, recall that $\mathbf{x}\equiv\mathbf{z}_{n}^{(i+1)}$ which
is solved for recursively from $i=0$ to $M-1$, and $\mathbf{g}^{(i)}$
is defined by the expression in Eq.~(\ref{eq:SteppingSol_g}). The
state-space vector $\mathbf{z}_{n}$ is then returned either by Eq.~(\ref{eq:SteppingFinalSol})
(single multiple root) or Eq.~(\ref{eq:PF-sol}) (distinct roots).

Equation~(\ref{eq:solution}) for both cases is solved by following
the procedure in \citep{Song2022}. The vectors $\mathbf{x}$ and
$\mathbf{g}^{(i)}$ can each be partitioned into two sub-vectors of
equal size 
\begin{equation}
\mathbf{x}=\begin{Bmatrix}\mathbf{x}_{1}\\
\mathbf{x}_{2}
\end{Bmatrix}\qquad\text{and}\qquad\mathbf{g}^{(i)}=\begin{Bmatrix}\mathbf{g}_{1}\\
\mathbf{g}_{2}
\end{Bmatrix}\,,\label{eq:RootPartition}
\end{equation}
then, making use of Eq.~(\ref{eq:matrix=000020A}), Eq.~(\ref{eq:solution})
is rewritten as
\begin{equation}
\begin{bmatrix}r\mathbf{I}+\Delta t\mathbf{M}^{-1}\mathbf{C} & \Delta t^{2}\mathbf{M}^{-1}\mathbf{K}\\
-\mathbf{I} & r\mathbf{I}
\end{bmatrix}\begin{Bmatrix}\mathbf{x}_{1}\\
\mathbf{x}_{2}
\end{Bmatrix}=\begin{Bmatrix}\mathbf{g}_{1}\\
\mathbf{g}_{2}
\end{Bmatrix}+\left\{ \begin{array}{c}
\Delta t^{2}\mathbf{M^{-1}}\mathbf{f}_{\mathrm{r}i}\\
\mathbf{0}
\end{array}\right\} \,.\label{eq:RootEq}
\end{equation}

The solution to $\mathbf{x}_{1}$ is given by 
\begin{equation}
\left(r^{2}\mathbf{M}+r\Delta t\mathbf{C}+\Delta t^{2}\mathbf{K}\right)\mathbf{x}_{1}=r\mathbf{M}\mathbf{g}_{1}-\Delta t^{2}\mathbf{K}\mathbf{g}_{2}+r\Delta t^{2}\mathbf{f}_{\mathrm{r}i}\,.\label{eq:realRootEq5}
\end{equation}
The solution to $\mathbf{x}_{2}$ follows as 
\begin{equation}
r\mathbf{x}_{2}=\mathbf{x}_{1}+\mathbf{g}_{2}\,.\label{eq:realRootEq4}
\end{equation}

Algorithms for executing the time-stepping procedure for partial fraction
expansions using distinct roots, or a single multiple root, are provided
in Tables~\ref{tab:Time-integration-solution_distinct} and \ref{tab:Time-integration-solution_multiple},
respectively. Sample MATLAB code for each algorithm are provided in
Fig.\,\ref{fig:MATLAB-function-TimeSolverTRhoPF} and Fig.\,\ref{fig:MATLAB-function-TimeSolverPF},
respectively.

\section{Calculation of acceleration}\label{sec:Calculation-of-acceleration}

Although the self-starting algorithm presented in this work does not
need acceleration in time-stepping solution detailed in the previous
section, acceleration can be computed to the same order of accuracy
as displacement and velocity using only vector operations, preserving
computational efficiency. Note that, for both the case of distinct
roots and single multiple root, the vectors $\mathbf{x}_{1}$ and
$\mathbf{\mathbf{x}}_{2}$ (and, similarly $\mathbf{g}_{1}$ and $\mathbf{g}_{2}$)
in Eq.~(\ref{eq:RootEq}) relate to velocity and displacement measures,
respectively. Therefore, considering dimensionless time-derivatives
of these measures, it holds that
\begin{equation}
\circdot{\mathbf{x}}_{2}=\mathbf{x}_{1},\,\,\,\,\,\,\circdot{\mathbf{g}}_{2}=\mathbf{g}_{1}\,.\label{eq:circdotx2}
\end{equation}

Consider taking the dimensionless time-derivative of the time-stepping
function Eq.~(\ref{eq:PF-sol}). The upper partition of this provides
the equation
\begin{equation}
\acco_{n}=\rho\acco_{n-1}+\sum_{i=1}^{M}a_{i}\circdot{\mathbf{x}}_{1}^{(i)}\,,
\end{equation}
with $\mathbf{x}_{1}^{(i)}$ referring to the $i$-th calculation
of vector $\mathbf{x}_{1}$. Consider now taking the dimensionless
time-derivative of Eq.~(\ref{eq:realRootEq4}) such that
\begin{equation}
r\circdot{\mathbf{x}}_{2}=\circdot{\mathbf{x}}_{1}+\circdot{\mathbf{g}}_{2}\,.
\end{equation}

Using Eq.~(\ref{eq:circdotx2}) it then follows that
\begin{equation}
\circdot{\mathbf{x}}_{1}=r\mathbf{x}_{1}-\mathbf{g}_{1}\,,
\end{equation}
where the quantities on the right-hand side of this expression are
already determined during the time-stepping procedure and not re-calculated.
Therefore, calculation of acceleration for any time step $n>0$ follows
\begin{equation}
\acco_{n}=\rho\acco_{n-1}+\sum_{i=1}^{M}a_{i}\left(r_{i}\mathbf{x}_{1}^{(i)}-\mathbf{g}_{1}^{(i)}\right)\,.\label{eq:acceleration-distinct}
\end{equation}

For the case of a single multiple root, an analogous argument is made,
and acceleration is determined according to
\begin{equation}
\acco_{n}=p_{\mathrm{r}M}\acco_{n-1}+\left(r\mathbf{x}_{1}^{(M)}-\mathbf{g}_{1}^{(M)}\right)\,,\label{eq:acceleration-multiple}
\end{equation}
which requires knowledge of only the final calculation of vectors
$\mathbf{x}_{1}$ and $\mathbf{g}_{1}$ during the time-stepping procedure.
The algorithm does not need any additional equation solutions, and
can be extended to find higher-order derivatives if desired. Significantly,
the accuracy of the acceleration is the same order as displacement
and velocity without excessive additional computational cost. It is
worthwhile to emphasize that:
\begin{enumerate}
\item The calculation of acceleration is an optional (effectively post-processing)
operation.
\item In the case where the user sets $\rho_{\infty}=0$ (hence, $\rho=0$
in Eq.\,(\ref{eq:acceleration-distinct}) or $p_{\mathrm{r}M}=0$
in Eq.~(\ref{eq:acceleration-multiple})), the calculation of acceleration
does not involve the acceleration of previous time steps, and the
initial acceleration $\acco_{0}$ therefore need not be known.
\item In many problems the the initial internal and external forces are
balanced. In these cases, acceleration is determined at all time steps
using only vector operations since the initial acceleration $\acco_{0}$
need not be known. Else, a single additional equation is solved prior
to the time-stepping solution procedure to find the initial acceleration.
\end{enumerate}

\section{Calculation of force derivatives}\label{sec:Calculation-of-force-deriv}

As outlined in Section \ref{subsec:Analytical-approximation-of},
the non-homogeneous term is treated in this work using a series expansion
of the (in general, nonlinear) force vector $\mathbf{f}$. The forces
are sampled at Gauss-Lobatto points $0=s_{1}\leq\dots\leq s_{N}=1$
to perform the numerical integration. The order of accuracy of the
Gauss-Lobatto quadrature with $N$ sampling points is equal to $2N-3$.
A sufficient number of points $N$ is selected as to ensure the order
of accuracy provided from the rational approximation of the matrix
exponential governs the accuracy of the scheme.  The force vectors
at these discrete times are expressed using Eq.~(\ref{eq:ForceFitting1})
and are grouped in the matrix
\begin{equation}
\mathbf{F}_{\mathrm{p}}=\left[\begin{array}{cccc}
\mathbf{f}(s_{1}) & \mathbf{f}(s_{2}) & \ldots & \mathbf{f}(s_{N})\end{array}\right]=\tilde{\mathbf{F}}\cdot\mathbf{V}\,,\label{eq:ForceFitting4}
\end{equation}
where
\begin{equation}
\mathbf{V}=[\begin{array}{cccc}
\mathbf{s}(s_{1}) & \mathbf{s}(s_{2}) & \dots & \mathbf{s}(s_{N})\end{array}]=\left[\begin{array}{cccc}
1 & 1 & \ldots & 1\\
(s_{1}-0.5) & (s_{2}-0.5) & \ldots & (s_{N}-0.5)\\
(s_{1}-0.5)^{2} & (s_{2}-0.5)^{2} & \ldots & (s_{N}-0.5)^{2}\\
\vdots & \vdots & \ddots & \vdots\\
(s_{1}-0.5)^{N-1} & (s_{2}-0.5)^{N-1} & \ldots & (s_{N}-0.5)^{N-1}
\end{array}\right]\label{eq:ForceFitting5}
\end{equation}
is constructed with the use of Eq.~(\ref{eq:ForceFitting3}). The
matrix of coefficients $\tilde{\mathbf{F}}$ of the polynomial expansion
in Eq.~(\ref{eq:ForceExpansion}) is obtained from the matrix of
force vectors at the discrete times $\mathbf{F}_{\mathrm{p}}$ as
\begin{equation}
\tilde{\mathbf{F}}=\mathbf{F}_{\mathrm{p}}\cdot\mathbf{T}_{\mathrm{p}}\,,\thinspace\thinspace\thinspace\mathbf{T}_{\mathrm{p}}=\mathbf{V}^{-\mathrm{T}}\thinspace.\label{eq:ForceFitting6}
\end{equation}
The coefficient of each term of the polynomial expansion in Eq.~(\ref{eq:ForceExpansion})
is then written as 
\begin{equation}
\tilde{\mathbf{f}}^{(k)}=\mathbf{F}_{\mathrm{p}}\cdot\mathbf{T}_{\mathrm{p}}^{(k)}\,,
\end{equation}
where $\mathbf{T}_{\mathrm{p}}^{(k)}$ is the $k$-th column of the
matrix $\mathbf{T}_{\mathrm{p}}$. For the case of distinct roots,
each vector $\mathbf{f}_{\mathrm{r}i}$ in Eq.~(\ref{eq:PF-nonhomo-fri-mtx})
can now alternatively be expressed in terms of the matrix of force-vectors
$\mathbf{F}_{\mathrm{p}}$ sampled at discrete times (see Eq.~(\ref{eq:ForceFitting4}))
\begin{equation}
\mathbf{f}_{\mathrm{r}i}=\mathbf{F}_{\mathrm{p}}\mathbf{T}_{\mathrm{p}}\mathbf{c}\mathbf{r}_{i}\,.\label{eq:PF-nonhomo-fri-discrete}
\end{equation}

For the case of a single multiple root, the equivalent vectors $\mathbf{f}_{\mathrm{r}i}$
were defined through Eq.~(\ref{eq:SteppingEq-fri-mtx}), and can
now be expressed in terms of the sampled forces $\mathbf{F}_{\mathrm{p}}$
by
\begin{equation}
\mathbf{f}_{\mathrm{r}i}=\mathbf{F}_{\mathrm{p}}\mathbf{T}_{\mathrm{p}}\mathbf{c}_{\mathrm{r}i}\,.\label{eq:fri_multiple_final}
\end{equation}

\subsection{Interpolation of displacements}\label{subsec:Interpolation-of-displacements}

For nonlinear problems, the sampled forces required to generate the
matrix $\mathbf{F}_{p}$ in Eq.~(\ref{eq:ForceFitting4}) are functions
of the state-variables $\mathbf{z}(t)$ at times $t(s)=t_{n-1}+s\Delta t$,
however these variables are only known at the beginning of the time
step (i.e. $\mathbf{z}_{n-1})$. During an iterative solution procedure,
successive estimates of the state at the end of the time step $\mathbf{z}_{n}$
are made. For the first iteration, a standard Taylor series extrapolation
is used to produce the initial estimate for $\mathbf{z}_{n}$. For
each subsequent iteration, using known values of displacement (and
additional time-derivatives) at each end of the time step, a Hermite
polynomial can be used to approximate states at the internal points
\begin{equation}
\underline{\mathbf{\mathbf{u}}}(t_{n-1}+s\Delta t)=\sum_{k=0}^{p_{n-1}}\alpha_{k}(s)\mathbf{u}_{n-1}^{(k)}+\sum_{k=0}^{p_{n}}\beta_{k}(s)\mathbf{u}_{n}^{(k)}\,,\label{eq:displacement-interpolation}
\end{equation}
with $\alpha_{k}(s),$$\beta_{k}(s)$ denoting the coefficients of
the interpolated polynomial for known dimensionless time-derivatives
of the displacement at the beginning and end of the time step, respectively.
The resulting polynomial is of order $p_{n-1}+p_{n}+1$, and truncation
error in the time-stepping procedure related to interpolation of state-variables
is $O\left(\Delta t^{p_{n-1}+p_{n}+3}\right)$. The velocity at $t(s)=t_{n-1}+s\Delta t$
is approximated using the derivative of the polynomial in Eq.~(\ref{eq:displacement-interpolation}).
In the proposed algorithms, acceleration is evaluated with negligible
computational effort as demonstrated in Section \ref{sec:Calculation-of-acceleration},
and therefore may be included in the Hermite interpolation scheme
without significant computational cost. Hence, we take $p_{n-1}=p_{n}=2$
, and a quintic polynomial interpolates the displacement field for
any given iteration. The time-stepping procedures presented in this
work therefore have maximum seventh-order accuracy using this interpolation
scheme. This restriction of order of accuracy only relates to nonlinear
problems; for linear transient problems sampled forces are not functions
of unknown state-variables and therefore no interpolation is required.
To increase accuracy beyond seventh-order for nonlinear problems,
higher-order derivatives are required at the beginning and end of
the time step, which may be evaluated through a simple extension of
the method used to calculate acceleration presented in Section \ref{sec:Calculation-of-acceleration},
though is left to a future study. Alternate interpolation schemes,
for example considering solved state-variables from previous time
steps, are also possible.

\subsection{Truncation error for proposed algorithms}

The time-stepping algorithms presented in this work possess error
from two main sources:
\begin{enumerate}
\item Rational approximation of the matrix exponential function,
\item Numerical integration of the non-homogeneous term.
\end{enumerate}
The aim is to integrate the non-homogeneous term to such accuracy
that the order of polynomials used in the rational approximation of
the matrix exponential governs the overall solution accuracy. For
the case of distinct roots, use of $M$ sub-steps provides an order
of accuracy of $2M$ in the absence of numerical dissipation ($\rho_{\infty}=1)$,
and $2M-1$ when numerical dissipation is introduced ($\rho_{\infty}<1)$.
Such schemes are labeled \emph{Pade-PF}\emph{:(M-1)M} in the examples
to follow. For the algorithms using a single multiple root, use of
$M$ sub-steps provides an order of accuracy of $M$ for the matrix
exponential. Such schemes are labeled \emph{M-PF:(M) }in the examples
to follow. For both the cases of distinct roots or single multiple
roots, it is found that the theoretical convergence rates governed
by the rational approximation of the matrix exponential are able to
be met if the series expansion of the nonlinear force-vector is taken
as $N=M+1$ .

\section{Time-stepping solution algorithms}

Table \ref{tab:Time-integration-solution_distinct} provides the algorithm
for the time-stepping solution procedure for nonlinear problems for
the case of distinct roots $r_{i}$ used in the factorization of the
divisor polynomial $\mathbf{Q}$. Note that lines 7, 21, 26, 28 are
optional and not required to complete the time-stepping solution procedure.
\begin{table}
\begin{tabular}{|c|l|l|}
\hline 
1 & \emph{define:} $\rho_{\infty}\in\left[0,1\right]$, order $M$, integration
points $N$, time step $\Delta t$ & \tabularnewline
2 & \emph{input:}\textbf{ }problem definition $\mathbf{M},\thinspace\mathbf{C},\thinspace\mathbf{K},\thinspace\mathbf{f}(t)$;
initial conditions $\mathbf{z}_{0}$ & \tabularnewline
3 & \emph{compute:} radius $\rho=p_{M}/q_{M}$; roots $r_{i}$; coefficients
$p_{\mathrm{L}i}$, $a_{i}$, $\mathbf{c}$ & Eqs.\,(\ref{eq:PF-expansionSln}, \ref{eq:Coefficient-pli}, \ref{eq:ScalarCoeff_ckj})\tabularnewline
4 & \emph{define:} interpolation points $s_{k}$, for $k=0,1,2,...,N$ & \tabularnewline
5 & \emph{compute:} transformation matrix $\mathbf{T}_{p}$ & Eq.\,(\ref{eq:ForceFitting6})\tabularnewline
6 & \emph{store:} effective stiffness $\tilde{\mathbf{K}}\leftarrow r^{2}\mathbf{M}+r\Delta t\mathbf{C}+\Delta t^{2}\mathbf{K}$ & \tabularnewline
7 & \emph{compute:} initial acceleration $\acco_{0}=\mathbf{M}^{-1}\left(\mathbf{f}_{E}(0)-\mathbf{f}_{I}(\mathbf{z}_{0})\right)$\,\,(only
if $\mathbf{\rho}_{\infty}\neq\mathbf{0}$) & Eq.\,(\ref{eq:EqofMotionNonlinear})\tabularnewline
8 & \emph{define:}\textbf{ $t_{n}=0$} & \tabularnewline
 & \textbf{for }each time step $n=1,2,3,...$ & \tabularnewline
9 & \emph{\,\,\,\,\,\,store:} effective stiffness $\tilde{\mathbf{K}}\leftarrow r^{2}\mathbf{M}+r\Delta t\mathbf{C}_{n-1}+\Delta t^{2}\mathbf{K}_{n-1}$ & \tabularnewline
10 & \,\,\,\,\,\,\emph{initialize:} target state-variables $\mathbf{z}_{n}\leftarrow\mathbf{z}_{n-1}+\circdot{\mathbf{z}}_{n-1}$ & \tabularnewline
11 & \,\,\,\,\,\,\emph{define:}\textbf{ }time vector $t_{k}=t_{n-1}+s_{k}\Delta t$ & \tabularnewline
 & \,\,\,\,\,\,\textbf{for} iterations $j=1,2,...,\mathrm{max\thinspace iter.}$ & \tabularnewline
12 & \emph{\,\,\,\,\,\,\,\,\,\,\,\,compute:} interpolated state-variables
$\mathbf{u}(t_{k})=\underline{\mathbf{u}}(t_{k})$, $\circdot{\mathbf{{u}}}(t_{k})=\circdot{\underline{{\mathbf{{u}}}}}(t_{k})$ & Eq.\,(\ref{eq:displacement-interpolation})\tabularnewline
13 & \,\,\,\,\,\,\,\,\,\,\,\,\emph{store:}\textbf{ }matrix
of \textbf{s}ampled forces $\mathbf{F}_{\mathrm{p}}\leftarrow\mathbf{f}(\mathbf{z}(t_{k}))$ & Eq.\,(\ref{eq:ForceFitting4})\tabularnewline
14 & \textbf{\,\,\,\,\,\,\,\,\,\,\,\,}\emph{compute:} force
vectors $\mathbf{f}_{\mathrm{r}i}\leftarrow\mathbf{F}_{\mathrm{p}}\cdot\mathbf{T}\cdot\mathbf{c}\cdot\mathbf{r}_{i}$ & Eq.\,(\ref{eq:PF-nonhomo-fri-discrete})\tabularnewline
15 & \,\,\,\,\,\,\,\,\,\,\,\,\emph{define:}\textbf{ }$\mathbf{z}_{n}\leftarrow\rho\mathbf{z}_{n-1}$,
and $\overset{\circ\circ}{\mathbf{u}}_{n}\leftarrow\rho\overset{\circ\circ}{\mathbf{u}}_{n-1}$ & Eqs.\,(\ref{eq:PF-sol}, \ref{eq:acceleration-distinct})\tabularnewline
 & \,\,\,\,\,\,\,\,\,\,\,\,\textbf{for:} $i=1$ to $n_{\mathrm{real}}$
(the real roots) & \tabularnewline
16 & \,\,\,\,\,\,\,\,\,\,\,\,\,\,\,\,\,\,\emph{store:}
$r_{i}\leftarrow\mathrm{Re}(r_{i})$ & \tabularnewline
17 & \,\,\,\,\,\,\,\,\,\,\,\,\,\,\,\,\,\,\emph{compute:}\textbf{
$\mathbf{g}=\left[\mathbf{g}_{1}\thinspace;\thinspace\thinspace\mathbf{g}_{2}\right]\leftarrow\mathrm{Re}\left(p_{\mathrm{L}i}\right)\mathbf{z}_{n-1}$} & Eq.\,(\ref{eq:DistinctRoot-stepg})\tabularnewline
18 & \,\,\,\,\,\,\,\,\,\,\,\,\,\,\,\,\,\,\emph{solve:}\textbf{
$\mathbf{x}_{1}\leftarrow\tilde{\mathbf{K}}^{-1}\cdot\left(r\mathbf{M}\cdot\mathbf{g}_{1}-\Delta t^{2}\mathbf{K}\cdot\mathbf{g}_{2}+r\Delta t^{2}\mathbf{f}_{\mathrm{r}i}\right)$} & Eq.\,(\ref{eq:realRootEq5})\tabularnewline
19 & \,\,\,\,\,\,\,\,\,\,\,\,\,\,\,\,\,\,\emph{solve:}\textbf{
$\mathbf{x}_{2}\leftarrow\left(\mathbf{x}_{1}+\mathbf{g}_{2}\right)/r$} & Eq.\,(\ref{eq:realRootEq4})\tabularnewline
20 & \,\,\,\,\,\,\,\,\,\,\,\,\,\,\,\,\,\,\emph{update:}\textbf{
$\mathbf{z}_{n}\leftarrow\mathbf{z}_{n}+\mathrm{Re}\left(a_{i}\right)\left[\mathbf{x}_{1}\thinspace;\thinspace\thinspace\mathbf{x}_{2}\right]$} & Eq.\,(\ref{eq:PF-sol})\tabularnewline
21 & \,\,\,\,\,\,\,\,\,\,\,\,\,\,\,\,\,\,\emph{update:}\textbf{
$\overset{\circ\circ}{\mathbf{u}}_{n}\leftarrow\overset{\circ\circ}{\mathbf{u}}_{n}+\mathrm{Re}\left(a_{i}\right)\left(r\mathbf{x}_{1}-\mathbf{g}_{1}\right)$} & Eq.\,(\ref{eq:acceleration-distinct})\tabularnewline
 & \,\,\,\,\,\,\,\,\,\,\,\,\textbf{end} & \tabularnewline
 & \,\,\,\,\,\,\,\,\,\,\,\,\textbf{for} $i=(n_{\mathrm{real}}+1)$
to $(n_{\mathrm{complex}}/2)$ (the sets of complex conjugate roots) & \tabularnewline
22 & \,\,\,\,\,\,\,\,\,\,\,\,\,\,\,\,\,\,\emph{compute:}\textbf{
$\mathbf{g}=\left[\mathbf{g}_{1}\thinspace;\thinspace\thinspace\mathbf{g}_{2}\right]\leftarrow p_{\mathrm{L}i}\mathbf{z}_{n-1}$} & Eq.\,(\ref{eq:DistinctRoot-stepg})\tabularnewline
23 & \,\,\,\,\,\,\,\,\,\,\,\,\,\,\,\,\,\,\emph{solve:}\textbf{
$\mathbf{x}_{1}\leftarrow\tilde{\mathbf{K}}^{-1}\cdot\left(r\mathbf{M}\cdot\mathbf{g}_{1}-\Delta t^{2}\mathbf{K}\cdot\mathbf{g}_{2}+r\Delta t^{2}\mathbf{F}_{\mathrm{r}[i]}\right)$} & Eq.\,(\ref{eq:realRootEq5})\tabularnewline
24 & \,\,\,\,\,\,\,\,\,\,\,\,\,\,\,\,\,\,\emph{solve:}\textbf{
$\mathbf{x}_{2}\leftarrow\left(\mathbf{x}_{1}+\mathbf{g}_{2}\right)/r$} & Eq.\,(\ref{eq:realRootEq4})\tabularnewline
25 & \,\,\,\,\,\,\,\,\,\,\,\,\,\,\,\,\,\,\emph{update:}\textbf{
$\mathbf{z}_{n}\leftarrow\mathbf{z}_{n}+2\mathrm{Re}\left(a_{i}\right)\left[\mathbf{x}_{1}\thinspace;\thinspace\thinspace\mathbf{x}_{2}\right]$} & Eqs.\,(\ref{eq:PF-sol}, \ref{eq:DistinctRoot-yconj})\tabularnewline
26 & \,\,\,\,\,\,\,\,\,\,\,\,\,\,\,\,\,\,\emph{update:}\textbf{
$\overset{\circ\circ}{\mathbf{u}}_{n}\leftarrow\overset{\circ\circ}{\mathbf{u}}_{n}+2\mathrm{Re}\left(a_{i}\left(r\mathbf{x}_{1}-\mathbf{g}_{1}\right)\right)$} & Eq.\,(\ref{eq:acceleration-distinct})\tabularnewline
 & \,\,\,\,\,\,\,\,\,\,\,\,\textbf{end} & \tabularnewline
 & \,\,\,\,\,\,\,\,\,\,\,\,\textbf{if }convergence met & \tabularnewline
27 & \,\,\,\,\,\,\,\,\,\,\,\,\,\,\,\,\,\,\emph{output:}
solution $\mathbf{z}_{n}$ & \tabularnewline
28 & \,\,\,\,\,\,\,\,\,\,\,\,\,\,\,\,\,\,\emph{output:}
acceleration $\acco_{n}$ & \tabularnewline
29 & \,\,\,\,\,\,\,\,\,\,\,\,\,\,\,\,\,\,\emph{update:}
time step $t_{n}\leftarrow t_{n}+\Delta t$ & \tabularnewline
 & \,\,\,\,\,\,\,\,\,\,\,\,\,\,\,\,\,\textbf{break} & \tabularnewline
 & \,\,\,\,\,\,\,\,\,\,\,\,\textbf{else} & \tabularnewline
30 & \,\,\,\,\,\,\,\,\,\,\,\,\,\,\,\,\,\,\emph{return
to }12 & \tabularnewline
 & \,\,\,\,\,\,\,\,\,\,\,\,\textbf{end} & \tabularnewline
 & \textbf{end} & \tabularnewline
\hline 
\end{tabular}\caption{Time integration solution algorithm for nonlinear analyses using partial
fraction expansion with distinct roots }\label{tab:Time-integration-solution_distinct}
\end{table}
Table \ref{tab:Time-integration-solution_multiple} provides the algorithm
for the time-stepping solution procedure for nonlinear problems using
a partial fraction expansion with a single multiple root $r$. Note
that lines 7 and 20 are optional and are not required to complete
the time-stepping procedure. These need only be executed if the user
desires acceleration to be output to the same order of accuracy as
displacement and velocity using only vector operations. Indicative
MATLAB code for execution of both algorithms for linear problems is
provided in \ref{sec:Sample-MATLAB-code}, which is simply extended
to the nonlinear algorithm in Table \ref{tab:Time-integration-solution_distinct}.

\begin{table}
\begin{tabular}{|c|l|l|}
\hline 
1 & \emph{define}\textbf{\emph{:}} $\rho_{\infty}\in\left[0,1\right]$,
order $M$, integration points $N$, time step $\Delta t$ & \tabularnewline
2 & \emph{input:}\textbf{ }problem definition $\mathbf{M},\thinspace\mathbf{C}(t),\thinspace\mathbf{K}(t),\thinspace\mathbf{f}(t)$,
initial conditions $\mathbf{z}_{0}$ & Eq.\,(\ref{eq:linear=000020eq=000020of=000020motion})\tabularnewline
3 & \emph{compute:} root $r$, shifted coefficients $p_{\mathrm{r}i}$,
$\mathbf{c}_{\mathrm{r}}$ & Eqs.\,(\ref{eq:Pr}, \ref{eq:ScalarCoeff_cri})\tabularnewline
4 & \emph{define:} interpolation points $s_{k}$, for $k=0,1,2,...,N$ & \tabularnewline
5 & \emph{compute:} transformation matrix $\mathbf{T}_{p}$ & Eq.\,(\ref{eq:ForceFitting6})\tabularnewline
7 & \emph{compute:} initial acceleration $\acco_{0}=\mathbf{M}^{-1}\left(\mathbf{f}_{E}(0)-\mathbf{f}_{I}(\mathbf{z}_{0})\right)$\,\,(only
if $\mathbf{\rho}_{\infty}\neq\mathbf{0}$) & Eq.\,(\ref{eq:EqofMotionNonlinear})\tabularnewline
8 & \emph{define:}\textbf{ $t_{n}\leftarrow0$}, $\mathbf{z}_{n}^{(0)}\leftarrow\mathbf{0}$ & \tabularnewline
 & \textbf{for }each time step $n=1,2,3,...$ & \tabularnewline
9 & \emph{\,\,\,\,\,\,store:} effective stiffness $\tilde{\mathbf{K}}\leftarrow r^{2}\mathbf{M}+r\Delta t\mathbf{C}_{n-1}+\Delta t^{2}\mathbf{K}_{n-1}$ & \tabularnewline
10 & \,\,\,\,\,\,\emph{initialize:} target state-variables $\mathbf{z}_{n}\leftarrow\mathbf{z}_{n-1}+\circdot{\mathbf{z}}_{n-1}$ & \tabularnewline
11 & \,\,\,\,\,\,\emph{define:}\textbf{ }time vector $t_{k}\leftarrow t_{n-1}+s_{k}\Delta t$ & \tabularnewline
 & \,\,\,\,\,\,\textbf{for} iterations $j=1,2,...,\mathrm{max\thinspace iter.}$ & \tabularnewline
12 & \,\,\,\,\,\,\,\,\,\,\,\,\emph{compute:} interpolated state-variables
$\mathbf{u}(t_{k})=\underline{\mathbf{u}}(t_{k})$, $\circdot{\mathbf{{u}}}(t_{k})=\circdot{\underline{{\mathbf{{u}}}}}(t_{k})$ & Eq.\,(\ref{eq:displacement-interpolation})\tabularnewline
13 & \,\,\,\,\,\,\,\,\,\,\,\,\emph{store:}\textbf{ }matrix
of \textbf{s}ampled forces $\mathbf{F}_{\mathrm{p}}\leftarrow\mathbf{f}(\mathbf{z}(t_{k}))$ & Eq.\,(\ref{eq:ForceFitting4})\tabularnewline
14 & \textbf{\,\,\,\,\,\,\,\,\,\,\,\,}\emph{compute:} force
vectors $\mathbf{f}_{\mathrm{r}i}\leftarrow\mathbf{F}_{\mathrm{p}}\cdot\mathbf{T}\cdot\mathbf{c}_{\mathrm{r}i}$ & Eq.\,(\ref{eq:fri_multiple_final})\tabularnewline
 & \,\,\,\,\,\,\,\,\,\,\,\,\textbf{for }each order from $i=0$
to $M-1$ & \tabularnewline
15 & \,\,\,\,\,\,\,\,\,\,\,\,\,\,\,\,\,\,\emph{store:}\textbf{
}$\left[\mathbf{g}_{1};\thinspace\thinspace\mathbf{g}_{2}\right]\leftarrow\mathbf{z}_{n}^{(i)}+p_{\mathrm{r}i}\mathbf{z}_{n-1}$ & Eq.\,(\ref{eq:SteppingSol_g})\tabularnewline
16 & \,\,\,\,\,\,\,\,\,\,\,\,\,\,\,\,\,\,\emph{solve:}\textbf{
}$\mathbf{x}_{1}\leftarrow\tilde{\mathbf{K}}^{-1}\cdot\left(r\mathbf{M}\cdot\mathbf{g}_{1}-\Delta t^{2}\mathbf{K}_{n-1}\cdot\mathbf{g}_{2}+r\Delta t^{2}\mathbf{f}_{\mathrm{r}i}\right)$ & Eq.\,(\ref{eq:realRootEq5})\tabularnewline
17 & \,\,\,\,\,\,\,\,\,\,\,\,\,\,\,\,\,\,\emph{solve:}\textbf{
$\mathbf{x}_{2}\leftarrow\left(\mathbf{x}_{1}+\mathbf{g}_{2}\right)/r$} & Eq.\,(\ref{eq:realRootEq4})\tabularnewline
18 & \,\,\,\,\,\,\,\,\,\,\,\,\,\,\,\,\,\,\emph{store:
$\mathbf{z}_{n}^{(i+1)}\leftarrow\left[\mathbf{x}_{1}\thinspace;\thinspace\thinspace\mathbf{x}_{2}\right]$} & \tabularnewline
 & \,\,\,\,\,\,\,\,\,\,\,\,\textbf{end} & \tabularnewline
19 & \,\,\,\,\,\,\,\,\,\,\,\,\emph{update:} time step solution
$\mathbf{z}_{n}\leftarrow p_{\mathrm{r}M}\mathbf{z}_{n-1}+\mathbf{z}_{n}^{(M)}$ & Eq.\,(\ref{eq:SteppingFinalSol})\tabularnewline
 & \,\,\,\,\,\,\,\,\,\,\,\,\textbf{if }convergence met & \tabularnewline
 & \,\,\,\,\,\,\,\,\,\,\,\,\,\,\,\,\,\,\emph{output:}
solution $\mathbf{z}_{n}$ & \tabularnewline
20 & \,\,\,\,\,\,\,\,\,\,\,\,\,\,\,\,\,\,\emph{output:}
acceleration $\acco_{n}\leftarrow p_{\mathrm{r}M}\acco_{n-1}+\left(r\mathbf{x}_{1}-\mathbf{g}_{1}\right)$ & Eq.\,(\ref{eq:acceleration-multiple})\tabularnewline
21 & \,\,\,\,\,\,\,\,\,\,\,\,\,\,\,\,\,\,\emph{update:}
time step $t_{n}\leftarrow t_{n}+\Delta t$ & \tabularnewline
 & \,\,\,\,\,\,\,\,\,\,\,\,\,\,\,\,\,\textbf{break} & \tabularnewline
 & \,\,\,\,\,\,\,\,\,\,\,\,\textbf{else} & \tabularnewline
22 & \,\,\,\,\,\,\,\,\,\,\,\,\,\,\,\,\,\,\emph{return
to }12 & \tabularnewline
 & \,\,\,\,\,\,\,\,\,\,\,\,\textbf{end} & \tabularnewline
 & \textbf{end} & \tabularnewline
\hline 
\end{tabular}\caption{Time integration solution algorithm for nonlinear analyses using partial
fraction expansion with single multiple root }\label{tab:Time-integration-solution_multiple}
\end{table}

\section{Numerical Examples}\label{sec:Numerical-Examples}

In this section, the algorithms designed using the partial fraction
decompositions with both single multiple roots and distinct roots
are employed to solve a range of numerical examples that are commonly
used in studies of implicit time integration methods. The numerical
examples selected focus on demonstrating the accuracy (and convergence)
of efficiently solving displacement, velocity and acceleration in
both linear and highly-nonlinear problems, as well as the effectiveness
of the high-order schemes in suppressing spurious high-frequency oscillations.
Computational performance of the high-order time-stepping schemes
reported in previous works (see \citep{Song2022,Song2024}) is still
valid.

Section~\ref{subsec:Single-degree-of-freedom-systems} performs a
convergence study on single-degree-of-freedom systems; first on a
linear spring subject to external excitation, then on a highly nonlinear
simple oscillating pendulum. The linear problem demonstrates that
the new algorithm derived in the preceding sections replicates numerical
results of the previous algorithms (\citep{Song2022,Song2024}), though
with additional capability to capture accelerations at each time step
to very high accuracy without additional equation solution. The pendulum
problem is the first time the accuracy of this family of high-order
time integration schemes is presented on nonlinear problems, and it
is shown that the same convergence rates for nonlinear problems as
linear problems are achieved. Section~\ref{subsec:A-three-degree-of-freedom-model}
assesses a three-degree-of-freedom model problem characterized by
a high stiffness ratio, where results are again compared to those
from our previous research, confirming the consistency and reliability
of the proposed algorithm. A new nonlinear three-degree-of-freedom
problem is also presented demonstrating the controllable numerical
dissipation properties and accuracy of the model over long time periods.

Comparison of the results from the proposed schemes against those
obtained using commercial software are carried out in Sections~\ref{subsec:One-dimensional-wave-propagation-rod}
and \ref{subsec:Two-dimensional-wave-propagation-Lamb}, emphasizing
their superior performance in suppressing high-frequency oscillations
and the high-accuracy in complex wave propagation models.

The proposed algorithm for calculating accelerations (Section~\ref{sec:Calculation-of-acceleration})
can be evaluated by comparing accelerations at each time step with
those obtained by directly solving the equation of motion (Eq.~(\ref{eq:EqofMotionNonlinear}))
using the consistent mass matrix with the displacement and velocity
solution. This is performed for the examples in Sections~\ref{subsec:One-dimensional-wave-propagation-rod}
and \ref{subsec:Two-dimensional-wave-propagation-Lamb} and several
other examples not reported in this paper. The relative difference
in $L$-infinity norms of all the tested cases are less than $10^{-8}$.
It is concluded that the difference is negligible.

\subsection{Single-degree-of-freedom (SDOF) systems}\label{subsec:Single-degree-of-freedom-systems}

\subsubsection{Linear error analysis}

To evaluate the accuracy of the proposed algorithm, a convergence
study of the single-degree-of-freedom (SDOF) system used in our previous
works \citep{Song,Song2024} is performed. The parameters in \citep{Kim2019}
are adopted: natural frequency $\omega=2\pi$, damping ratio $\xi=0$,
initial displacement $u_{0}=2.0\,\unit{m}$, initial velocity $\dot{u}_{0}=\pi/3\,\unit{m/s}$.
The external excitation is given as a harmonic function:
\begin{equation}
f_{1}(t)=10\cos\left(\frac{2\sqrt{5}}{5}t\right)+70\sin\left(2\sqrt{10}t\right)\,,\label{eq:f1}
\end{equation}
and the analysis is performed for $t_{\mathrm{sim}}=10\,\mathrm{s}$.
To assess the accuracy of the time integration schemes, errors are
calculated based on the $L_{2}$-norm ($\epsilon_{L_{2}}$) 
\begin{equation}
\epsilon_{L_{2}}=\cfrac{\int\limits_{t=0}^{t_{\mathrm{sim}}}(\ddot{u}_{\mathrm{exact}}(t)-\ddot{u}_{\mathrm{numerical}}(t))^{2}\;\mathrm{d}t}{\int\limits_{t=0}^{t_{\mathrm{sim}}}\ddot{u}_{\mathrm{exact}}(t)^{2}\;\mathrm{d}t}\times100[\%]\,.\label{eq:ErrorL2}
\end{equation}

For the sake of conciseness, only the error plots and convergence
rates for accelerations are included here. The errors for displacements
and velocities were presented in \citep{Song2024} and show identical
qualitative behaviour.

\begin{figure}
\centering \subfloat{\includegraphics[width=0.6\textwidth]{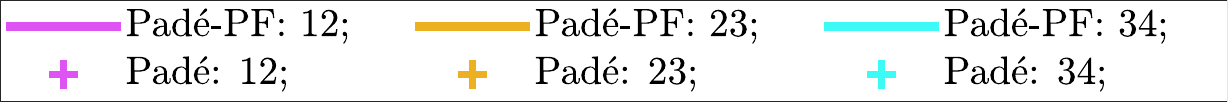}}\medskip{}
\\
 \setcounter{subfigure}{0}

\includegraphics[width=0.5\textwidth]{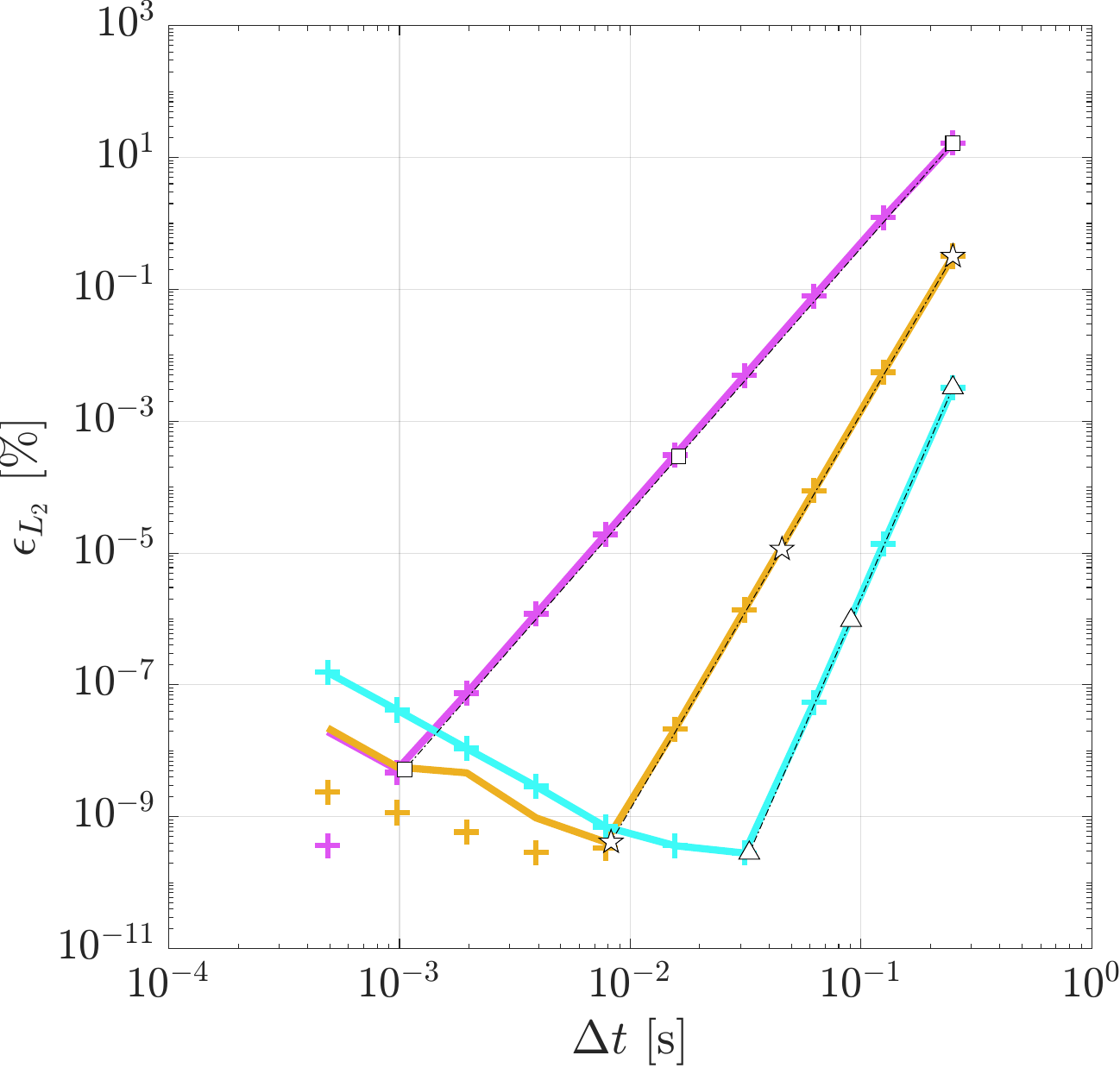}\hfill{}\includegraphics[width=0.5\textwidth]{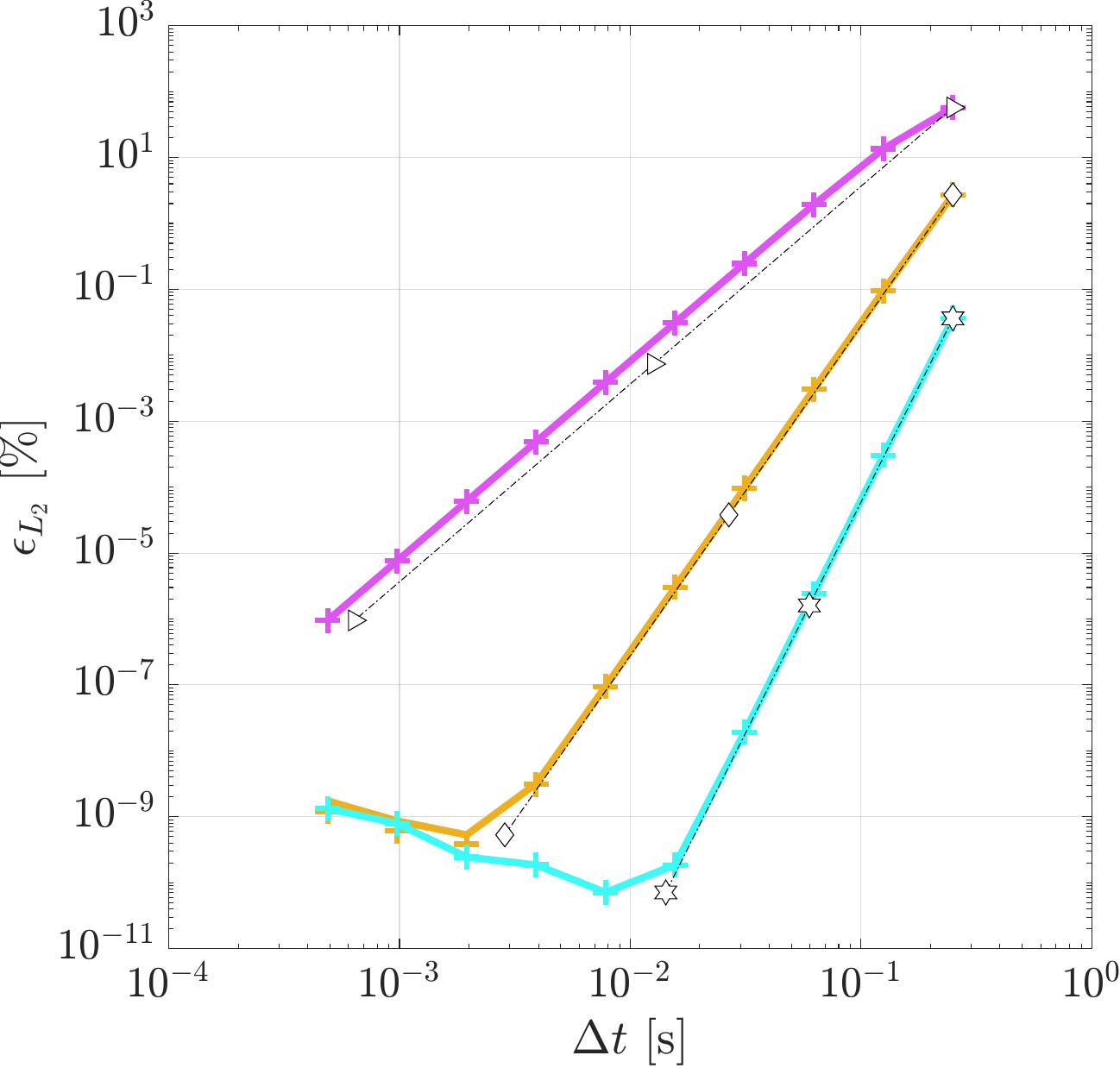}

\caption{Acceleration error in the $L_{2}$-norm for the single degree of freedom
system obtained by the algorithm for mixed-order Padé-based schemes
in \citep{Song} with distinct roots. Left: $\rho_{\infty}=1$ ; Right:
$\rho_{\infty}=0$. The dash-dotted lines indicate the optimal rates
of convergence corresponding to slopes of 3 (right-pointing triangle),
4 (square), 5 (diamond), 6 (pentagram), 7 (hexagram), and 8 (upward-pointing
triangle), respectively. }\label{fig:SDOF_Error_multiple}
\end{figure}

\begin{figure}
\centering \subfloat{\includegraphics[width=0.6\textwidth]{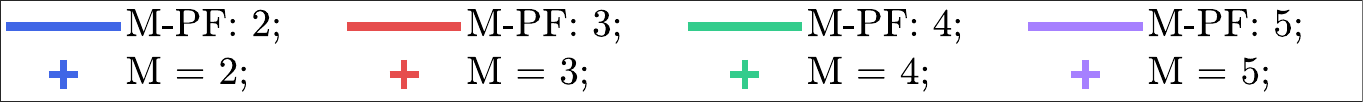}}\medskip{}
\\
 \setcounter{subfigure}{0}

\includegraphics[width=0.5\textwidth]{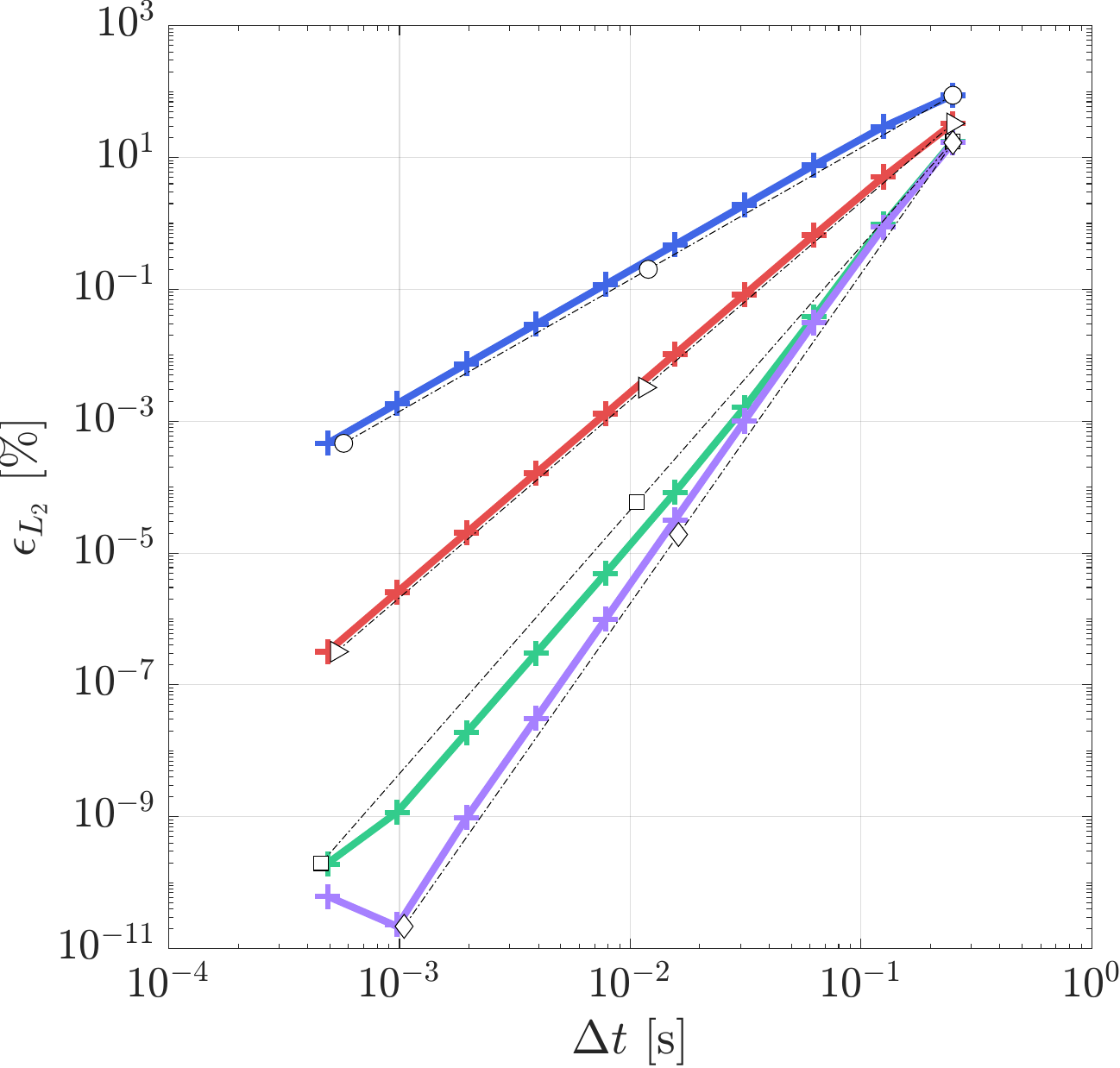}\hfill{}\includegraphics[width=0.5\textwidth]{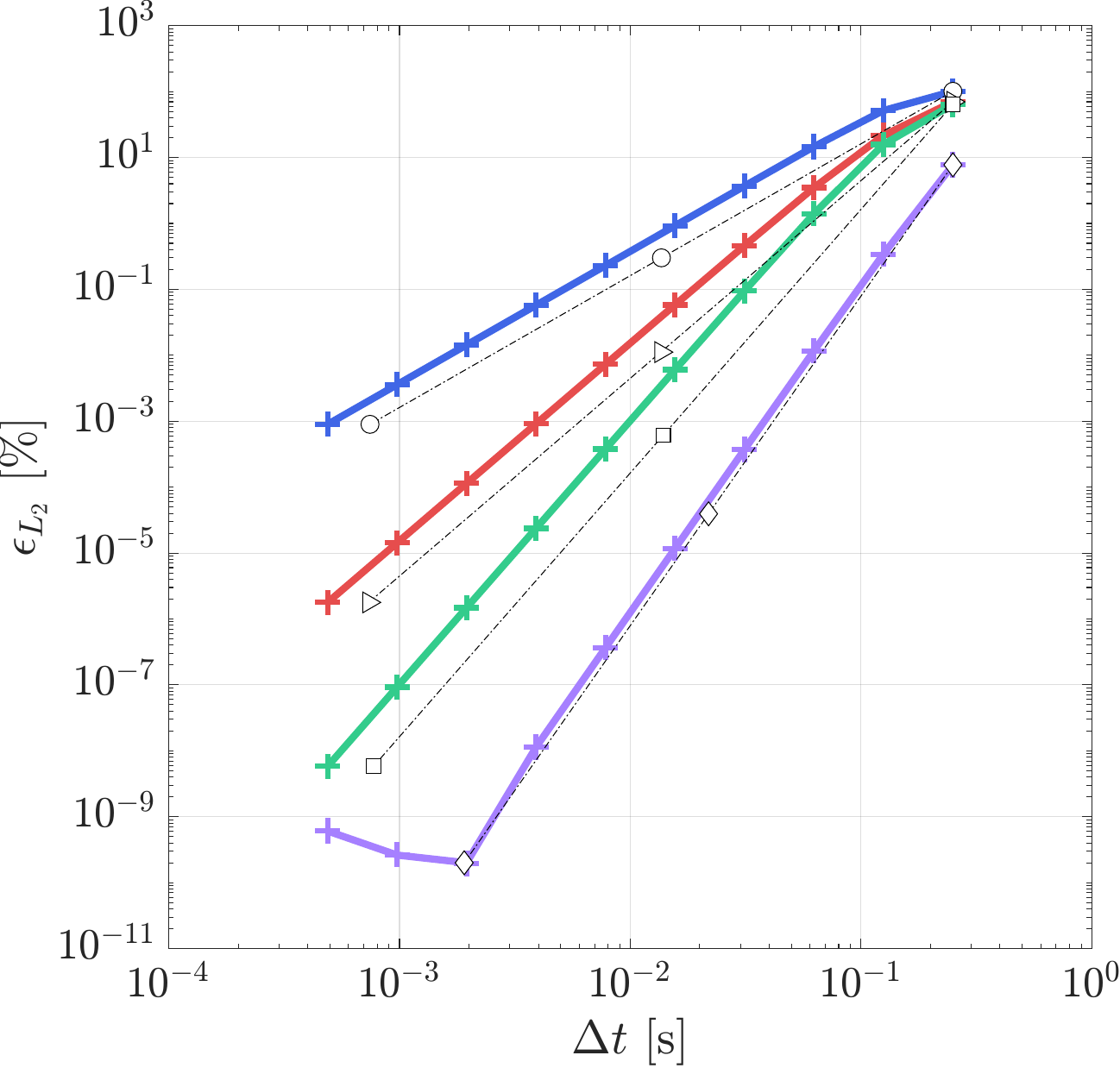}\caption{Acceleration error in the $L_{2}$-norm for the single degree of freedom
system obtained by the algorithm for $M-$schemes in \citep{Song2024}
with single multiple roots. Left: $\rho_{\infty}=1$ ; Right: $\rho_{\infty}=0$.
The dash-dotted lines indicate the optimal rates of convergence corresponding
to slopes of 2 (circle), 3 (right-pointing triangle), 4 (square) and
5 (diamond), respectively. }\label{fig:SDOF_Error_single}
\end{figure}

The convergence behaviors are shown at two typical values of the user-specified
parameters: $\rho_{\infty}=1$ (no numerical damping) and $0$ (maximum
numerical damping). Results corresponding to intermediate values of
$\rho_{\infty}$, which lie between these two extremes, closely align
with those seen at $\rho_{\infty}=0$ and are omitted for conciseness.
The error plots for the cases of distinct roots and single multiple
roots are depicted with solid lines in Fig.~\ref{fig:SDOF_Error_multiple}
and Fig.~\ref{fig:SDOF_Error_single}, respectively. Results obtained
using the composite $M-$schemes in~\citep{Song2024} and using corresponding
mixed-order Padé-based schemes in \citep{Song} are indicated by the
plus markers in corresponding colors in those two figures respectively.
It is observed that the present algorithms based on partial fractions
are numerically identical to the algorithm in \citep{Song2024} for
the composite $M-$schemes with single multiple roots, and to the
algorithm in \citep{Song} for the mixed-order Padé-based schemes
with distinct roots.

\subsubsection{Nonlinear error analysis}

\begin{figure}
\begin{centering}
\includegraphics[width=0.3\textwidth]{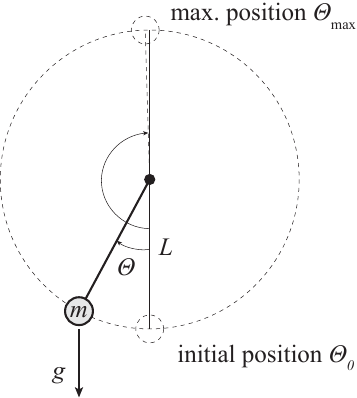}
\par\end{centering}
\caption{Nonlinear simple pendulum}\label{fig:pendulum}
\end{figure}

\begin{figure}
\centering \subfloat{\includegraphics[width=0.55\textwidth]{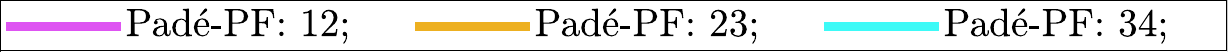}}\medskip{}
\\
 \setcounter{subfigure}{0}

\includegraphics[width=0.5\textwidth]{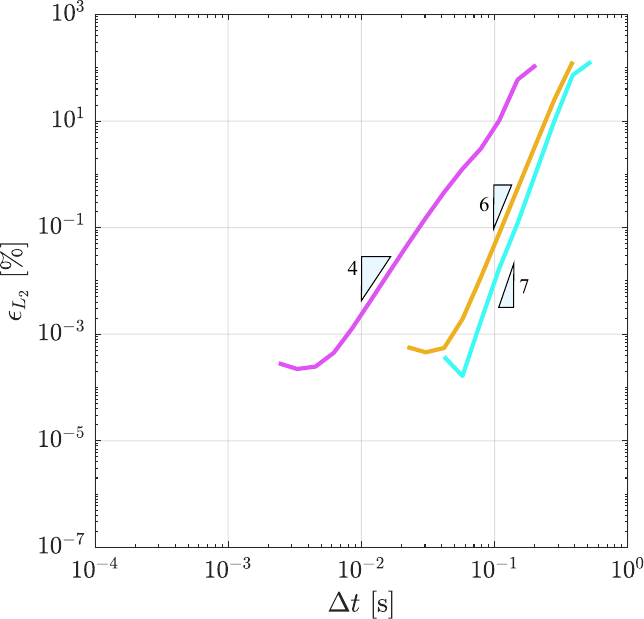}\hfill{}\includegraphics[width=0.5\textwidth]{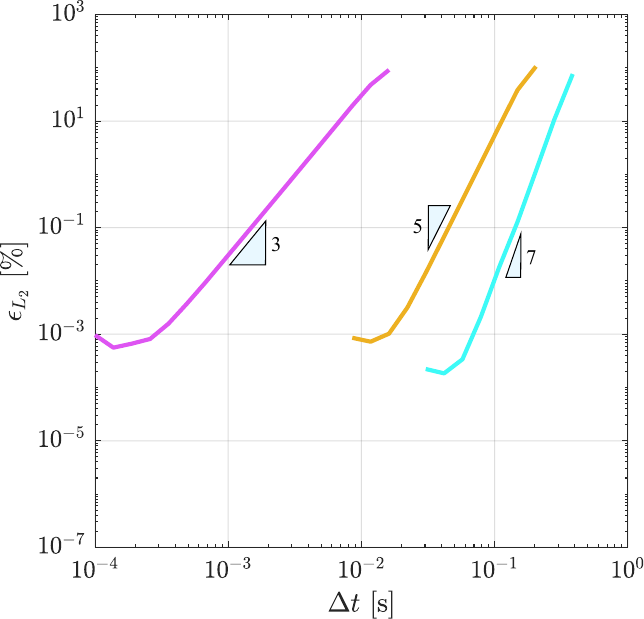}

\caption{Acceleration error in the $L_{2}$-norm for the nonlinear oscillating
simple pendulum obtained by the algorithm for mixed-order Padé-based
schemes with distinct roots. Left: $\rho_{\infty}=1$ ; Right: $\rho_{\infty}=0$.
}\label{fig:pendulum_Error_PADE}
\end{figure}

\begin{figure}
\centering \subfloat{\includegraphics[width=0.8\textwidth]{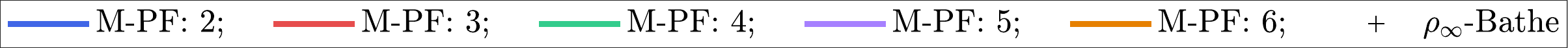}}\medskip{}
\\
 \setcounter{subfigure}{0}

\includegraphics[width=0.5\textwidth]{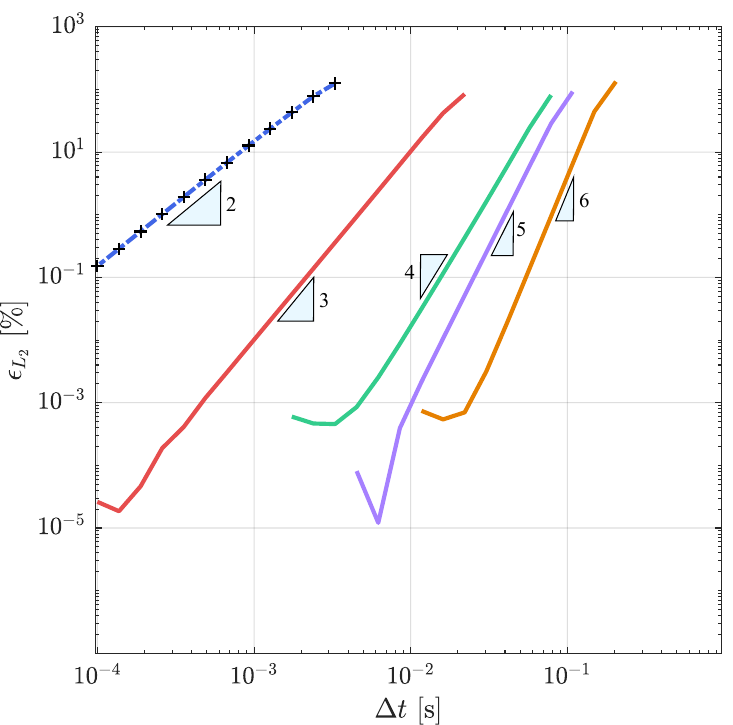}\hfill{}\includegraphics[width=0.5\textwidth]{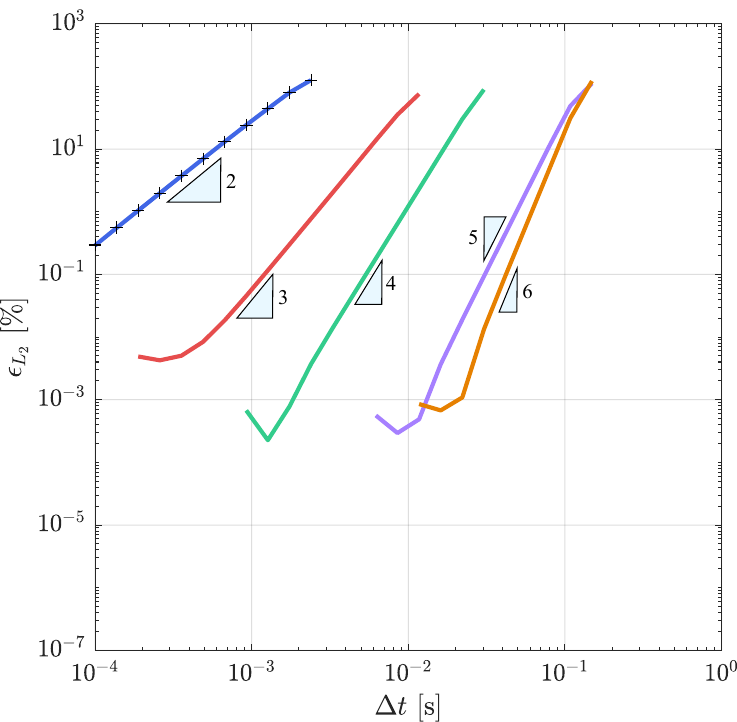}\caption{Acceleration error in the $L_{2}$-norm for the nonlinear simple pendulum
obtained by the algorithm for $M-$schemes with single multiple roots.
Left: $\rho_{\infty}=1$ ; Right: $\rho_{\infty}=0$. }\label{fig:Pendulum_Error_M}
\end{figure}

To demonstrate the performance of the proposed algorithms on problems
with strong nonlinearity, the nonlinear (SDOF) oscillating simple
pendulum shown in Fig.~\ref{fig:pendulum} is considered here. The
governing equation of motion is given as
\[
\ddot{\theta}+\omega^{2}\sin\theta=0,\thinspace\thinspace\thinspace\theta_{0}=0\thinspace\mathrm{rad},\thinspace\thinspace\dot{\theta_{0}=1.999999238456499\thinspace\mathrm{rad/s}}\thinspace,
\]
with $\omega=\sqrt{g/L}$ and the dimensionless case of $L=1,\thinspace g=1$
is used. The initial velocity is chosen such that the pendulum possesses
the energy to almost reach its maximum position ($\theta_{\mathrm{max}}=\pm179.9^{\circ}$)
though not complete full rotations~\citep{Kim2018}. 
The period of oscillation is approximately $T=33.712$ seconds, and
the pendulum is simulated for two full periods. The errors shown in
Fig.~\ref{fig:pendulum_Error_PADE} and \ref{fig:Pendulum_Error_M},
for distinct root models and single root models, respectively, are
based on the $L_{2}$-norm defined earlier in Eq.~(\ref{eq:ErrorL2}).
The reference solution follows the exact solution provided in \citep{Kim2018}.
Once more for conciseness, only the error plots and convergence rates
for accelerations are included here; the errors for displacements
and velocities are of similar performance. It can be observed that
the theoretical convergence rates are able to be met for this nonlinear
problem. In Fig.~\ref{fig:pendulum_Error_PADE}a, the theoretical
convergence rate is $2M$, however, there is an upper bound of 7 due
to the interpolation scheme in use for state-space variables. For
Fig.~\ref{fig:pendulum_Error_PADE}, the theoretical error rate is
$2M-1$. For the $M$-schemes in Fig.~\ref{fig:Pendulum_Error_M},
the theoretical error rate is $M$, and is achieved for all models.
Numerical results using the $\rho_{\infty}$-Bathe model \citep{Malakiyeh2019}
overlap on the plots with the second-order model in this paper.

\begin{figure}
\begin{centering}
\includegraphics[width=0.45\textwidth]{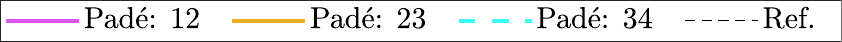}\smallskip{}
\par\end{centering}
\includegraphics[width=1\textwidth]{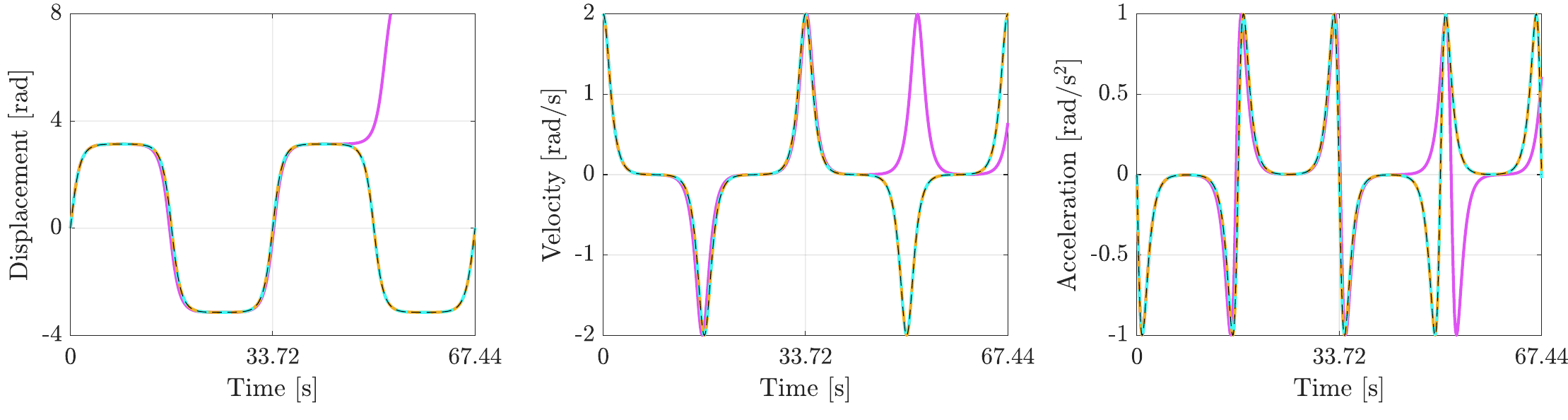}

\caption{Time-history responses of the nonlinear oscillating pendulum using
mixed-order Padé-based schemes, with $\Delta t=0.1686$ and $\rho_{\infty}=1$.
Note that models relating to the (2,3) and (3,4) Padé expansions essentially
overlap.}\label{fig:pendulum_Pade}
\end{figure}

\begin{figure}
\begin{centering}
\includegraphics[width=0.3\textwidth]{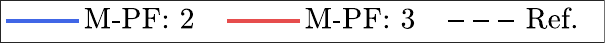}\smallskip{}
\par\end{centering}
\includegraphics[width=1\textwidth]{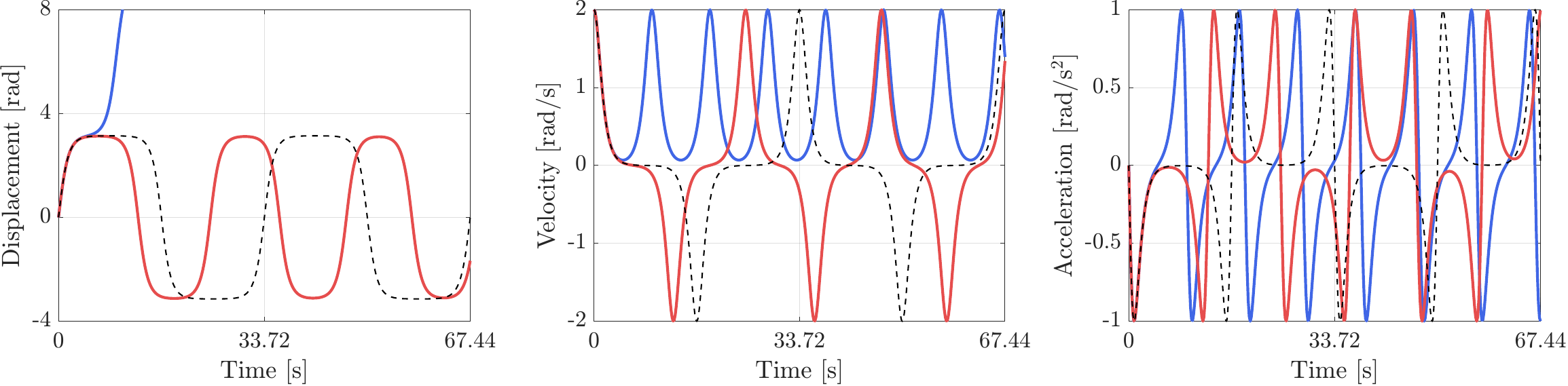}

\caption{Time-history responses of the nonlinear oscillating pendulum using
second- and third-order $M$-based schemes, with $\Delta t=0.1686$
and $\rho_{\infty}=1$. }\label{fig:pendulum_M-low}
\end{figure}

\begin{figure}
\begin{centering}
\includegraphics[width=0.45\textwidth]{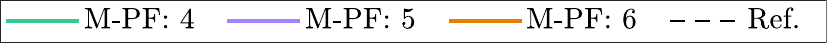}\smallskip{}
\par\end{centering}
\includegraphics[width=1\textwidth]{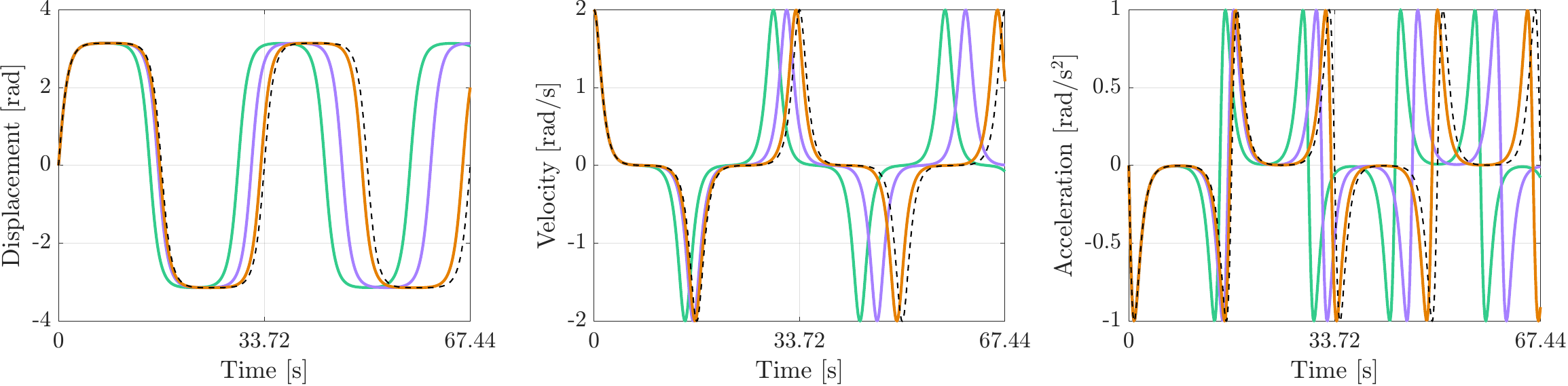}

\caption{Time-history responses of the nonlinear oscillating pendulum using
higher-order $M$-based schemes, with $\Delta t=0.1686$ and $\rho_{\infty}=1$.
}\label{fig:pendulum_M}
\end{figure}

Figure~\ref{fig:pendulum_Pade} presents the time-history for models
using the distinct root algorithm presented in Table~\ref{tab:Time-integration-solution_distinct},
with no numerical dissipation introduced ($\rho_{\infty}=1$). A large
time step of $\Delta t=T/200=0.1686$ seconds is used to illustrate
the difference in the responses obtained at various orders. As can
be seen, during the second period of oscillation the (1,2) model erroneously
completes a full revolution, whilst the (2,3) and (3,4) models are
near indistinguishable from the reference solution. Figures~\ref{fig:pendulum_M-low}
and \ref{fig:pendulum_M} show the time history for single root algorithms
presented in Table~\ref{tab:Time-integration-solution_multiple}.
The second-order model  completes full revolutions, whilst the third-order
model presents large period errors. For higher-order schemes, Fig~\ref{fig:pendulum_M}
demonstrates how the period error reduces as the number of sub-steps
increases, converging on the reference solution. 

Once more, it is observed from Figs.~\ref{fig:pendulum_Error_PADE}
and \ref{fig:Pendulum_Error_M} that:
\begin{enumerate}
\item All schemes provide identical qualitative convergence as for linear
problems.
\item All schemes converge at their optimal rates, up to a maximum of seven
due to the Hermite interpolation polynomial used in the nonlinear
solution algorithm.
\item For a given scheme, the error increases with decreasing $\rho_{\infty}$,
i.e., increasing the amount of numerical damping. This difference
is significant here due to the high nonlinearity of the problem.
\item The order of accuracy for schemes with single multiple roots is not
affected by the value of $\rho_{\infty}$, while the order of accuracy
for distinct roots cases will decrease by one if numerical damping
is introduced.
\item High-order models achieve high levels of accuracy ($<10^{-3}\%$)
using time steps several orders of magnitude larger than second- or
third-order models.
\end{enumerate}

\subsection{Three-degree-of-freedom model problems}\label{subsec:A-three-degree-of-freedom-model}

\begin{figure}
\begin{centering}
\includegraphics[width=0.6\textwidth]{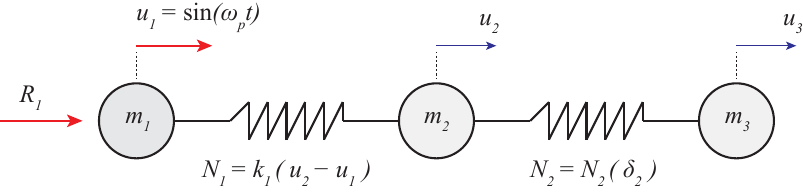}
\par\end{centering}
\caption{A three-degree-of-freedom model problem for linear and nonlinear analyses
}\label{fig:A-linear-three-degree-of-freedom}
\end{figure}

This section details the analysis of a linear and a nonlinear three-degree-of-freedom
(3DOF) problem. The 3DOF problem is shown in Fig.~\ref{fig:A-linear-three-degree-of-freedom}.
Each mass is taken to be $m_{1}=0$, $m_{2}=m_{3}=1$, and the system
is initially at rest: $u_{2}(0)=\dot{u}_{2}(0)=u_{3}(0)=\dot{u}_{3}(0)=0$.
The displacement of mass 1 is prescribed as $u_{1}=\sin(\omega_{p}t$)
with $\omega_{p}=1.2$, corresponding to a period of vibration of
$T_{p}=5.236$. Spring 1 is taken to be a stiff linear spring with
coefficient $k_{1}=10^{7}$. The behaviour of Spring 2 is governed
by the internal force function $N_{2}=N_{2}\left(\delta_{2}\right)$,
a function of the elongation $\delta_{2}=u_{3}-u_{2}$. Following
static condensation of the known displacement $u_{1}$, the equation
of motion of the system as per Eq.~(\ref{eq:EqofMotionNonlinear})
is expressed as:
\begin{equation}
\left[\begin{array}{cc}
m_{2} & 0\\
0 & m_{3}
\end{array}\right]\left\{ \begin{array}{c}
\ddot{u}_{2}\\
\ddot{u}_{3}
\end{array}\right\} +\left\{ \begin{array}{c}
k_{1}u_{2}-N_{2}\\
N_{2}
\end{array}\right\} =\left\{ \begin{array}{c}
k_{1}u_{1}\\
0
\end{array}\right\} ,
\end{equation}
which can be written in the form of Eq.~(\ref{eq:linear=000020eq=000020of=000020motion})
by taking as the mass, tangent damping and tangent stiffness, respectively:
\begin{equation}
\mathbf{M}=\left[\begin{array}{cc}
m_{2} & 0\\
0 & m_{3}
\end{array}\right],\thinspace\thinspace\mathbf{C}_{n-1}=\mathbf{0},\thinspace\thinspace\mathbf{K}_{n-1}=\left[\begin{array}{cc}
k_{1}+\frac{{\rm d}N_{2}}{{\rm d}\delta_{2}} & -\frac{{\rm d}N_{2}}{{\rm d}\delta_{2}}\\
-\frac{{\rm d}N_{2}}{{\rm d}\delta_{2}} & \frac{{\rm d}N_{2}}{{\rm d}\delta_{2}}
\end{array}\right]\thinspace,\label{eq:MCK_3DOF}
\end{equation}
where the derivative ${\rm d}N_{2}/{\rm d}\delta_{2}$ is evaluated
at $t_{n-1}$ in $\mathbf{K}_{n-1}$. The nonlinear force vector $\mathbf{f}(t)$
is given by
\begin{equation}
\mathbf{f}=\left\{ \begin{array}{c}
\left(k_{1}u_{1}+N_{2}\right)-\frac{{\rm d}N_{2}}{{\rm d}\delta_{2}}\delta_{2}\\
-N_{2}+\frac{{\rm d}N_{2}}{{\rm d}\delta_{2}}\delta_{2}
\end{array}\right\} \thinspace.\label{eq:force_3DOF}
\end{equation}

The reaction force at $m_{1}$ is $R_{1}=m_{1}\ddot{u}_{1}+k_{1}(u_{1}-u_{2})$.

\subsubsection{Linear analysis}

A linear analysis is performed to compare the proposed algorithm based
on partial fraction decomposition with the original algorithms in
\citep{Song,Song2024}. Here, the linear 3DOF problem studied in~\citep{Bathe2012,Noh2018}
is considered, with $N_{2}=k_{2}\delta_{2}$, and $k_{2}=1$. As expected,
$\mathbf{K}_{n-1}$ in Eq.~(\ref{eq:MCK_3DOF}) is constant-valued,
and the force-vector from Eq.~(\ref{eq:force_3DOF}) is independent
of unknown displacements and velocities. To suppress the high-frequency
oscillations, $\rho_{\infty}=0$ is chosen. The time step size $\Delta t=0.14$
is selected as was used in \citep{Choi2022}. The analysis is performed
over an extended period of $t=5000$, approximately equivalent to
$967\thinspace T_{p}$, to ensure a comprehensive analysis.

\begin{figure}
\begin{centering}
\includegraphics[width=0.6\textwidth]{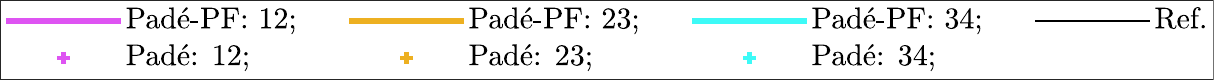}
\par\end{centering}
\includegraphics[width=1\textwidth]{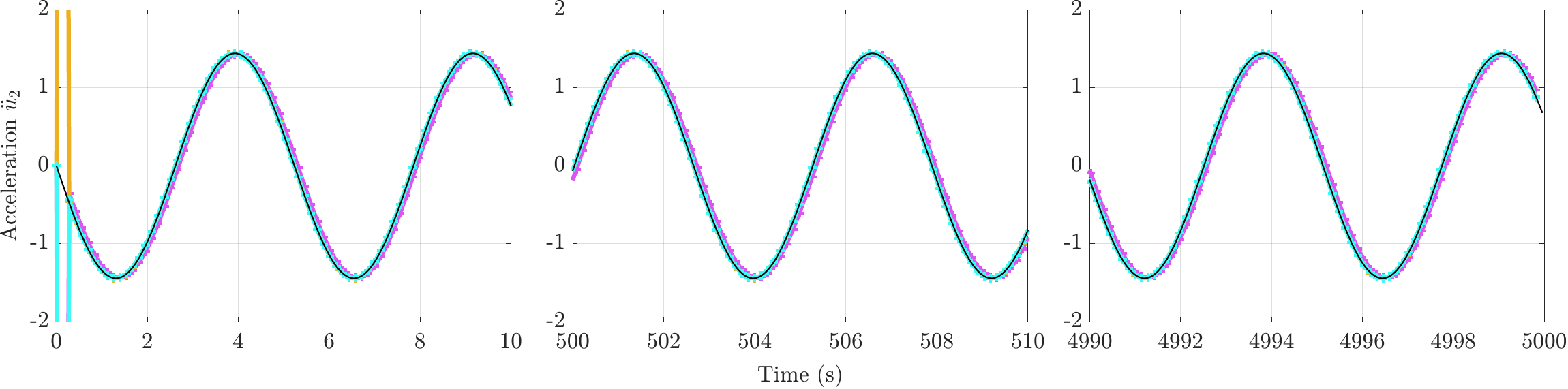}

\includegraphics[width=1\textwidth]{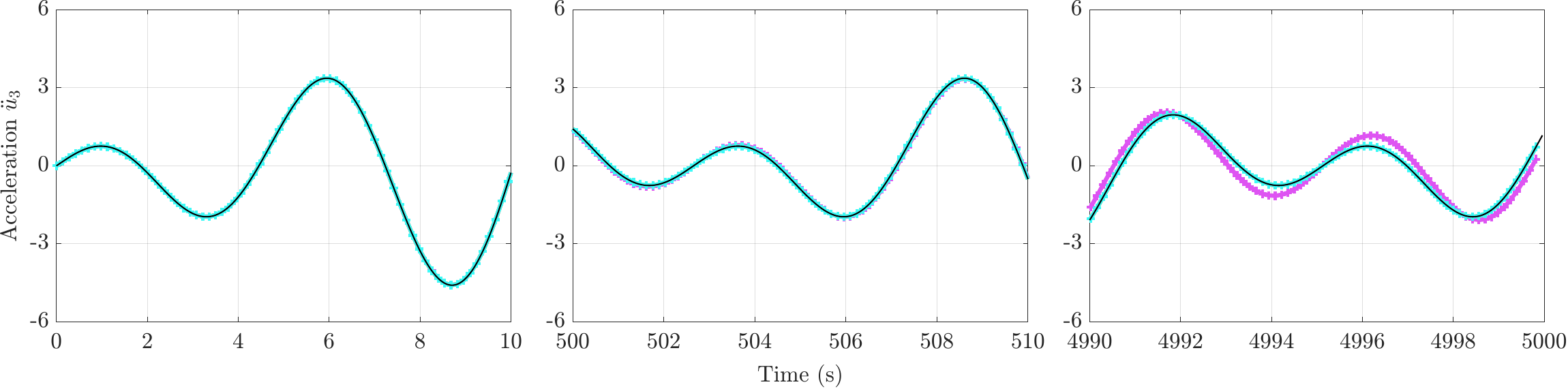}

\includegraphics[width=1\textwidth]{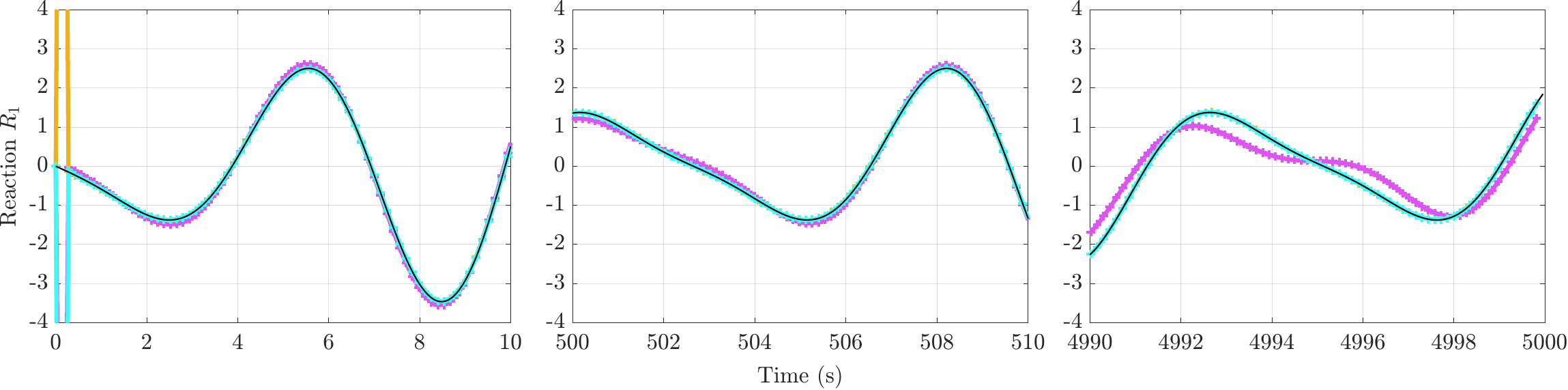}

\caption{Acceleration responses of $m_{2}$ (top row), acceleration responses
of $m_{3}$ (middle row), and reaction force response (bottom row)
for distinct root cases during $0\protect\leq t\protect\leq10$ (left
column), $500\protect\leq t\protect\leq510$ (middle column) and $4990\protect\leq t\protect\leq5000$
(right column). }\label{fig:3dofs_Distinct}
\end{figure}

In the time-stepping schemes utilizing distinct roots, the acceleration
responses of $m_{2}$, $m_{3}$, and the reaction force $R_{1}(t)$
are illustrated in Fig.~\ref{fig:3dofs_Distinct}. The solid lines
represent results obtained with the proposed algorithm using partial
fractions, while the plus markers with the same color denote results
using the previous algorithm in \citep{Song}. The three columns of
each figure show the responses at three different time intervals:
$0\leq t\leq10$ (left column), $500\leq t\leq510$ (middle column),
and $4990\leq t\leq5000$ (right column). It is observed that both
algorithms are numerically identical. Although the proposed method
at order (1,2) exhibits a slight difference from the reference solution
as time increases, the results converge to the reference solution
as the order increases. At higher orders, specifically (2,3) and (3,4),
the results overlap with the reference solution throughout the whole
duration. 

\begin{figure}
\begin{centering}
\includegraphics[width=0.6\textwidth]{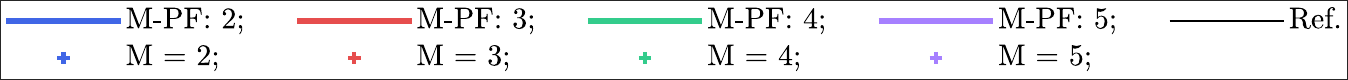}
\par\end{centering}
\includegraphics[width=1\textwidth]{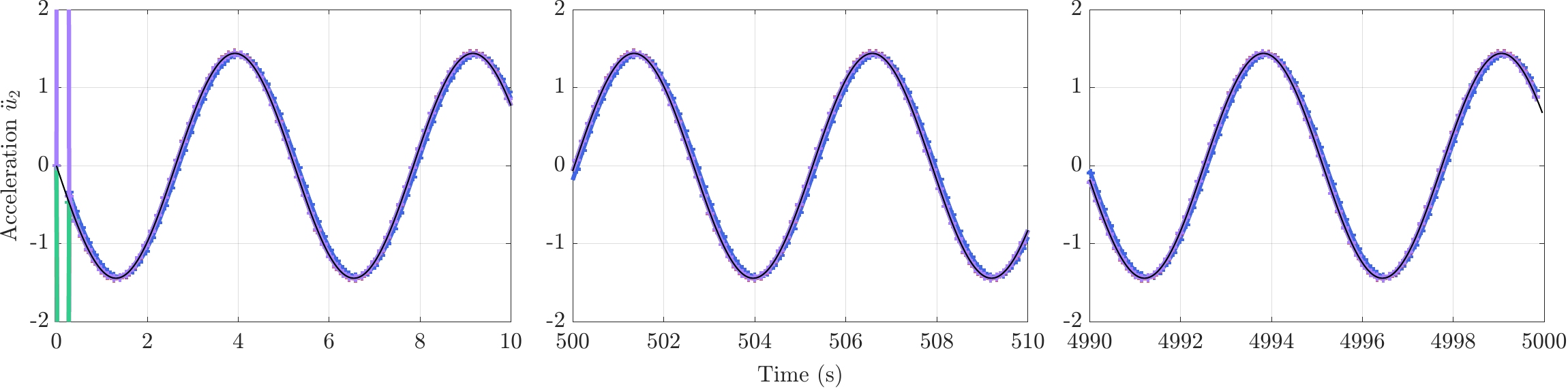}

\includegraphics[width=1\textwidth]{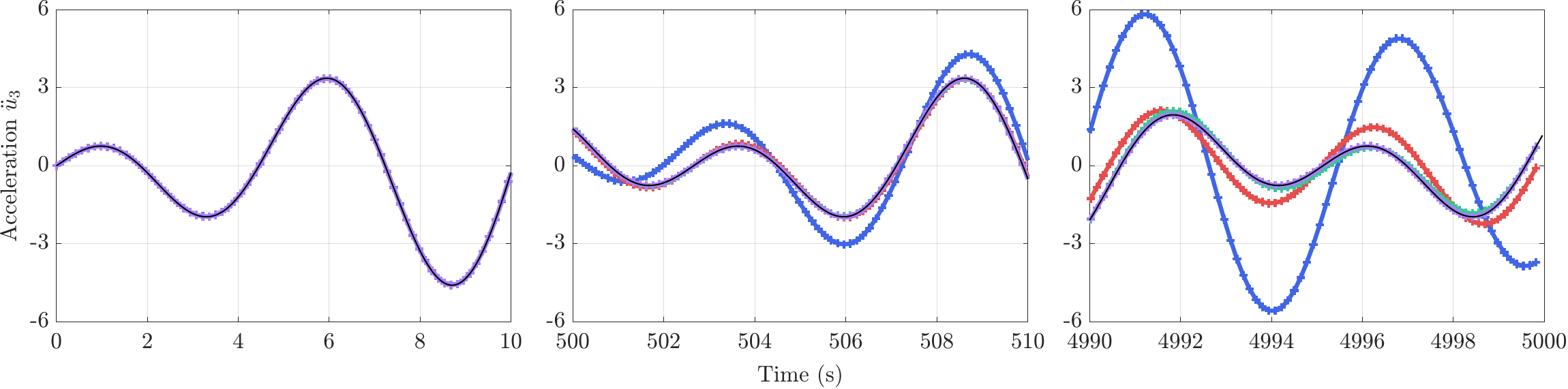}

\includegraphics[width=1\textwidth]{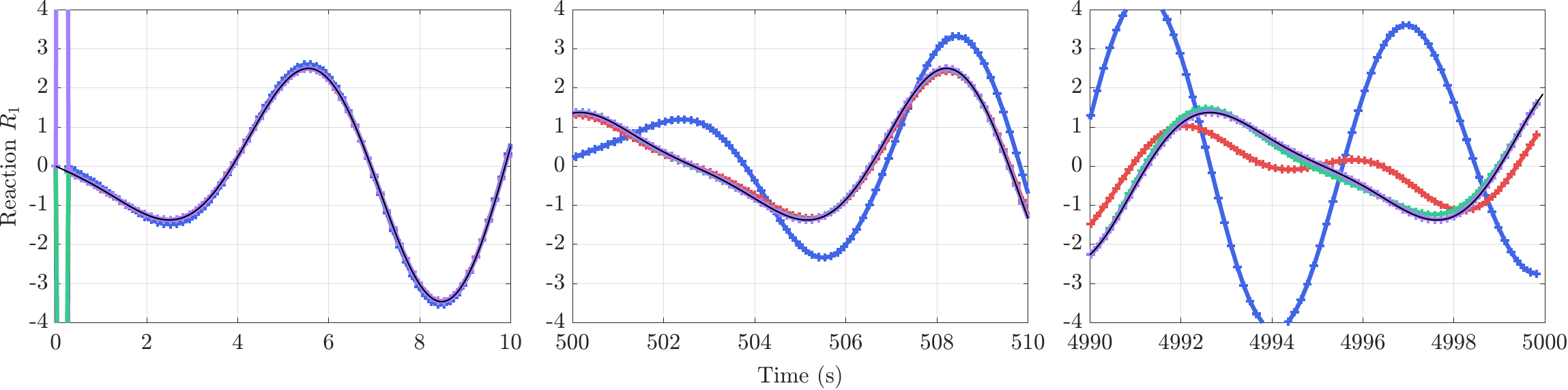}

\caption{Acceleration responses of $m_{2}$ (top row), acceleration responses
of $m_{3}$ (middle row), and reaction force response (bottom row)
for the $M-$schemes with single multiple root during $0\protect\leq t\protect\leq10$
(left column), $500\protect\leq t\protect\leq510$ (middle column)
and $4990\protect\leq t\protect\leq5000$ (right column). }\label{fig:3dofs_Single}
\end{figure}

For the $M-$schemes in \citep{Song2024} with single multiple roots,
the velocity and acceleration responses of $m_{2}$, $m_{3}$, and
the reaction force $R_{1}(t)$ obtained with the present algorithm
are plotted in Fig.~\ref{fig:3dofs_Single}. The solid lines represent
results obtained with the proposed algorithm using partial fractions,
and results from the previous algorithm in \citep{Song2024} are denoted
by plus markers in corresponding colors. It is observed that the present
algorithm and the algorithm in \citep{Song2024} yield numerically
identical results. While all schemes exhibit satisfactory performance
in the responses of $m_{2}$, significant deviations from the reference
solutions are observed in the responses of $m_{3}$ and reaction force
over time with the second-order scheme. Increasing the order leads
to decreasing errors, and by fifth-order ($M=5$) the response is
indistinguishable from the reference solution.

\subsubsection{Nonlinear analysis}

For the nonlinear case, Spring 2 is modelled by a nonlinear stiffening
spring with internal force
\begin{equation}
N_{2}=k_{2}\sinh\delta_{2}\thinspace,
\end{equation}
and $k_{2}=1$ is taken. The force-vector of Eq.~(\ref{eq:force_3DOF})
is nonlinear with respect to displacements. 

\begin{figure}
\begin{centering}
\includegraphics[width=0.35\textwidth]{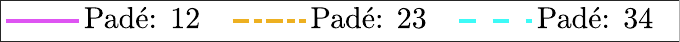}\smallskip{}
\par\end{centering}
\includegraphics[width=1\textwidth]{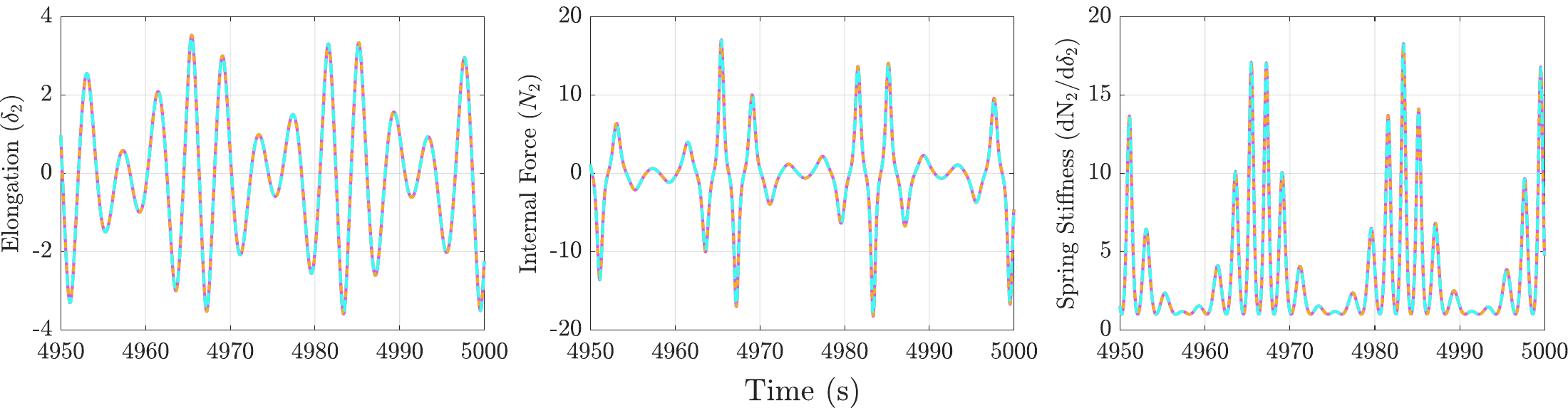}

\caption{Response of nonlinear spring connecting masses $m_{2}$ and $m_{3}$
for the Padé-based schemes during $4950\protect\leq t\protect\leq5000$
}\label{fig:3dofs_Pade-nonlinear-stiffness}
\end{figure}

\begin{figure}
\begin{centering}
\includegraphics[width=0.35\textwidth]{pade-legend}\smallskip{}
\par\end{centering}
\includegraphics[width=1\textwidth]{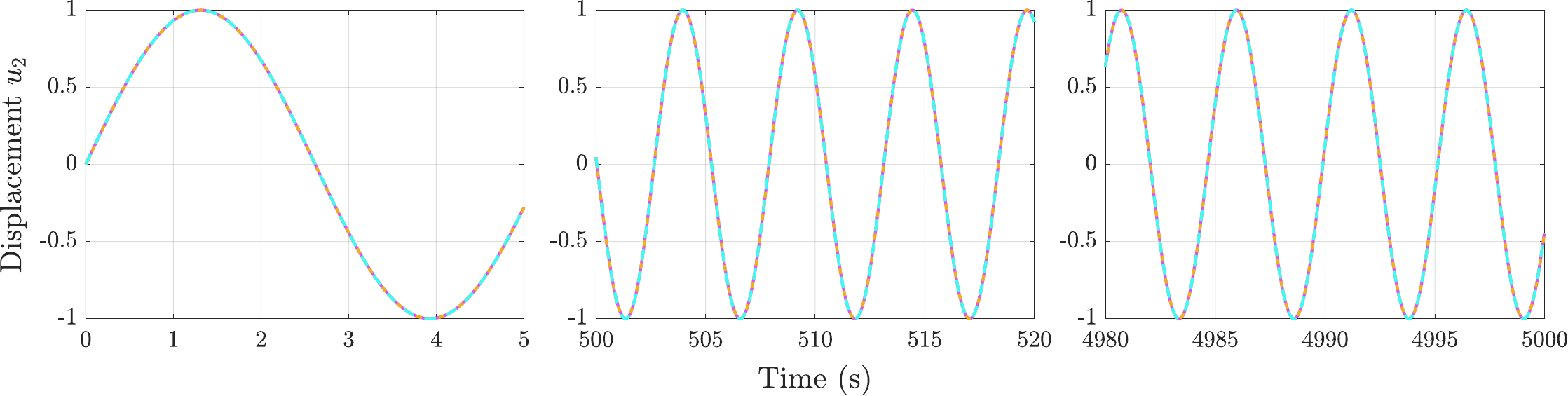}

\includegraphics[width=1\textwidth]{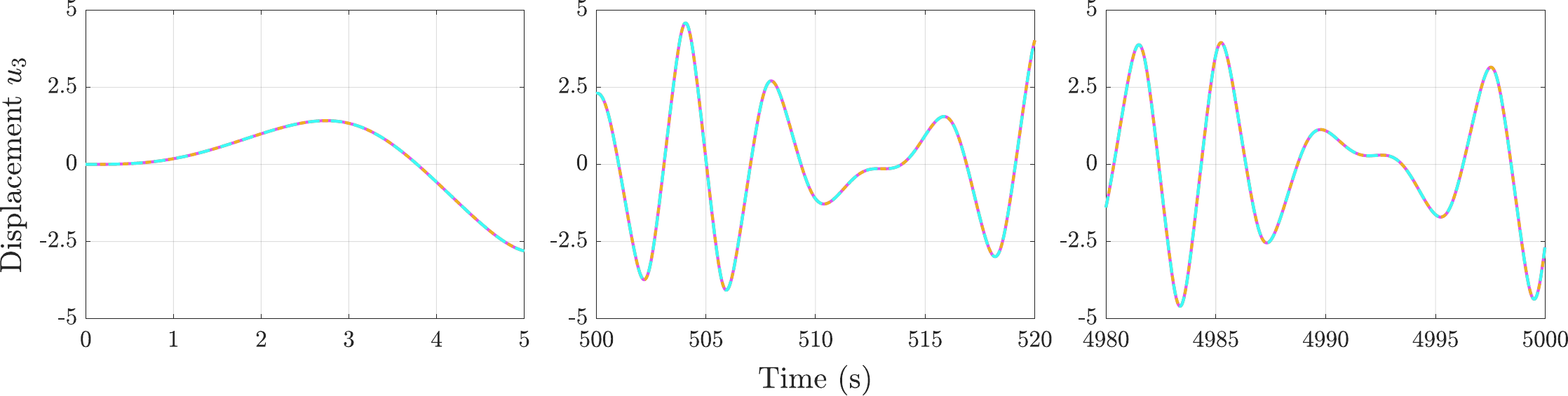}

\caption{Displacement responses of $m_{2}$ (top row), and $m_{3}$ (bottom
row), for the Padé-based schemes during $0\protect\leq t\protect\leq5$
(left column), $500\protect\leq t\protect\leq520$ (middle column)
and $4980\protect\leq t\protect\leq5000$ (right column). }\label{fig:3dofs_Pade-nonlinear-dsp}
\end{figure}

\begin{figure}
\begin{centering}
\includegraphics[width=0.35\textwidth]{pade-legend}\smallskip{}
\par\end{centering}
\includegraphics[width=1\textwidth]{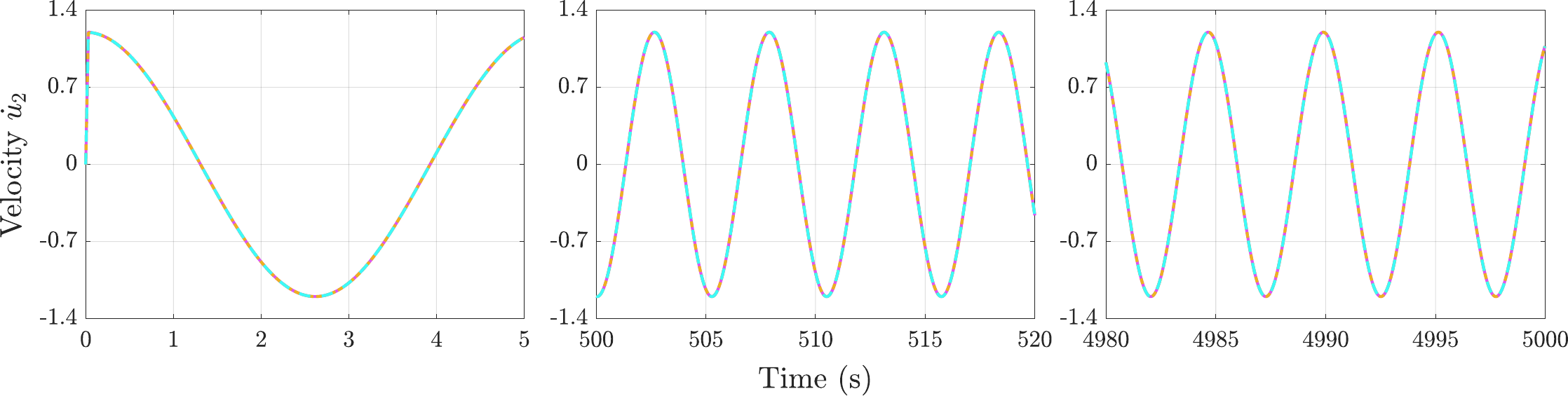}

\includegraphics[width=1\textwidth]{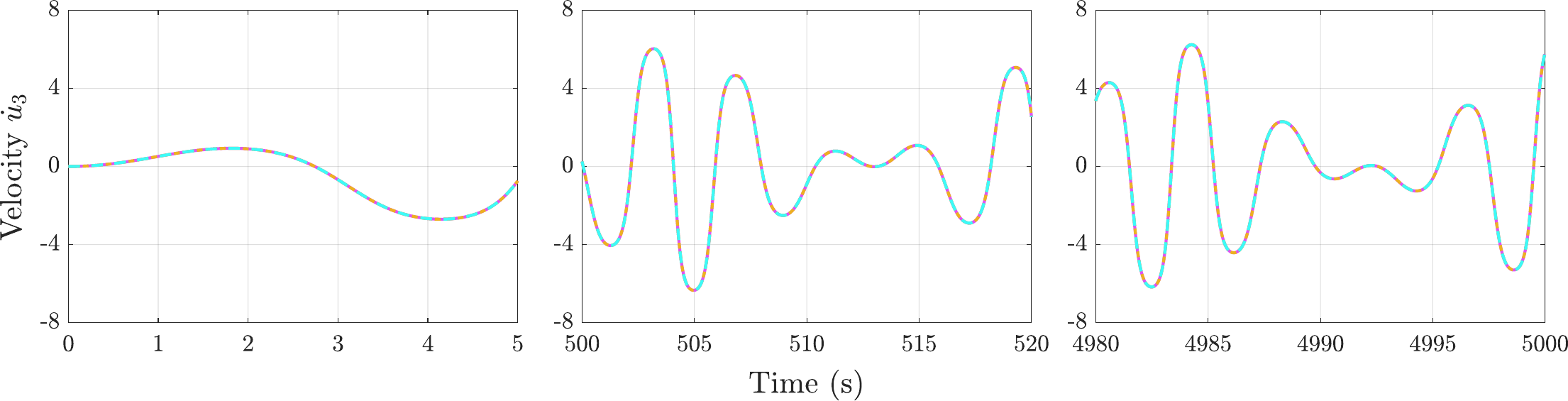}

\caption{Velocity responses of $m_{2}$ (top row) and $m_{3}$ (bottom row),
for the Padé-based schemes during $0\protect\leq t\protect\leq5$
(left column), $500\protect\leq t\protect\leq520$ (middle column)
and $4980\protect\leq t\protect\leq5000$ (right column). }\label{fig:3dofs_Pade-nonlinear-vel}
\end{figure}

\begin{figure}
\begin{centering}
\includegraphics[width=0.35\textwidth]{pade-legend}\smallskip{}
\par\end{centering}
\includegraphics[width=1\textwidth]{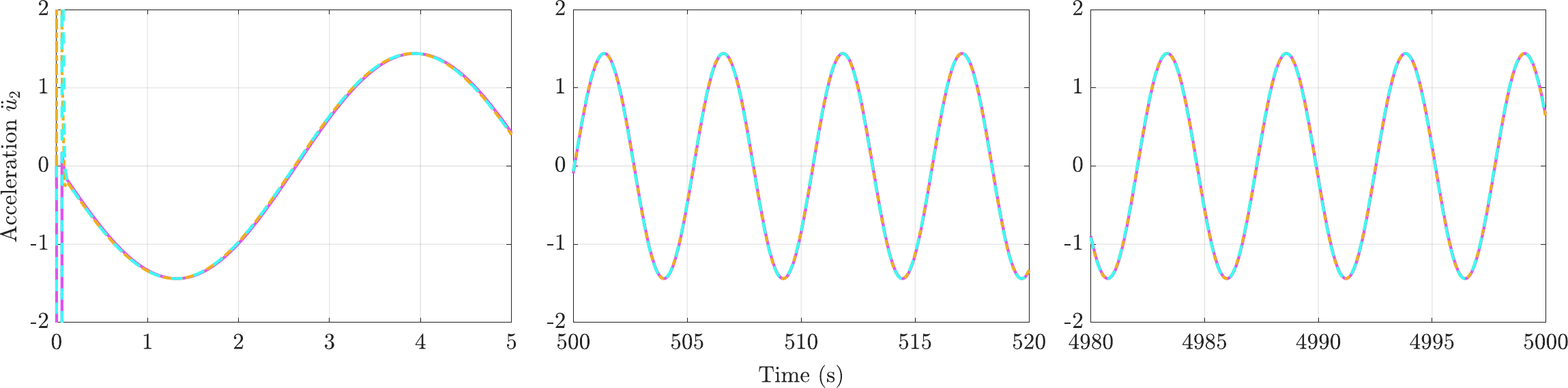}

\includegraphics[width=1\textwidth]{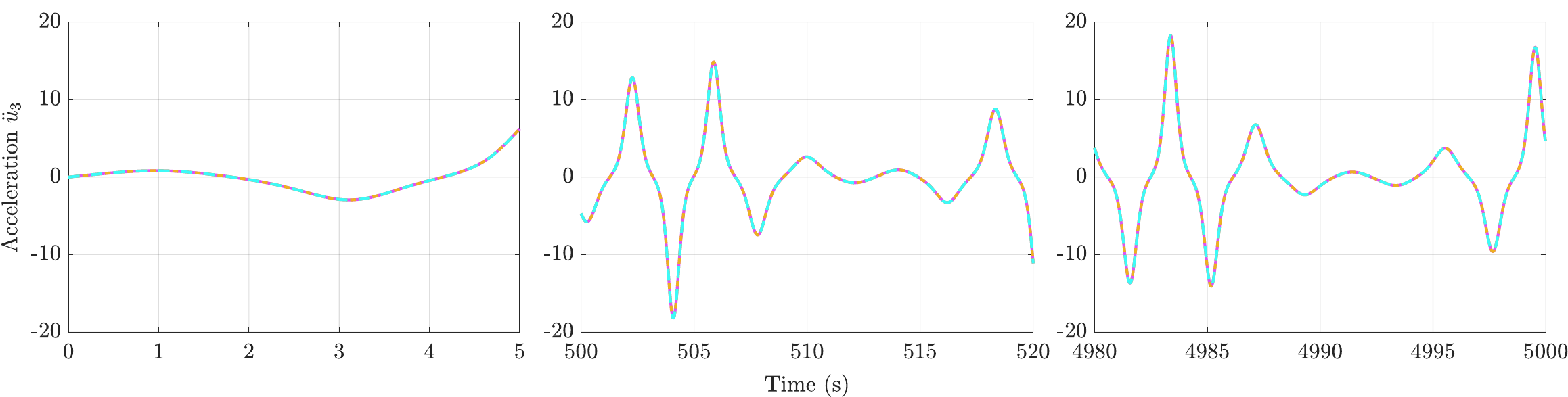}

\caption{Acceleration response of $m_{2}$ (top row), and $m_{3}$ (bottom
row), for the Padé-based schemes during $0\protect\leq t\protect\leq5$
(left column), $500\protect\leq t\protect\leq520$ (middle column)
and $4980\protect\leq t\protect\leq5000$ (right column). }\label{fig:3dofs_Pade-nonlinear-1-1}
\end{figure}

Again, to suppress the high-frequency oscillations, $\rho_{\infty}=0$
is chosen and the analysis is performed over the extended period of
$t=5000$. To demonstrate suitable numerical dissipation, a time step
of $\Delta t=0.03$ is selected for all models. The elongation, internal
force and stiffness of Spring 2 for $4950\leq t\leq5000$ are shown
in Fig.~\ref{fig:3dofs_Pade-nonlinear-stiffness} for Padé-based
schemes of order (1,2), (2,3), and (3,4), where all models effectively
overlap for the duration of the test. The nonlinear stiffening spring
ranges from 1 to approximately 18.32 in stiffness. Displacement response
of masses $m_{2}$ and $m_{3}$ are presented in Fig.~\ref{fig:3dofs_Pade-nonlinear-dsp},
with velocity responses shown in Fig.~\ref{fig:3dofs_Pade-nonlinear-vel}.
Acceleration of each mass is presented in Fig.~\ref{fig:3dofs_Pade-nonlinear-1-1}.
The Padé-based models are of order $2M-1$ due to the use of numerical
dissipation. The three columns of each figure show the responses at
three different time intervals: $0\leq t\leq5$ (left column), $500\leq t\leq520$
(middle column), and $4980\leq t\leq5000$ (right column). For displacement,
velocity and acceleration, the Padé-based schemes tested overlap on
the plots, converging on a solution for the entire duration.

\begin{figure}
\begin{centering}
\includegraphics[width=0.6\textwidth]{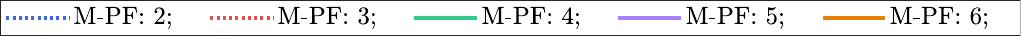}\smallskip{}
\par\end{centering}
\includegraphics[width=1\textwidth]{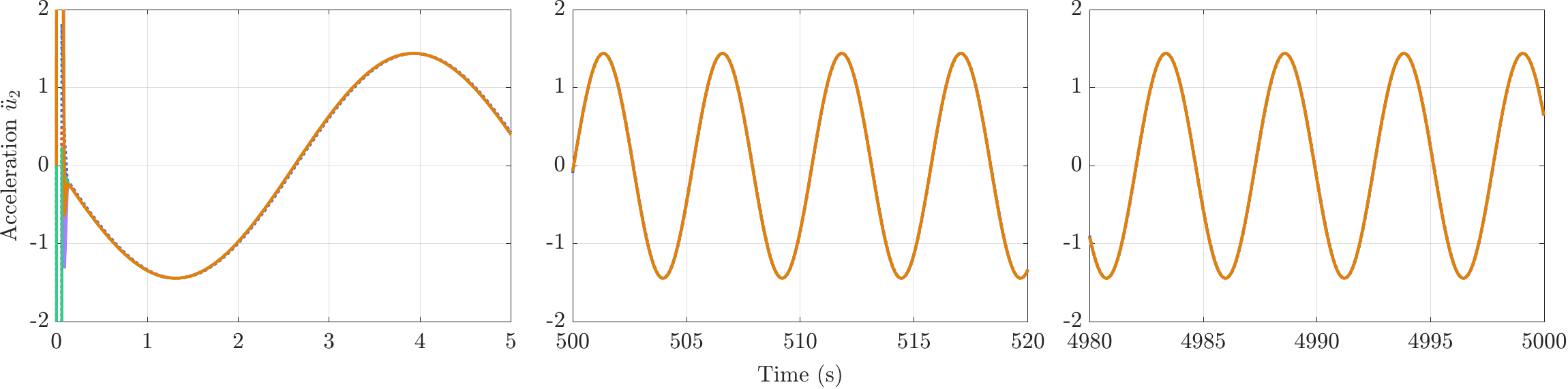}

\includegraphics[width=1\textwidth]{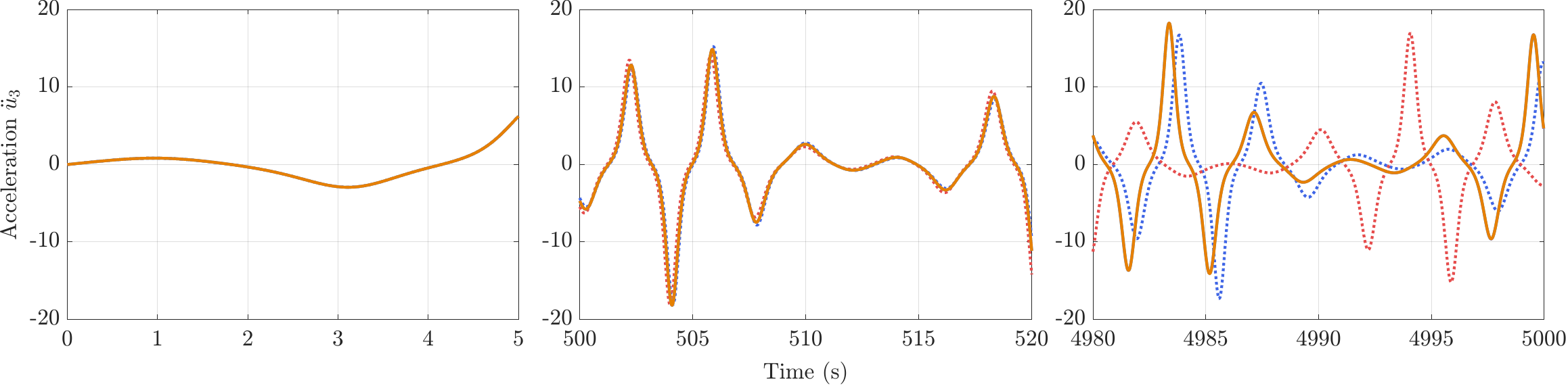}

\caption{Acceleration responses of $m_{2}$ (top row), and $m_{3}$ (bottom
row), for the $M-$schemes with single multiple root during $0\protect\leq t\protect\leq5$
(left column), $500\protect\leq t\protect\leq520$ (middle column)
and $4980\protect\leq t\protect\leq5000$ (right column). }\label{fig:3dofs_M-nonlinear}
\end{figure}

 Using $M$-based schemes (of accuracy $M$), Fig.~\ref{fig:3dofs_M-nonlinear}
shows that the second- and third-order models deviate from high-order
model solution as the simulation progresses, though all models of
fourth-order or higher converge to the same solution as the high-order
Padé-based models. Using $\rho_{\infty}=0$, all Padé-based and $M$-based
schemes dampen the oscillations in the first time step.

\begin{figure}
\begin{centering}
\includegraphics[width=0.35\textwidth]{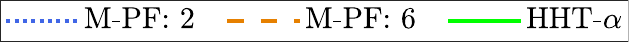}\smallskip{}
\par\end{centering}
\includegraphics[width=1\textwidth]{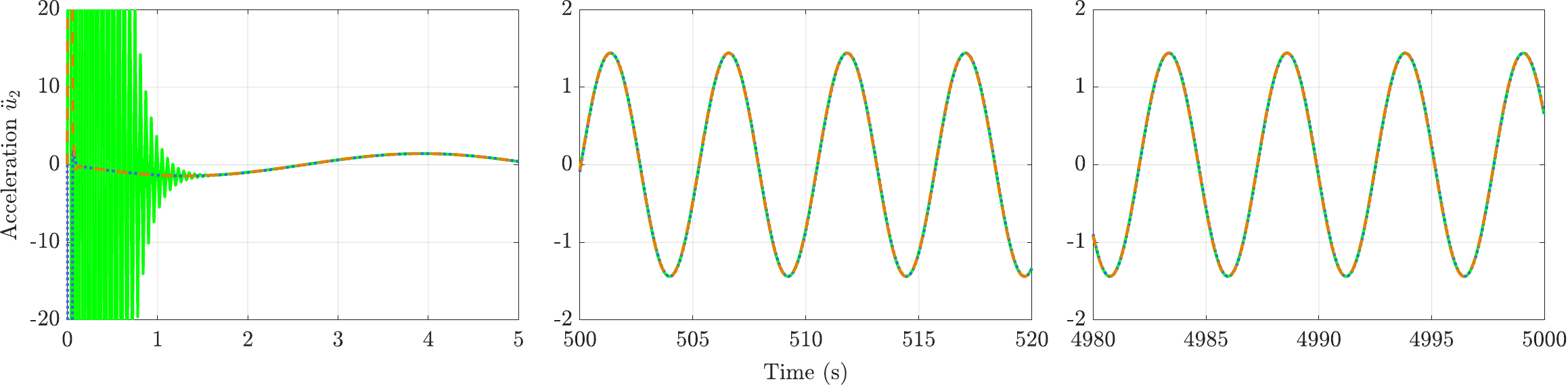}

\includegraphics[width=1\textwidth]{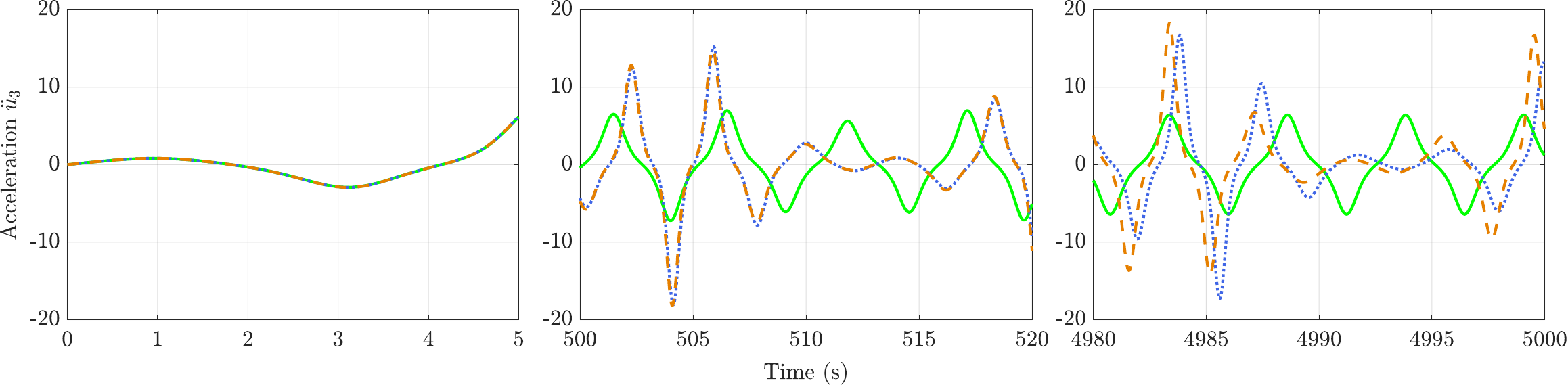}

\caption{Acceleration responses of $m_{2}$ (top row), and $m_{3}$ (bottom
row), for the HHT-$\alpha$ method ($\alpha=-0.1$), and second- and
sixth-order $M-$schemes with single multiple root for comparison;
for time durations $0\protect\leq t\protect\leq5$ (left column),
$500\protect\leq t\protect\leq520$ (middle column) and $4980\protect\leq t\protect\leq5000$
(right column). }\label{fig:3dofs_HHT-nonlinear-1}
\end{figure}

Figure~\ref{fig:3dofs_HHT-nonlinear-1} presents the response of
the HHT-$\alpha$ algorithm ($\alpha=-0.1$), with results using the
second- and sixth-order $M$-schemes also shown for comparison. In
addition to the large spurious oscillations produced for $m_{2}$,
it is noted from that the HHT-$\alpha$ algorithm demonstrates poor
accuracy in the response of $m_{3}$ over the duration of the simulation.

\subsection{One-dimensional wave propagation in a prismatic rod}\label{subsec:One-dimensional-wave-propagation-rod}

\begin{figure}
a)\includegraphics[width=0.33\textwidth]{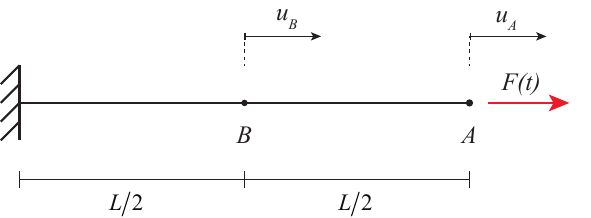}\hfill{}b)\includegraphics[viewport=0bp 0bp 388bp 311bp,clip,width=0.3\textwidth]{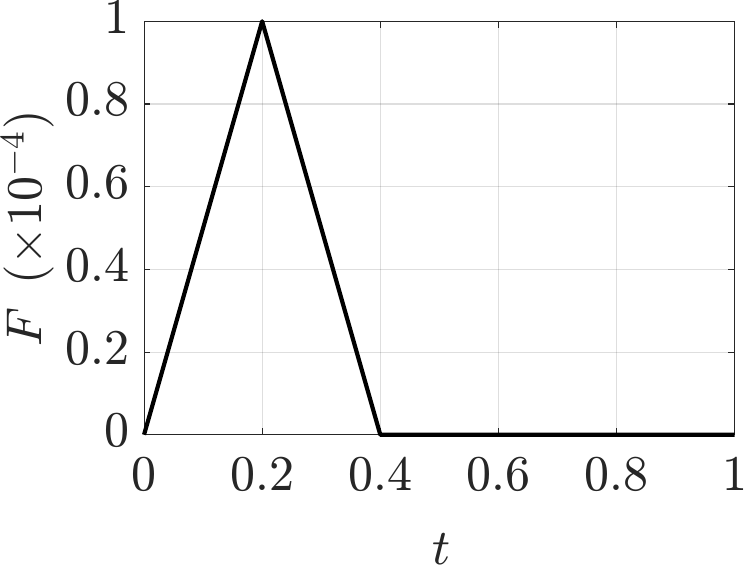}\hfill{}c)\includegraphics[viewport=0bp 0bp 388bp 311bp,clip,width=0.3\textwidth]{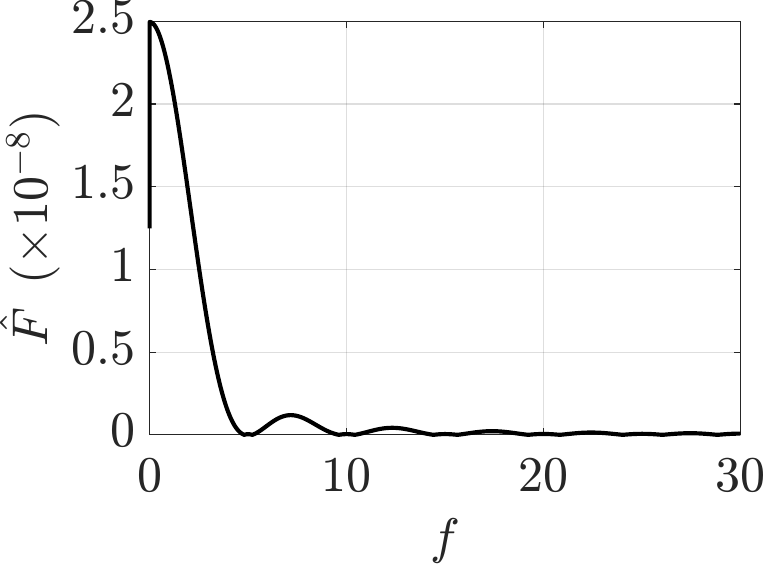}

\caption{A prismatic rod subjected to triangular impulse: a) Geometry and boundary
conditions; b) Excitation signal in the time domain; c) Excitation
signal in the frequency domain. }\label{fig:1DRod}
\end{figure}

A prismatic rod subjected to a triangular impulse of loading, as illustrated
in Fig.~\ref{fig:1DRod}, is frequently used to evaluate the performance
of time integration methods. Most studies are limited to displacement
and velocity. In this section, the acceleration response is also investigated.
A consistent set of units is used: the length $L$ and cross-section
area $A$ of the rod are chosen as unit values, i.e., $L\,{=}\,1$
and $A\,{=}\,1$. The left end of the rod is fixed. The triangular
impulse force $F(t)$ is applied at the right (free) end with a peak
value of $0.0001$ reached at $t=0.2$ (Fig.~\ref{fig:1DRod}b).
The material properties are chosen as: Young's modulus $E\,{=}\,1$,
Poisson's ratio $\nu\,{=}\,0$, and mass density $\rho\,{=}\,1$.
The d'Alembert's solution of wave equation is derived and used as
the exact solution in the following investigation. The responses at
two specific points marked by the dots in Fig.~\ref{fig:1DRod}a),
Point A at the loaded (right) end and Point B at the midpoint, are
considered.

The spatial discretization of the prismatic rod consists of $2{,}000$
linear finite elements uniformly spaced along the length (element
size of $\Delta x=5\times10^{-4}$). The Courant-Friedrichs-Lewy (CFL)
number is calculated as $\mathrm{CFL}=c\,(\Delta t/\Delta x)$ with
the longitudinal wave speed defined by $c=\sqrt{E/\rho}=1$. During
each line segment of the triangle impulse, i.e., $t\in[0,0.2]$ and
$t\in[0.2,$0.4{]}, the wave propagates through $400$ elements. The
frequency spectrum of the loading is shown in Fig.~\ref{fig:1DRod}c.
The value is negligible beyond the frequency of $20$, which is considered
as the maximum frequency of interest. At this frequency, the wavelength
is $0.05$, i.e. about the length of $100$ elements. The mesh is
considered to be very fine for earthquake engineering applications.

\begin{figure}
\centering{}\includegraphics[width=1\textwidth]{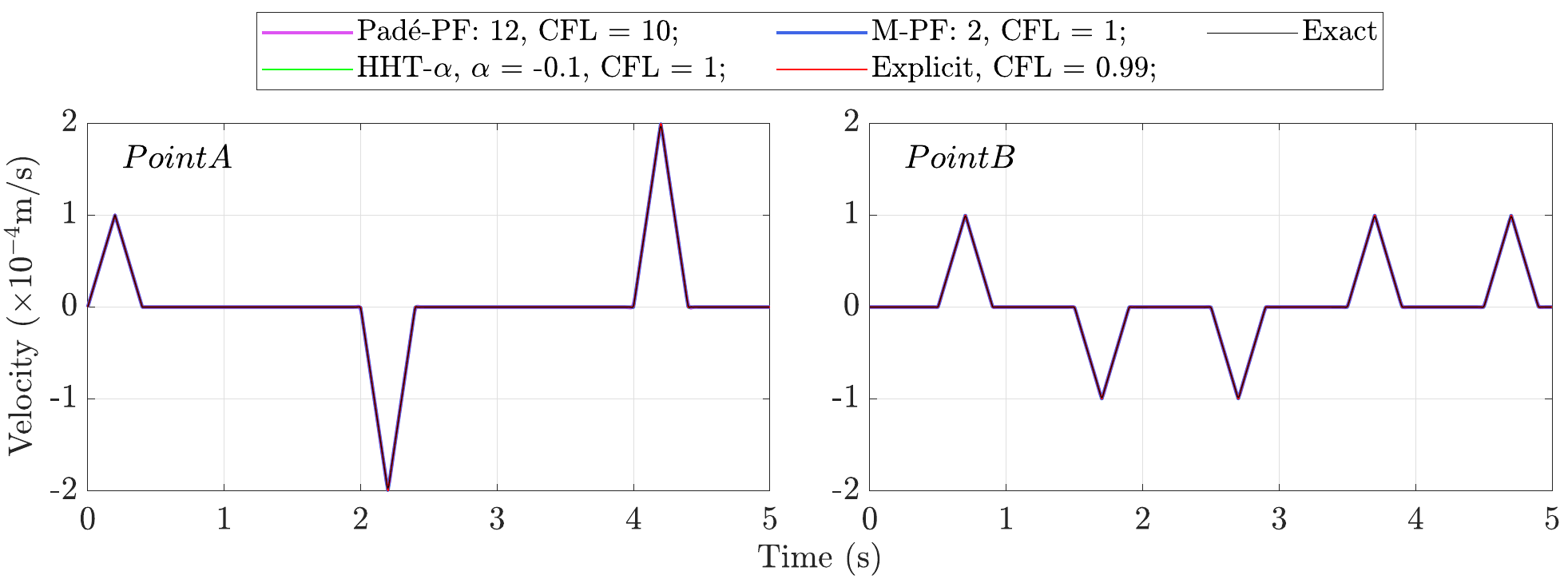}\caption{Velocity response of the rod.}\label{fig:1DRod-vel}
\end{figure}

For comparison, the HHT-$\alpha$ method and the explicit second-order
central-difference time integration method built in the commercial
finite element package Abaqus are also employed to analyze this wave
propagation problem. For the HHT-$\alpha$ method, the parameters
are chosen as $\alpha=-0.1$ and $\mathrm{CFL}=1$. For the explicit
method, automatic time-stepping in Abaqus is used which yields a $\mathrm{CFL}=0.99$.
All time integration methods perform very well in predicting the displacement
and velocity responses. For brevity, Fig.~\ref{fig:1DRod-vel} only
shows the velocity responses of a few selected schemes.

\begin{figure}
\begin{centering}
\includegraphics[width=1\textwidth]{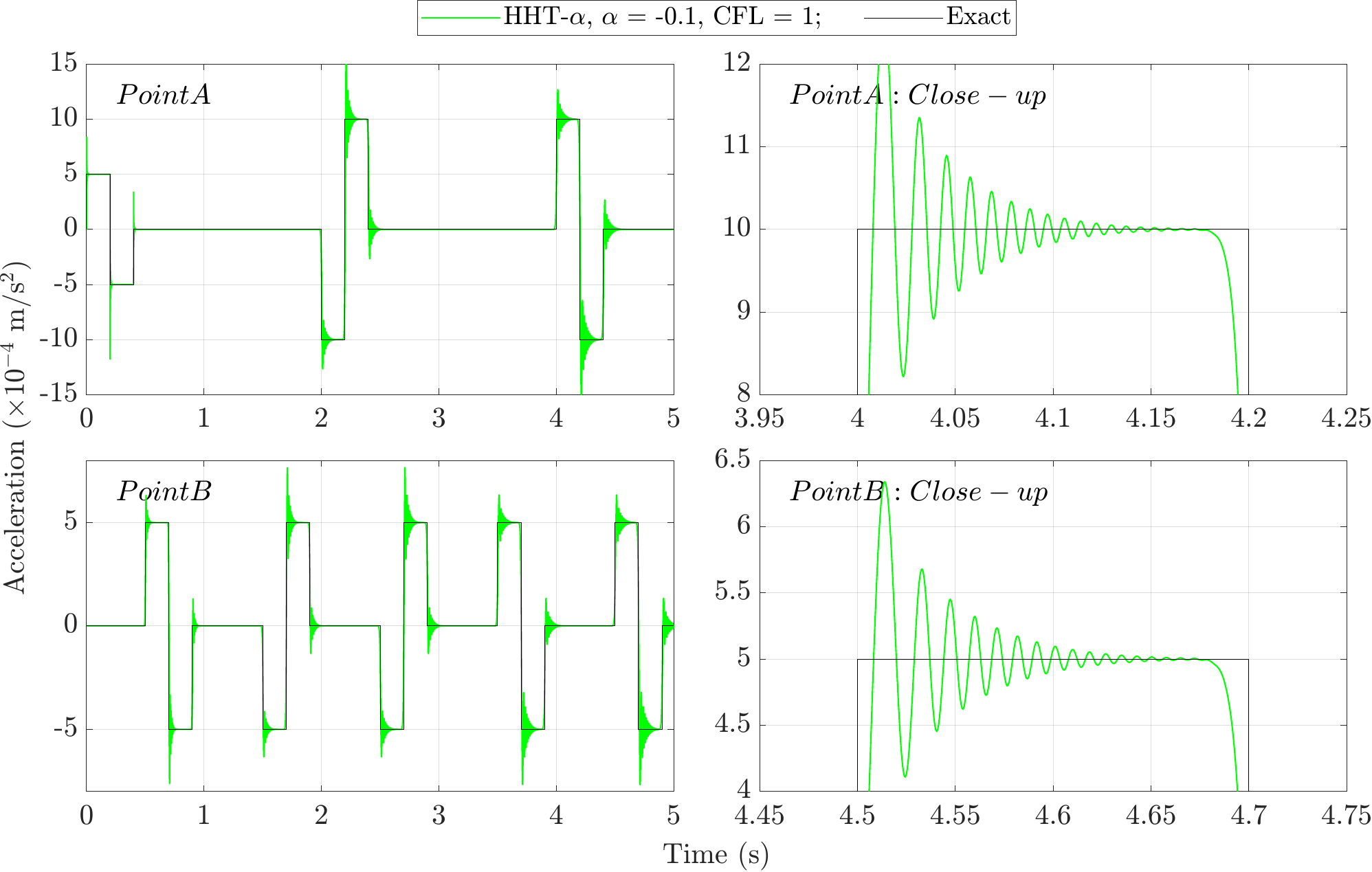}
\par\end{centering}
\caption{Acceleration response of the prismatic rod obtained with HHT-$\alpha$
method in Abaqus. }\label{fig:Rod-acc-HHT}
\end{figure}

\begin{figure}
\includegraphics[width=1\textwidth]{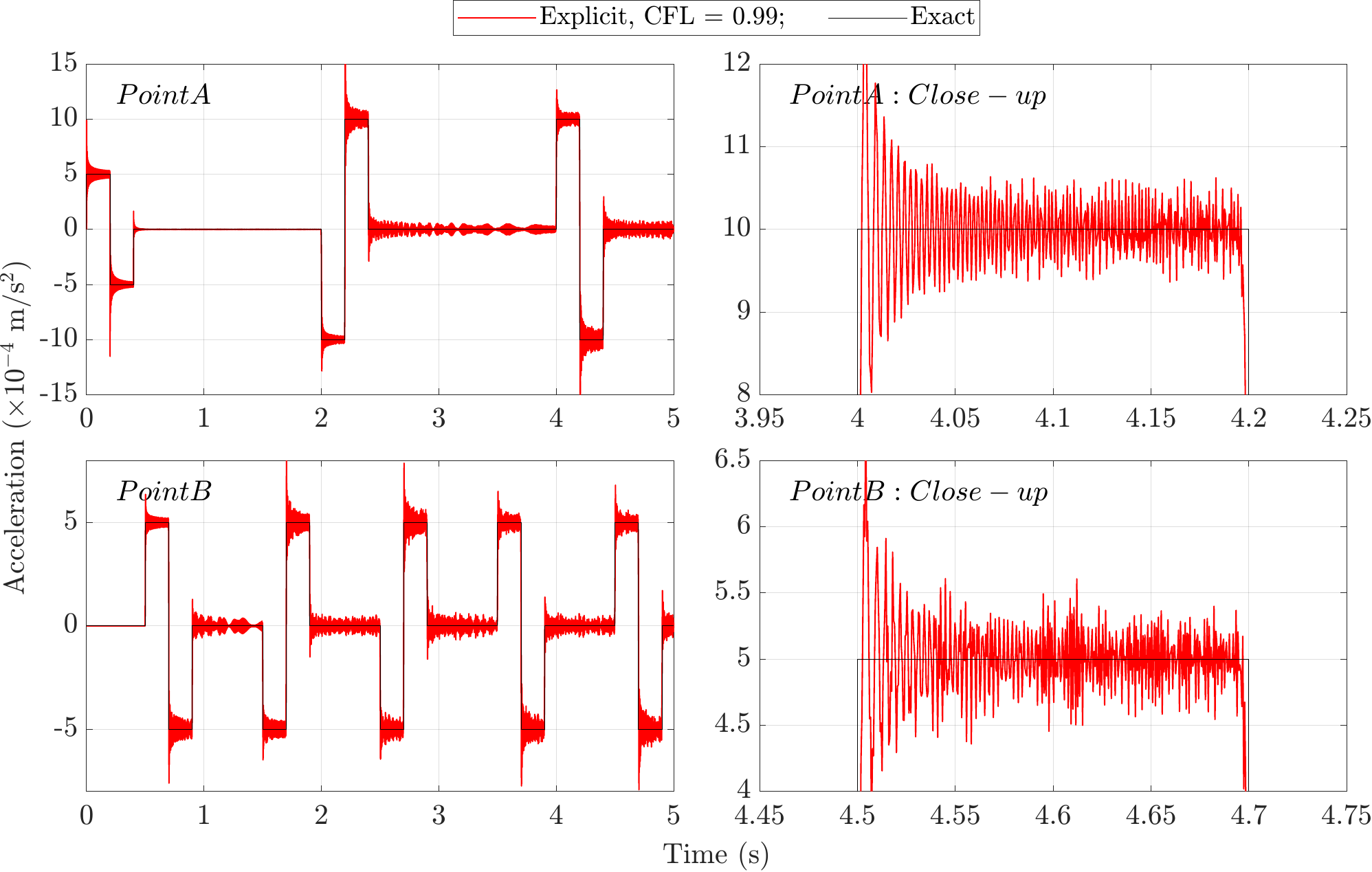}

\caption{Acceleration response of the prismatic rod obtained with explicit
method in Abaqus. }\label{fig:Rod-acc-Explicit}
\end{figure}

The acceleration responses reveal a notable difference in the performance
of the time integration methods. The results of the HHT-$\alpha$
method and the explicit method of Abaqus are plotted in Fig.~\ref{fig:Rod-acc-HHT}
and Fig.~\ref{fig:Rod-acc-Explicit}, respectively. Significant spurious
high-frequency oscillations are observed. The explicit method (Fig.~\ref{fig:Rod-acc-Explicit})
leads to noticeably stronger oscillations than the implicit method
does (Fig.~\ref{fig:Rod-acc-HHT}). 

\begin{figure}
\includegraphics[width=1\textwidth]{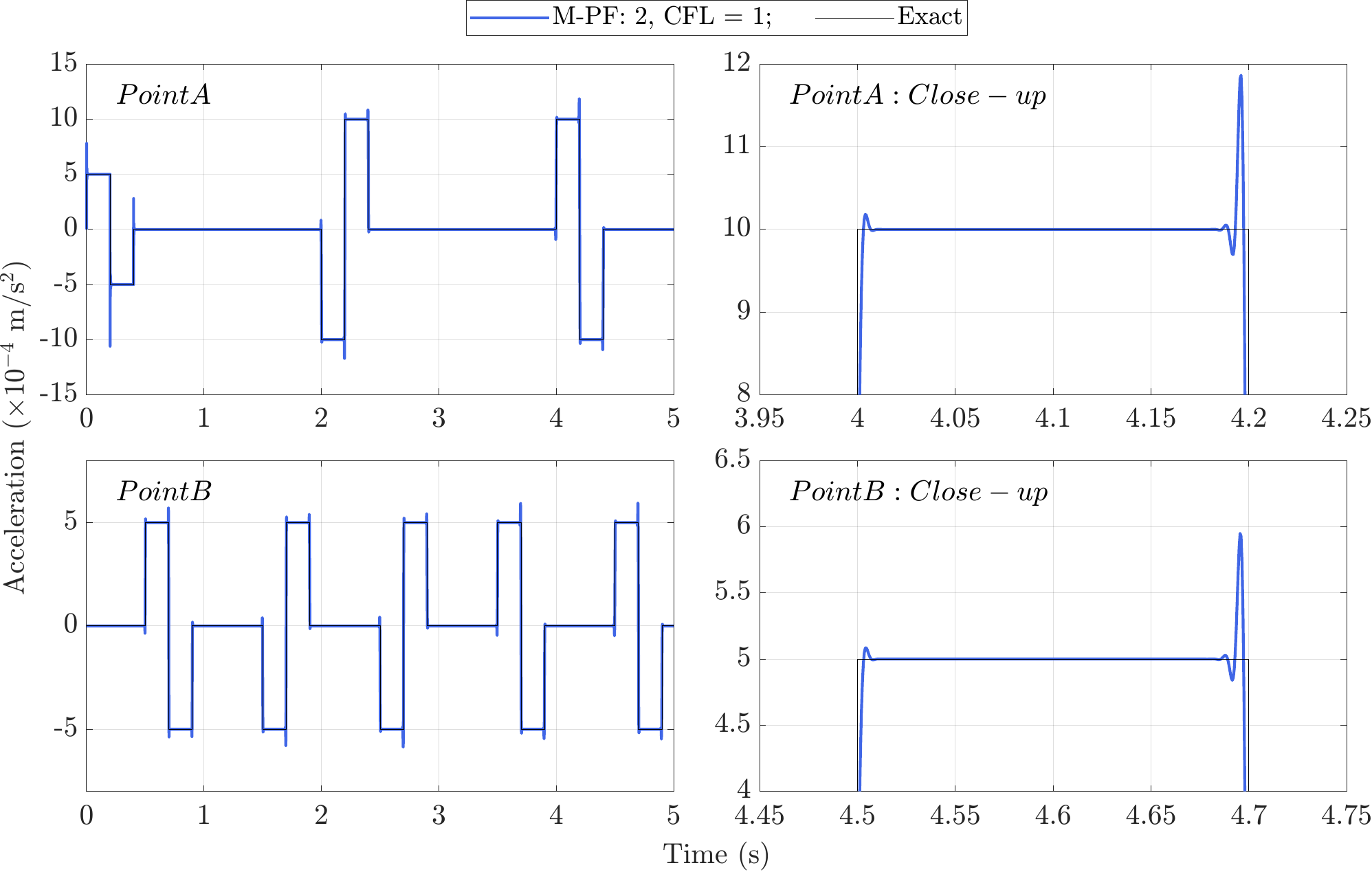}

\caption{Acceleration response of the prismatic rod obtained with the second-order
$M=2$ cases. }\label{fig:Rod-acc-M2}
\end{figure}

The Bathe method \citep{Noh2019,Kim2021} has been reported to have
better performance than the HHT-$\alpha$ method in suppressing spurious
oscillations. It is shown in \citep{Song2024} and verified for this
example that the $M=2$ composite scheme is numerically equivalent
to the second-order $\rho_{\infty}$-Bathe method. This scheme is
considered with the parameters $\rho_{\infty}=0$ and $\mathrm{CFL}=1$.
Using the proposed algorithm for single multiple roots, the acceleration
responses are obtained and shown in Fig.~\ref{fig:Rod-acc-M2}. The
spurious oscillations observed in the HHT-$\alpha$ and explicit methods
are largely suppressed. At Point A, overshoots appear within the duration
of the triangle impulse.

\begin{figure}
\includegraphics[viewport=0bp 0bp 939bp 590bp,width=1\textwidth]{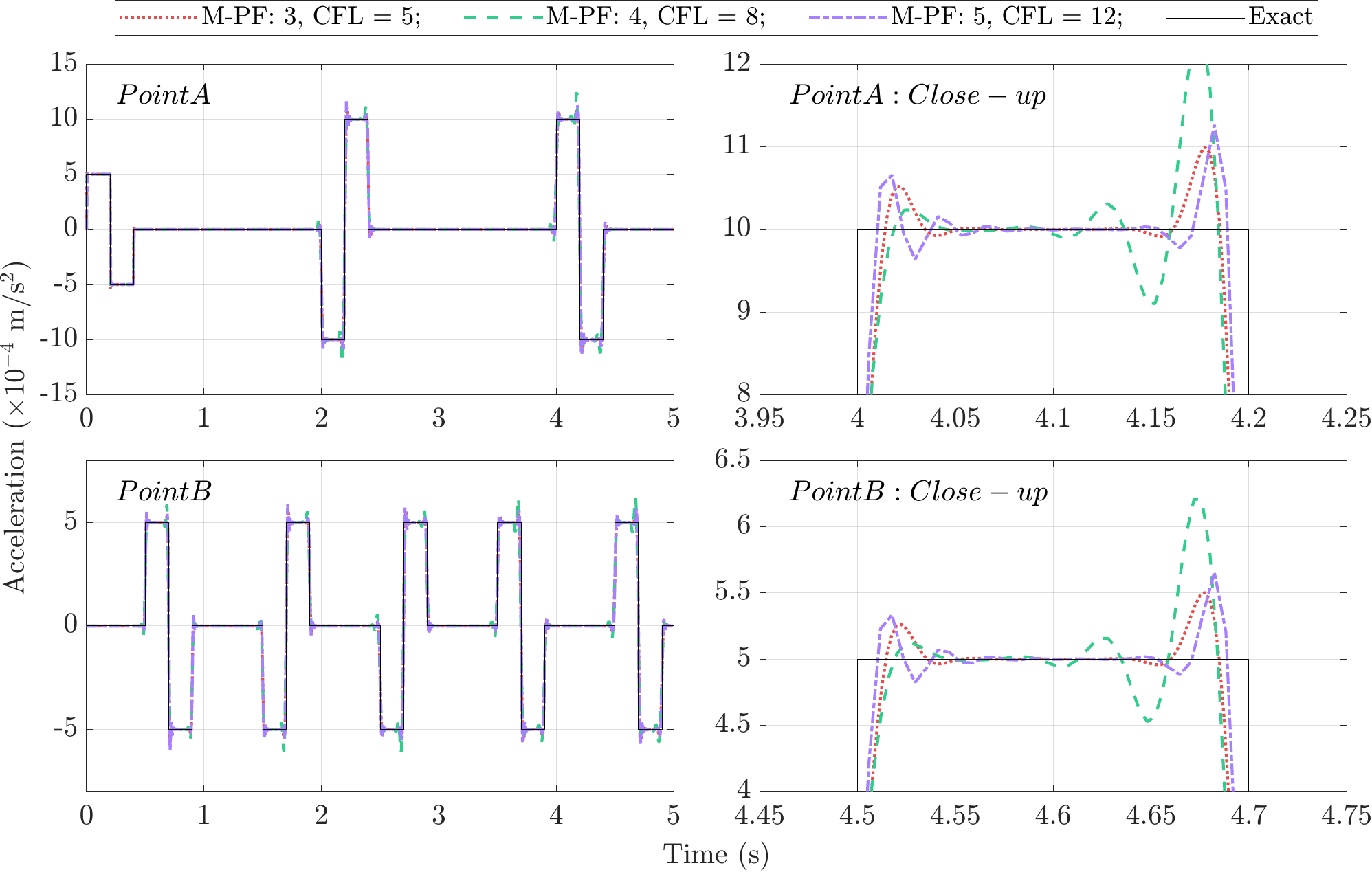}

\caption{Acceleration response of the prismatic rod obtained with the single
multiple roots cases. }\label{fig:Rod-acc-single}
\end{figure}

The high-order schemes from \citep{Song2024} are considered with
the parameter $\rho_{\infty}=0$. The proposed algorithm for the case
of single multiple roots is applied. The CFL numbers are selected
according to the order of the schemes by following the recommendation
in \citep{Song2024} as: $\mathrm{CFL}=5$ for $M=3$ (third-order),
$\mathrm{CFL}=8$ for $M=4$ (fourth-order) and $\mathrm{CFL}=12$
for $M=5$ (fifth-order). The acceleration responses are plotted in
Fig.~\ref{fig:Rod-acc-single}. It can be observed that the spurious
oscillations, including the overshoots at Point A, are effectively
suppressed here. 

\begin{figure}
\includegraphics[viewport=0bp 0bp 939bp 590bp,width=1\textwidth]{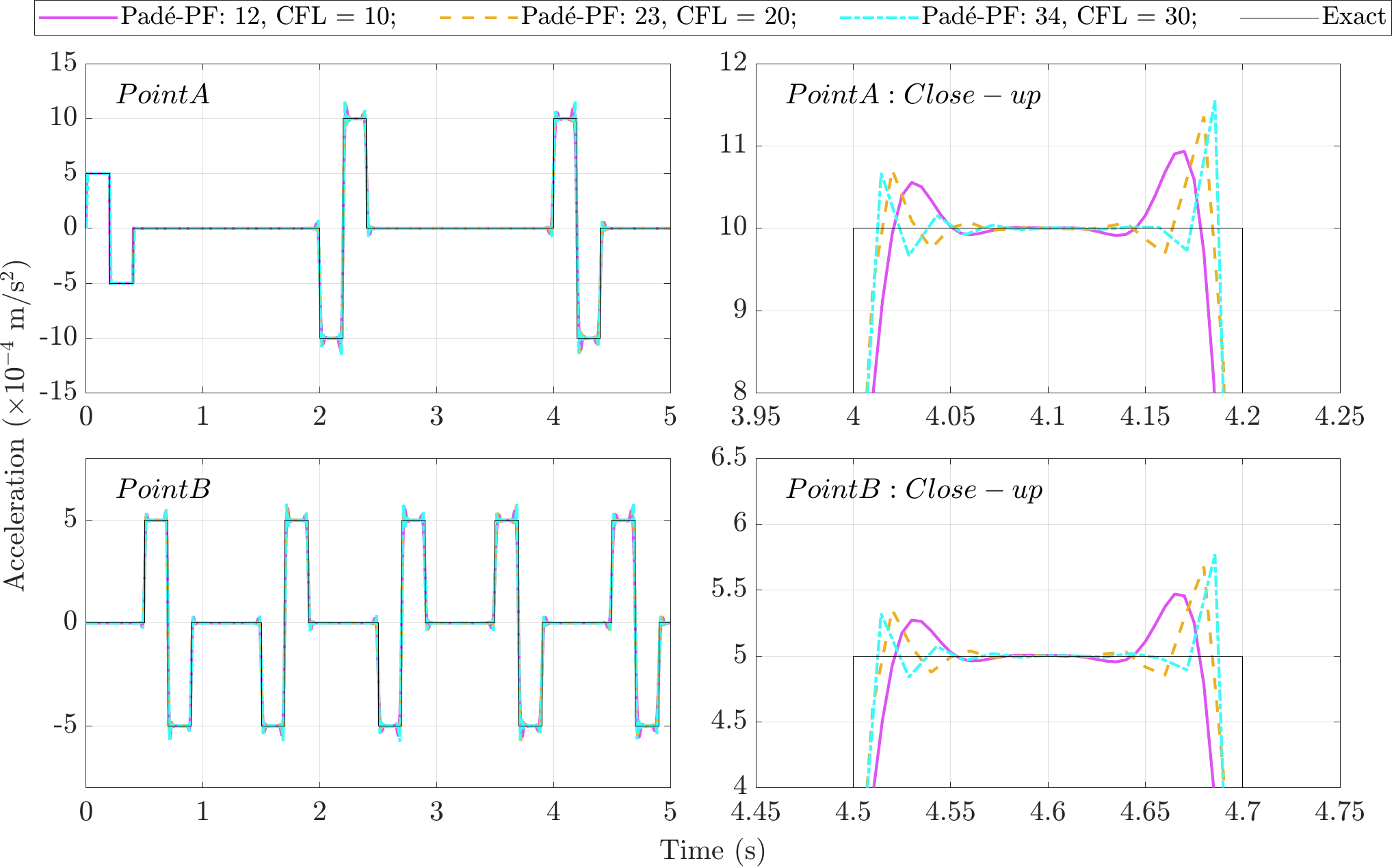}

\caption{Acceleration response of the prismatic rod obtained with distinct
roots cases. }\label{fig:Rod-acc-distinct}
\end{figure}

The high-order schemes based on the Padé expansion of the matrix exponential
in \citep{Song} are considered by applying the current proposed algorithm
for distinct roots. Again, the parameter $\rho_{\infty}=0$ is used.
The CFL numbers are chosen as $\mathrm{CFL}=10$ for the third-order
scheme (PadéPF:12), $\mathrm{CFL}=20$ for the fifth-order scheme
(PadéPF:23), and $\mathrm{CFL}=30$ for the seventh-order scheme (PadéPF:34).
The acceleration responses are shown in Fig.~\ref{fig:Rod-acc-distinct}.
It is observed that the spurious oscillations have been effectively
suppressed. The overshoots are also smaller than those in other schemes.

\subsection{Two-dimensional wave propagation in a semi-infinite elastic plane
- Lamb problem}\label{subsec:Two-dimensional-wave-propagation-Lamb}

\begin{figure}
a)\includegraphics[width=0.28\textwidth]{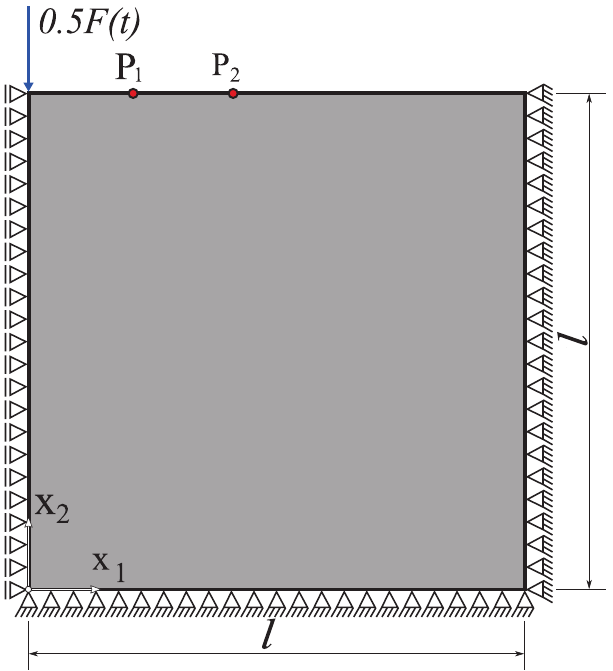}\hfill{}b)\includegraphics[width=0.3\textwidth]{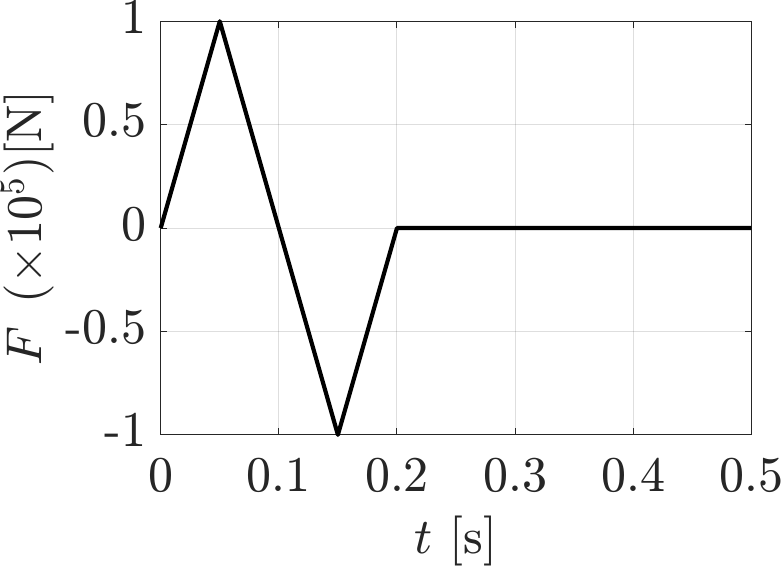}\hfill{}c)\includegraphics[width=0.3\textwidth]{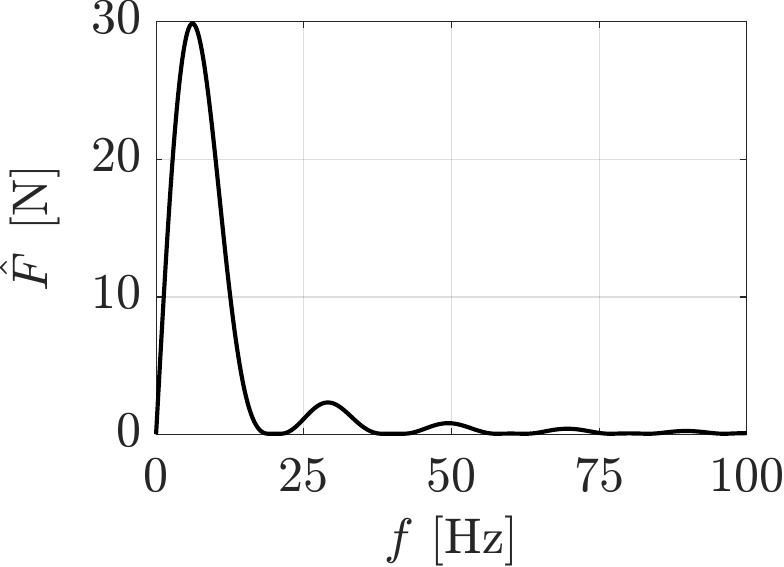}

\caption{A semi-infinite elastic domain in plane strain conditions: a) Geometry
and boundary conditions; b) Excitation signal in the time domain;
c) Excitation signal in the frequency domain. }\label{fig:geometry=000020Lamb}
\end{figure}

Lamb's problem with a vertical point load is examined under plain
strain conditions in this section. This a typical example of wave
propagation in earthquake engineering where the acceleration response
is of practical importance in a seismic design. As depicted in Fig.~\ref{fig:geometry=000020Lamb},
considering symmetry, only the right side relative to the point load
$F(t)$ is considered using the boundary conditions shown. The investigation
is conducted within a square domain of dimension $l\times l$. The
material parameters of the linear elastic half-plane are chosen the
same as in \citep{Kwon2020} and \citep{Kim2021}: Young's modulus
$E=\unit[18.77\times10^{9}]{\,Pa},$ Poisson's ratio $\nu=0.25$,
and mass density $\rho=\unit[2{,}200]{\,kg/m^{3}}$. Accordingly,
the P-wave, S-wave and Rayleigh-wave speeds are equal to $c_{p}=\unit[3{,}200\,]{m/s}$,
$c_{s}=\unit[1{,}847.5]{\,m/s}$, and $c_{R}=\unit[1{,}698.6\,]{m/s}$,
respectively. The time history of the point load is obtained by integrating
the series of step functions in~\citep{Kim2021} as
\[
F(t)=\begin{cases}
2\times10^{6}\unit[t]{N} & 0\leq t<\unit[0.05]{s}\\
\unit[10^{5}-4\times10^{6}\times(t-0.05)]{N} & \unit[0.05]{s}\leq t<\unit[0.1]{s}\\
\unit[-10^{5}+2\times10^{6}\times(t-0.1)]{N} & \unit[0.1]{s}\leq t\leq\unit[0.15]{s}
\end{cases}
\]
as shown in Fig.~\ref{fig:geometry=000020Lamb}b. Reference solutions
are given in~\citep{Kim2021} for the displacements at the two observation
points, $P_{1}\,(640,2{,}800)$ and $P_{2}\,(1{,}280,2{,}800)$, indicated
in Fig.~\ref{fig:geometry=000020Lamb}a. Under the present loading
case, these reference solutions correspond to the solutions of the
velocity responses.

To mitigate the effect of wave reflections from the boundaries, the
length of the extended square domain is chosen as $l=\unit[2{,}800]{\,m}$,
which is large enough to ensure no reflections from the artificially
fixed boundaries reach the observation points within $\unit[1\,]{s}$.
The domain is discretized into a uniform mesh consisting of $2{,}800\times2{,}800$
linear finite elements, leading to a total of $15{,}682{,}799$ degrees
of freedom. Each element is a square with a side length $\Delta x=\unit[1]{m}$.
The CFL number is calculated as $\mathrm{CFL}=c_{p}\Delta t/\Delta x$.
The mesh is sufficiently fine for earthquake engineering applications.
During one line segment of the loading, i.e. $\unit[0.05]{s}$, the
Rayleigh-wave travels about $\unit[85]{m}$, i.e. 85 elements. The
frequency spectrum of the loading is plotted in Fig.~\ref{fig:geometry=000020Lamb}b.
Considering the maximum frequency of interest as $\unit[75]{Hz}$,
the Rayleigh wavelength is about $\unit[22.6]{m}$. 

\begin{figure}
\includegraphics[width=1\textwidth]{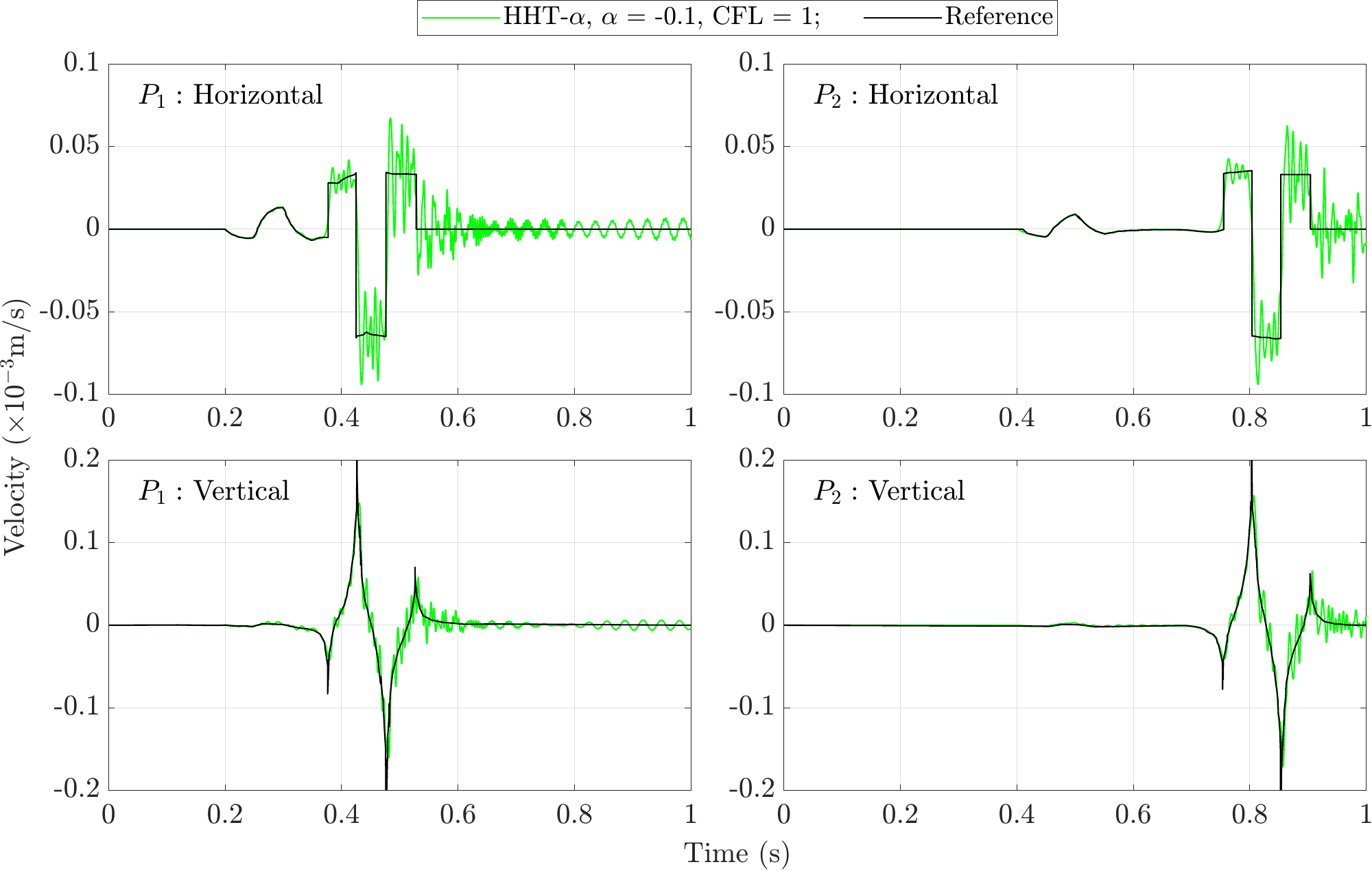}

\bigskip{}

\includegraphics[width=1\textwidth]{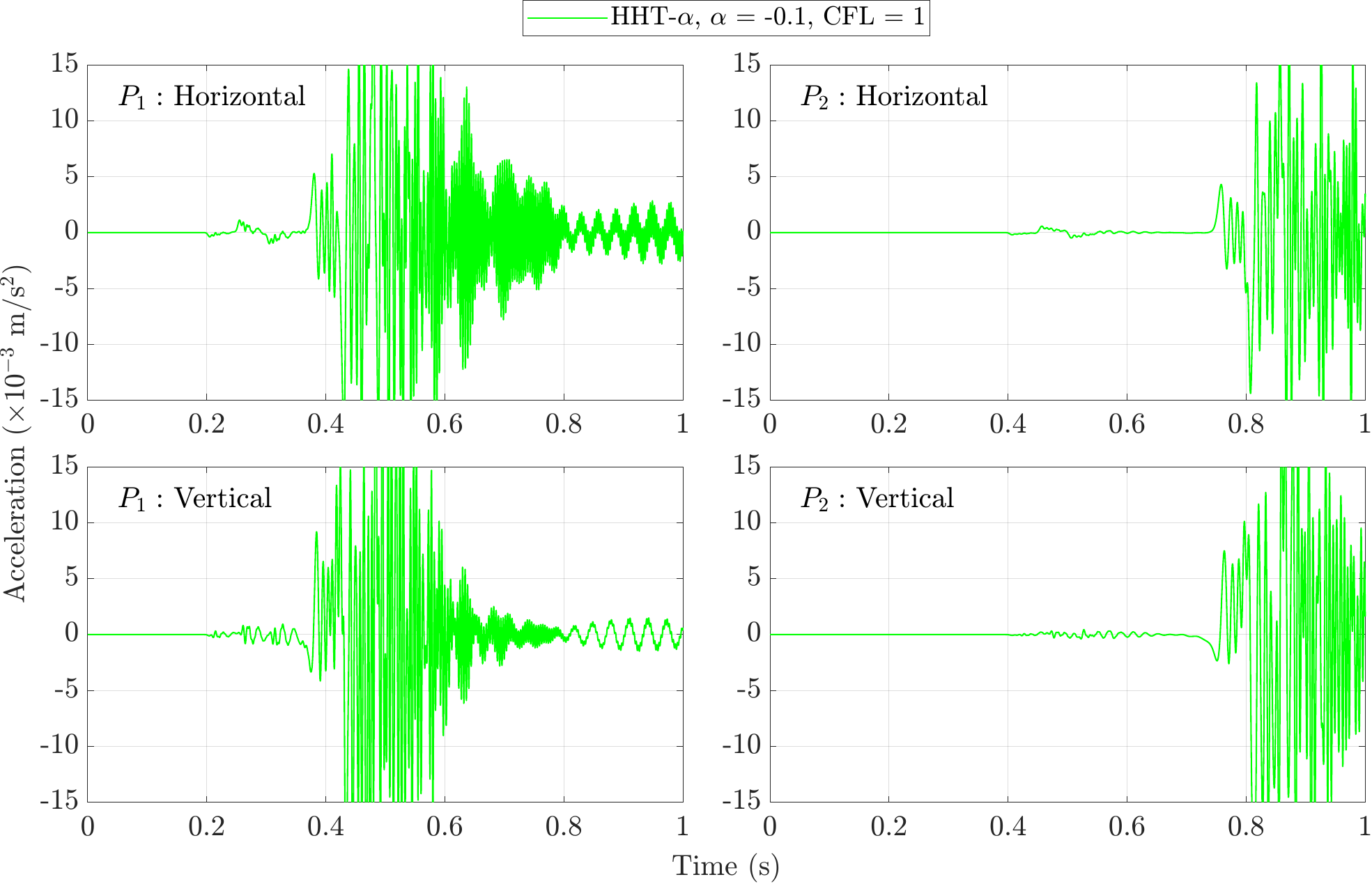}

\caption{Velocity and acceleration response of the Lamb problem obtained with
the HHT-$\alpha$ method in Abaqus. }\label{fig:Lamb-HHT}
\end{figure}

\begin{figure}
\includegraphics[width=1\textwidth]{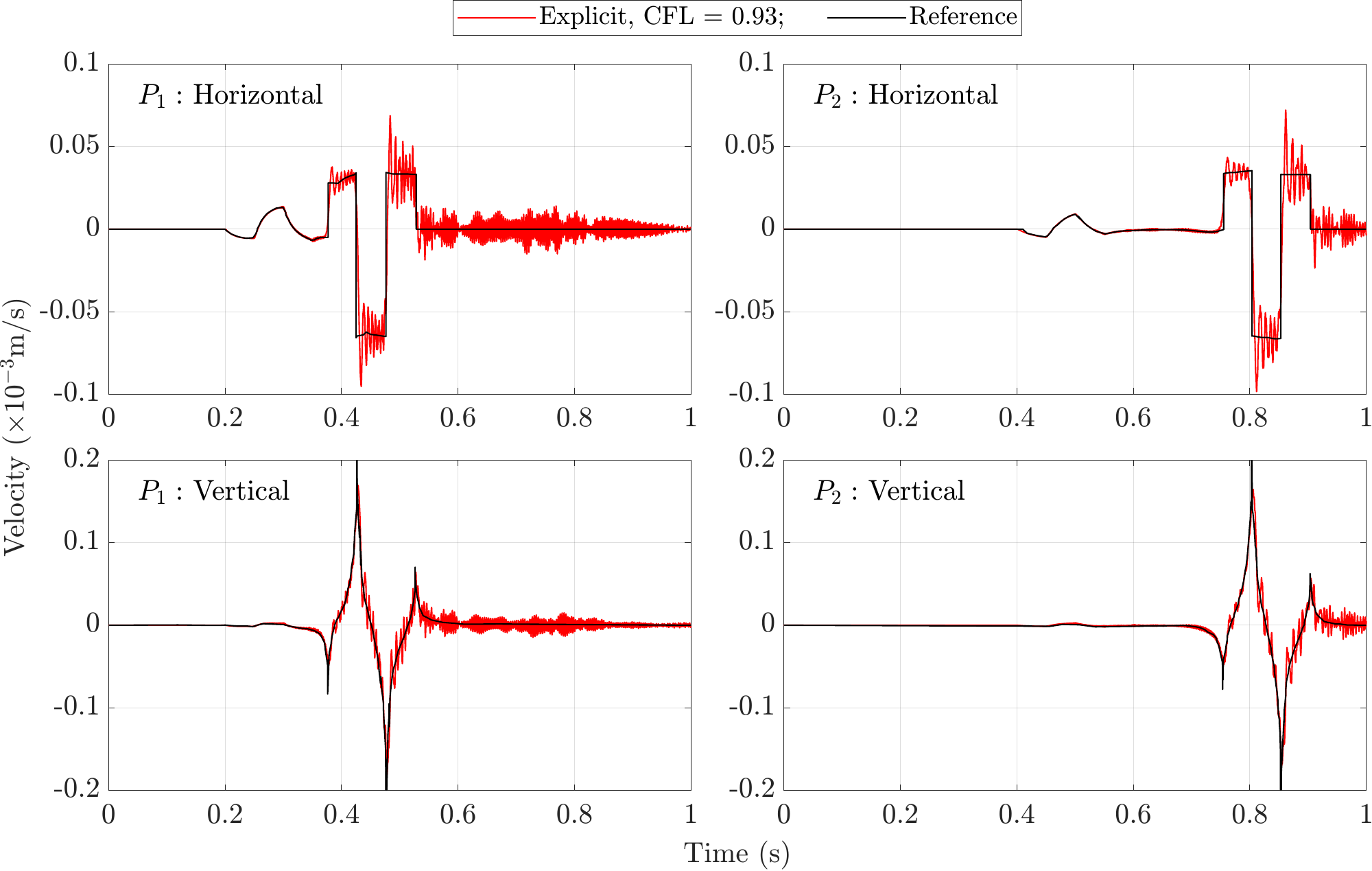}

\bigskip{}

\includegraphics[viewport=0bp 0bp 939bp 590bp,width=1\textwidth]{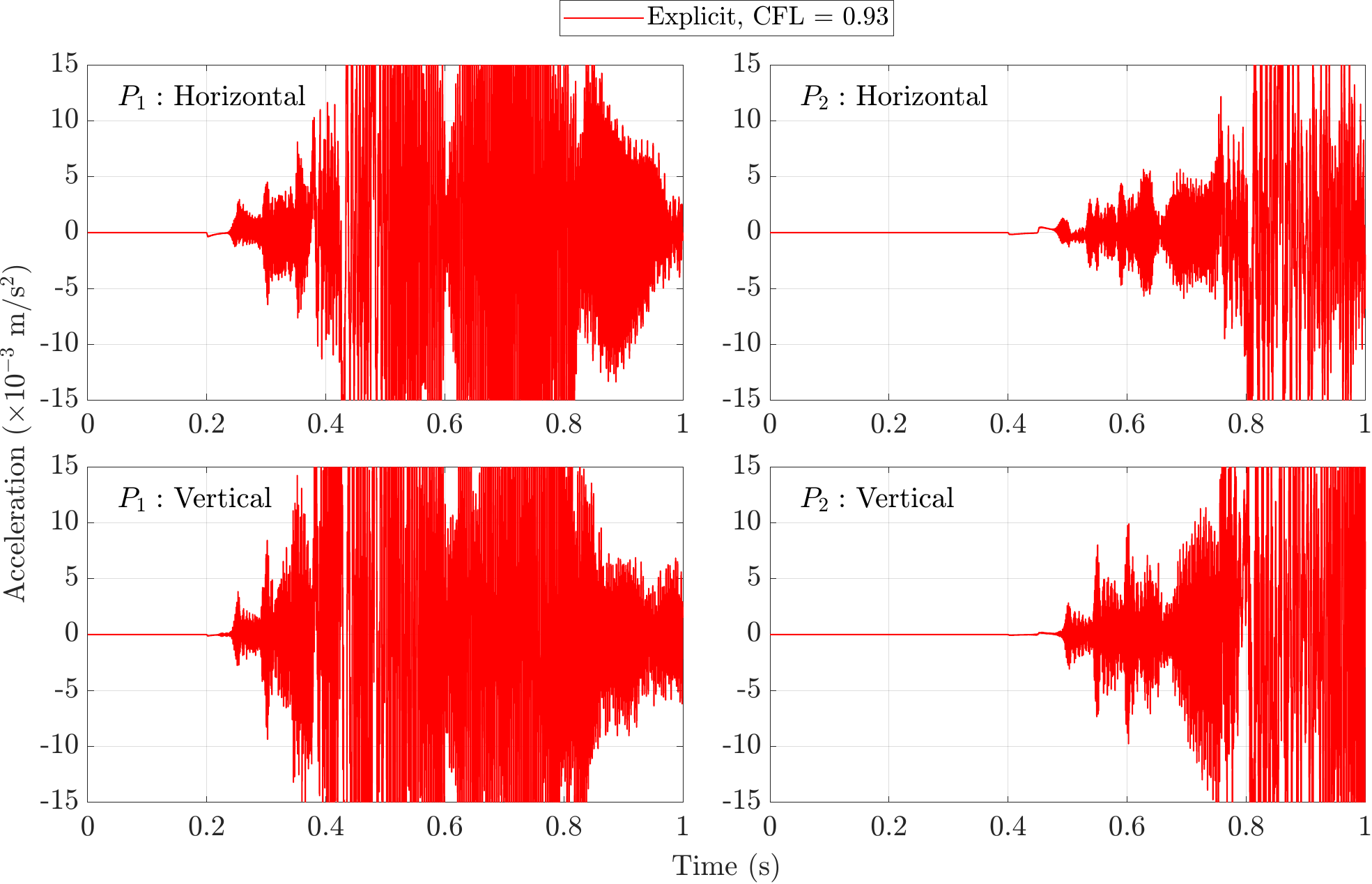}

\caption{Velocity and acceleration response of the Lamb problem obtained with
the explicit method in Abaqus. }\label{fig:Lamb-Explicit}
\end{figure}

The commercial finite element package Abaqus is first employed to
analyze this wave propagation problem. The velocity and acceleration
responses obtained using the HHT-$\alpha$ method with $\alpha=-0.1$
and $\mathrm{CFL}=1$ and the explicit method with $\mathrm{CFL}=0.93$
are shown in Figs.~\ref{fig:Lamb-HHT} and \ref{fig:Lamb-Explicit},
respectively. Note that the vertical response is larger than the horizontal
response, and different scales are used for each in the plots. Some
spurious high-frequency oscillations are observed in the velocity
responses. The explicit method exhibits stronger oscillations than
the implicit method, especially after the waves have passed. The acceleration
responses from both explicit and implicit methods are heavily polluted
by the high-frequency oscillations to the extent that the results
are not directly usable.

\begin{figure}
\includegraphics[width=1\textwidth]{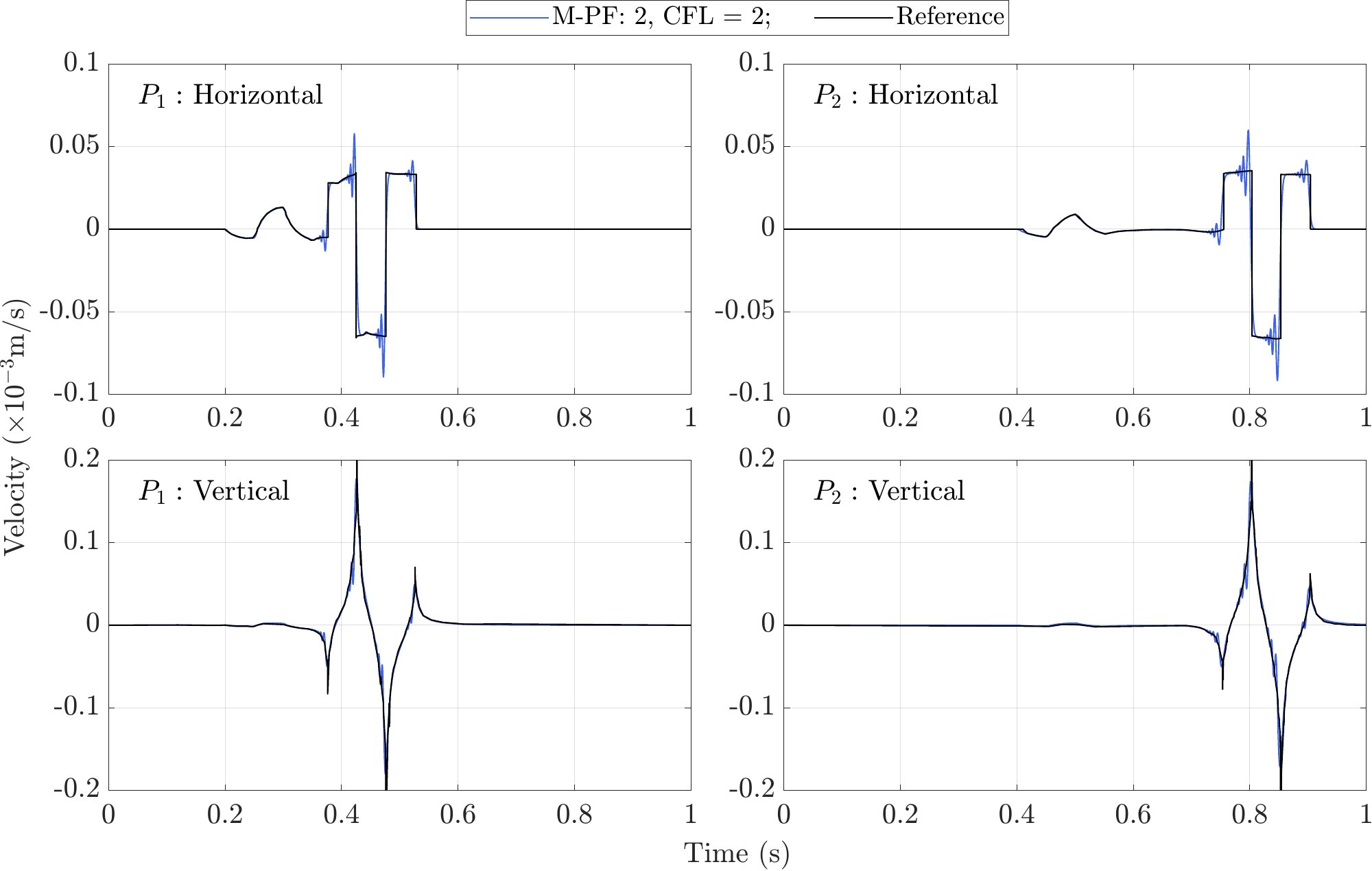}

\bigskip{}

\includegraphics[viewport=0bp 0bp 939bp 590bp,width=1\textwidth]{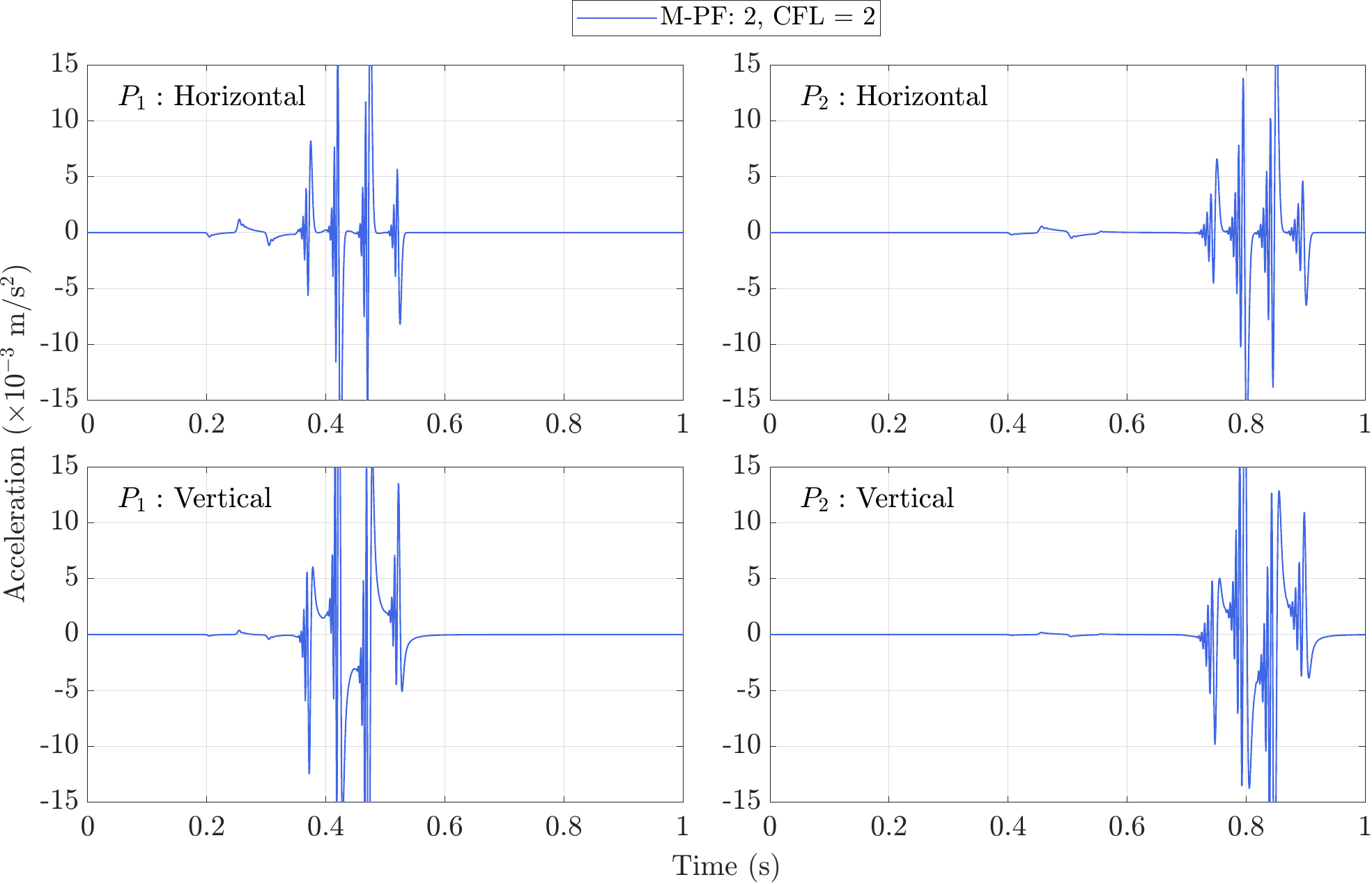}

\caption{Velocity and acceleration response of the Lamb problem obtained by
the $M=2$ composite scheme with single multiple roots. }\label{fig:Lamb-M2}
\end{figure}

The $M=2$ composite scheme, which is numerically equivalent to the
second-order $\rho_{\infty}$-Bathe method, is considered with the
parameter $\rho_{\infty}=0$. The velocity and acceleration responses
obtained using the proposed algorithm for single multiple roots are
shown in Fig.~\ref{fig:Lamb-M2} (denoted as \textquotedbl M-PF:
2\textquotedbl ). The spurious high-frequency oscillations are significantly
smaller in the velocity responses than those observed in the results
of HHT-$\alpha$ and explicit methods. Noticeable oscillations still
exist in the acceleration responses. 

\begin{figure}
\includegraphics[width=1\textwidth]{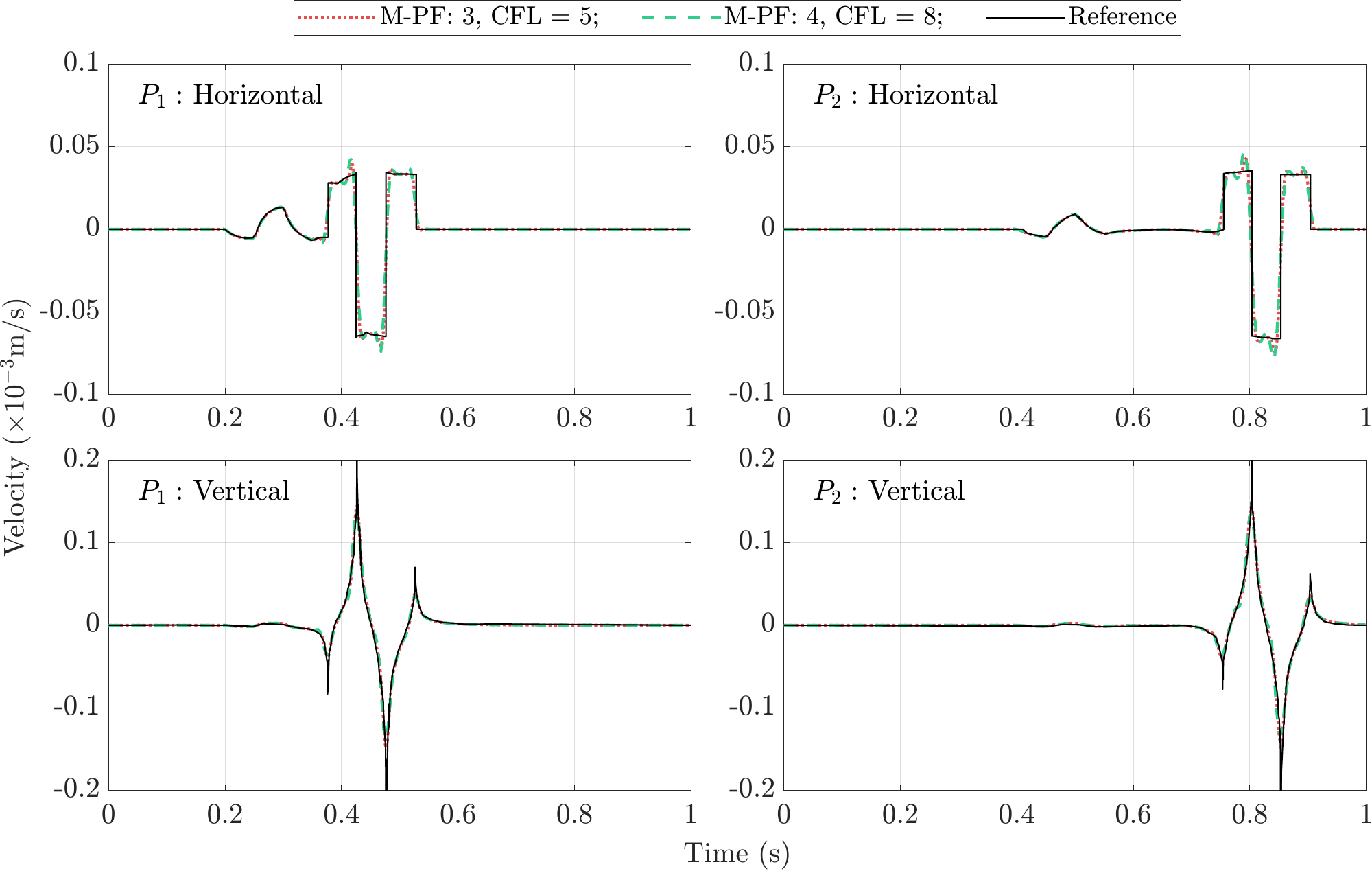}

\bigskip{}

\includegraphics[viewport=0bp 0bp 939bp 590bp,width=1\textwidth]{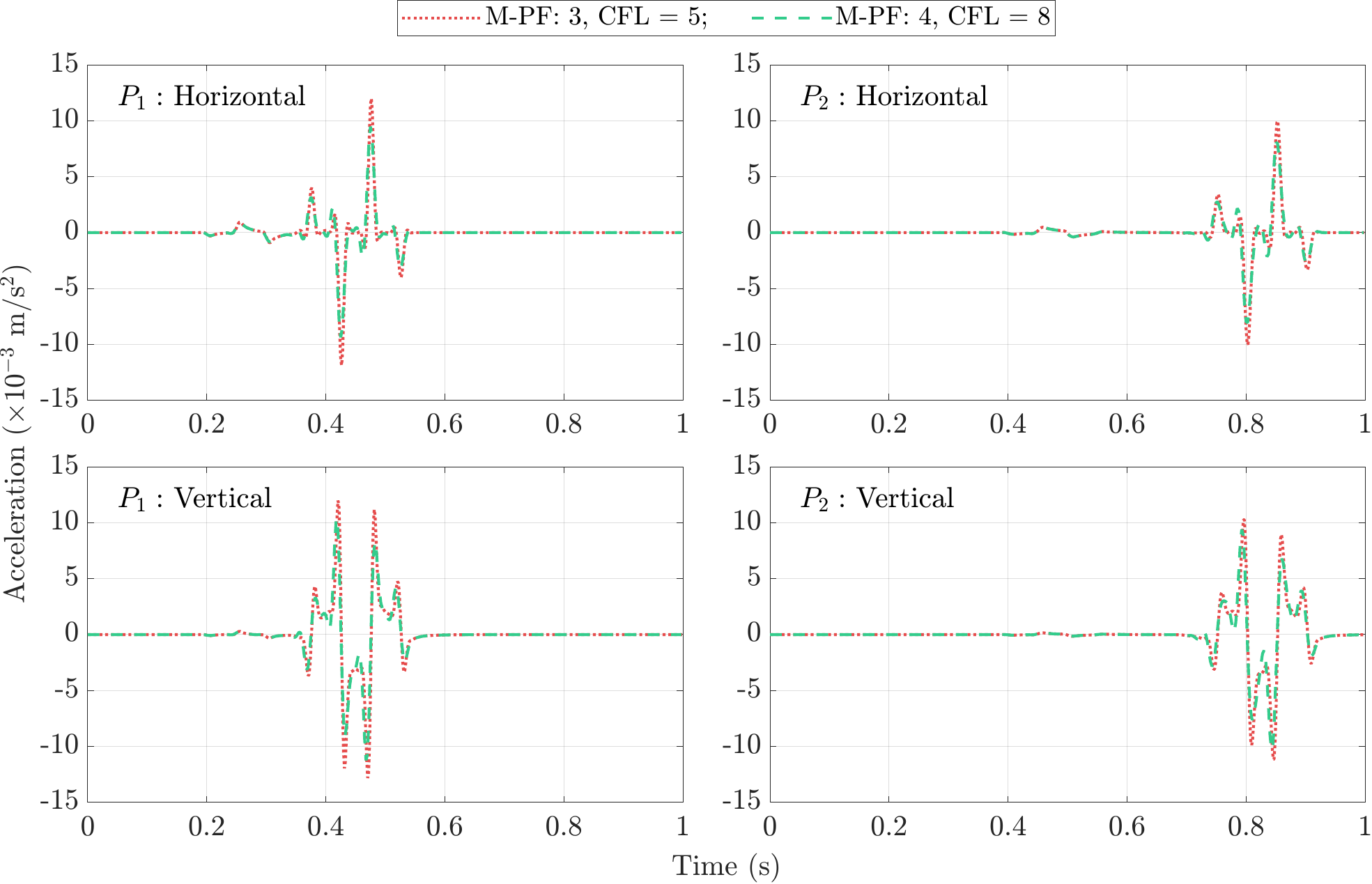}

\caption{Velocity and acceleration response of the Lamb problem obtained with
the $M=3$ and $4$ composite scheme with single multiple roots. }\label{fig:Lamb-single}
\end{figure}

The proposed algorithm is then applied with the high-order composite
schemes of orders $M=3$ and $4$. The parameter $\rho_{\infty}=0$
is used. The CFL numbers are selected as $\mathrm{CFL}=5$ for $M=3$
and $\mathrm{CFL}=8$ for $M=4$ schemes according to the recommendation
outlined in \citep{Song2024}. The proposed algorithm with single
multiple roots is applied. The velocity and acceleration responses
are shown in Fig.~\ref{fig:Lamb-single}. Good agreement with the
reference solution is observed for velocity responses. Spurious oscillations
observed in the acceleration responses predicted by the second-order
schemes are largely suppressed here. The results obtained with orders
$M=3$ and $4$ are very similar. 

\begin{figure}
\includegraphics[width=1\textwidth]{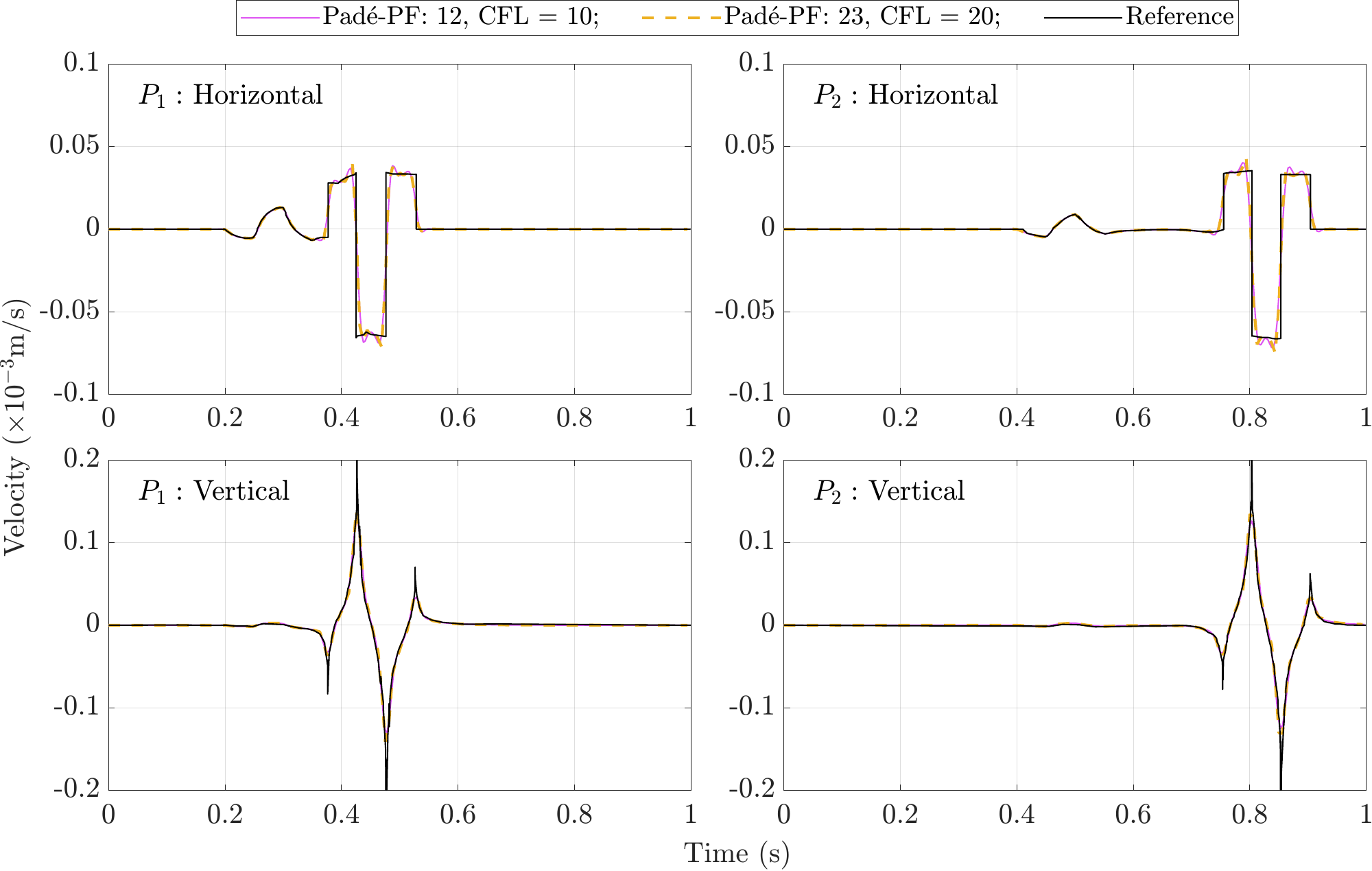}

\bigskip{}

\includegraphics[viewport=0bp 0bp 939bp 590bp,width=1\textwidth]{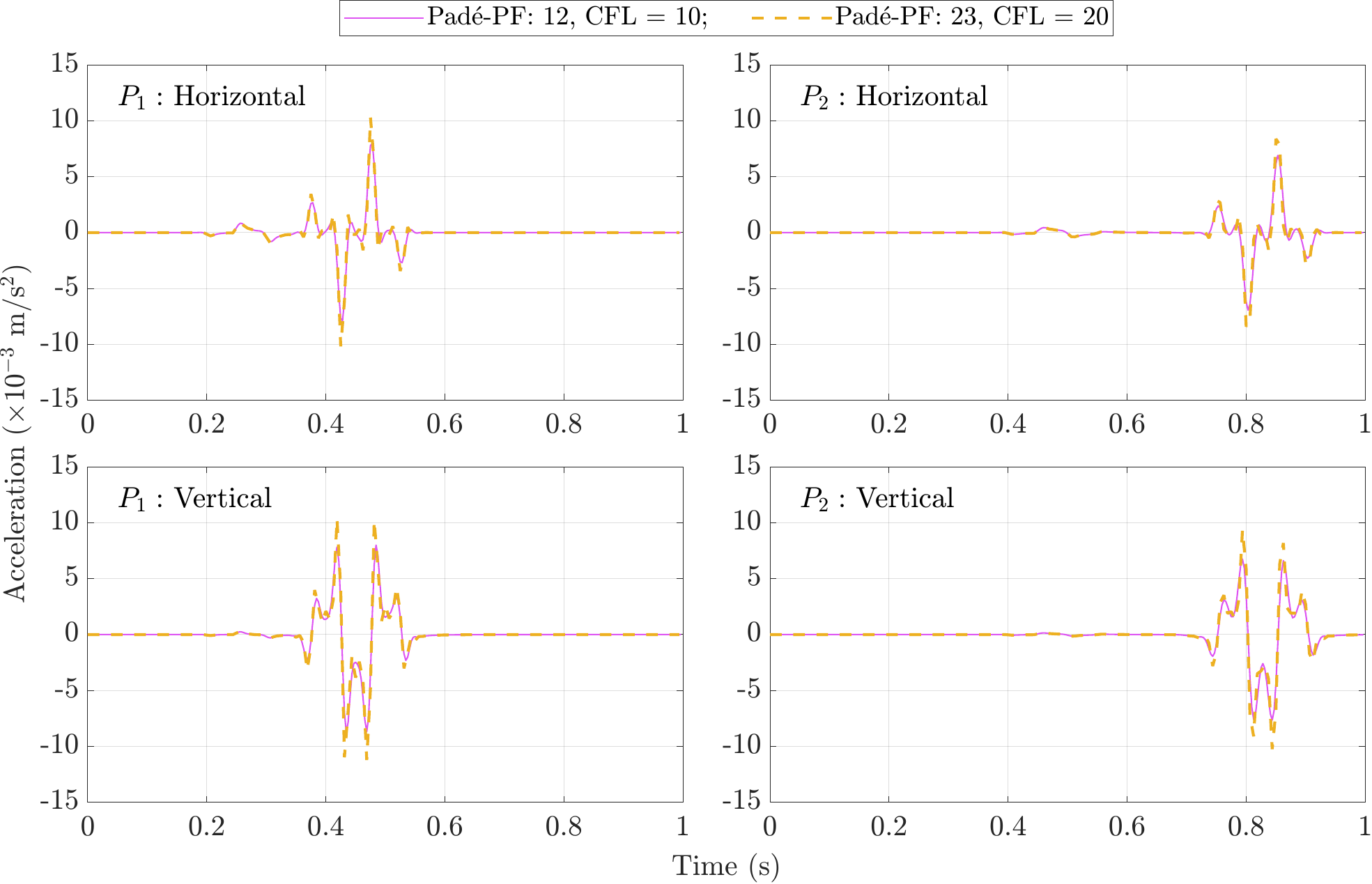}

\caption{Velocity and acceleration response of the Lamb problem obtained by
distinct roots cases. }\label{fig:Lamb-distinct}
\end{figure}

The proposed algorithm for the distinct roots case is applied with
the Padé-expansion based schemes in \citep{Song}. The parameter $\rho_{\infty}=0$
is used. The velocity and acceleration response are plotted in Fig.~\ref{fig:Lamb-distinct}.
The CFL numbers are chosen as $\mathrm{CFL}=10$ for the third-order
scheme (PadéPF:12) and $\mathrm{CFL}=20$ for the fifth-order scheme
(PadéPF:23). Again, good agreement with the reference solution of
velocity is observed. No significant spurious oscillations exist in
both velocity and acceleration responses. The results of both high-order
models are close to each other and also to the results obtained with
the high-order single multiple root composite schemes in Fig.~\ref{fig:Lamb-single}.

\section{Conclusions}\label{sec:Conclusions}

A new algorithm for a family of high-order implicit time integration
schemes based on the rational approximation of matrix exponential
is proposed. Compared with our previous works, the new algorithm removes
the need for the factorization of the mass matrix and is particularly
advantages for nonlinear problems. A new tecnique to compute the acceleration
without additional equation solutions is proposed. The acceleration
is obtained at the same order of accuracy as the displacement. As
only vector operations are used, the additional computational cost
is minimal.

It is observed from the numerical examples that the implicit HHT-$\alpha$
method, explicit second-order central-difference method, and implicit
second-order composite methods may heavily pollute the acceleration
responses for some wave propagation problems. The high-order schemes
are highly effective in suppressing the spurious high-frequency oscillations
when the order is equal to or higher than three. For nonlinear problems,
the proposed algorithms show the same convergence rates for linear
problems and very high accuracy.

\section*{Acknowledgments}

The work presented in this paper is partially supported by the Australian
Research Council through Grant Number DP200103577. 


\appendix

\section{Illustrative examples for determining coefficients}\label{sec:Illustrative-examples}

An important aspect of the algorithm presented in this work is initializing
the coefficients used in the partial fraction decomposition of the
matrix exponential, and therefore executing the numerical integration
to the desired degree of accuracy. This Appendix provides sample MATLAB
code for initializing the scheme for two cases: 1). distinct roots,
and 2). single multiple root. In both cases, the algorithms require
only the order of the numerical scheme $M$ and the user-defined parameter
$\rho_{\infty}\in[0,1]$ to be defined. Examples of the coefficients
determined for third-order ($M=3$) schemes with arbitrarily chosen
$\rho_{\infty}=0.125$ are provided for each of the two cases.

\subparagraph{Distinct roots}\label{subsec:Appendix-Distinct-roots}

Figure \ref{fig:MATLAB-code-initial_distinct} presents MATLAB code
that executes the algorithm required to determine the roots $r_{i}$,
the coefficients $p_{\mathrm{L}i}$, $a_{i}$, and matrix $\mathbf{c}$
required to initialize the time-stepping scheme for the case of distinct
roots, and therefore be able to solve Eq.\,(\ref{eq:PF-sol}). The
Padé-expansion based scheme with order $(L=2,\,M=3)$ is considered
as an example, and the spectral radius $\rho_{\infty}=0.125$ is used
(which has fifth-order accuracy when numerical dissipation is introduced).
This is run by executing $\texttt{InitSchemePadePF(3, 0.125)}.$

\begin{figure}

\texttt{function {[}rho, plr, rs, a, cfr{]} = function {[}rho, plr,
rs, a, cfr{]} = InitSchemePadePF(M,rhoInfty,pf)}

\texttt{{[}pcoe, qcoe, rs{]} = PadeExpansion(M,rhoInfty);}

\texttt{rho = pcoe(end)/qcoe(end);}

\texttt{plcoe = pcoe(1:end-1) - rho{*}qcoe(1:end-1);}

\texttt{a = polyPartialFraction(qcoe, rs);}

\texttt{plr = plcoe{*}reshape(rs,1,{[}{]}).\textasciicircum ((0:M-1)');}

\texttt{cf = TimeIntgCoeffForce(pcoe,qcoe,pf);}

\texttt{cfr = cf{*}reshape(rs,1,{[}{]}).\textasciicircum ((0:M-1)');}

\caption{MATLAB code for initializing the time-stepping scheme for case of
single multiple root}\label{fig:MATLAB-code-initial_distinct}
\end{figure}

The first line uses a function $\texttt{PadeExpansion}$ to find coefficients
of the polynomials in Eq.~(\ref{eq:Pade}). This function is provided
in Fig.\,\ref{fig:MATLAB-function-PadeExpansion}. 
\begin{figure}

\texttt{function {[}pcoe, qcoe, r{]} = PadeExpansion(M,rhoInfty)}

\texttt{L = M-1;}

\texttt{{[}p1, q1{]} = PadeCoeff(M, M);}

\texttt{{[}p2, q2{]} = PadeCoeff(M, L);}

\texttt{pcoe = rhoInfty{*}p1 + (1-rhoInfty){*}{[}p2 zeros(M-L){]};}

\texttt{qcoe = rhoInfty{*}q1 + (1-rhoInfty){*}q2;}

\texttt{r = roots(fliplr(qcoe));}

\texttt{r(imag(r)<-1.d-6) = {[}{]};}

\texttt{{[}\textasciitilde , idx{]} = sort(imag(r));}

\texttt{r = r(idx);}

\texttt{end}

\texttt{function {[}p, q{]} = PadeCoeff(M,L)}

\texttt{fc = @(x) factorial(x);}

\texttt{ii=0:L;}

\texttt{p=(fc(M+L-ii))./(fc(ii).{*}fc(L-ii));}

\texttt{ii=0:M;}

\texttt{q=(fc(M+L-ii)).{*}((-1).\textasciicircum ii)./(fc(ii).{*}fc(M-ii)){*}fc(M)/fc(L);}

\texttt{end}

\caption{MATLAB function to determine the coefficients of the Padé expansion
for case of distinct roots, and find the roots of $Q(\mathbf{A})$.
}\label{fig:MATLAB-function-PadeExpansion}
\end{figure}
Following \citep{Song}, each polynomial is defined in the following
mixed form
\begin{align}
\mathbf{P} & =\rho_{\infty}\mathbf{P}_{M/M}+\left(1-\rho_{\infty}\right)\mathbf{P}_{L/M}\,,\\
\mathbf{Q} & =\rho_{\infty}\mathbf{Q}_{M/M}+\left(1-\rho_{\infty}\right)\mathbf{Q}_{L/M}\,,
\end{align}
where the coefficients are determined through use of
\begin{align}
\mathbf{P}_{L/M} & =\sum_{i=0}^{L}\frac{\left(M+L-i\right)!}{i!\left(L-i\right)!}\mathbf{A}^{i}\,,\\
\mathbf{Q}_{L/M} & =\frac{M!}{L!}\sum_{i=0}^{M}\frac{\left(M+L-i\right)!}{i!\left(M-i\right)!}\left(-\mathbf{A}\right)^{i}\,.
\end{align}

In the current work, the partial fraction expansion uses $L=M-1$,
and for the demonstrative example case $M=3$ and $\rho_{\infty}=0.125$,
the coefficients of the following polynomials 
\begin{equation}
\mathbf{P}=67.5\mathbf{I}+28\mathbf{A}+4.125\mathbf{A}^{2}+0.125\mathbf{A}^{3}\,,\label{eq:M1scheme-PolyP}
\end{equation}
\begin{equation}
\mathbf{Q}=67.5\mathbf{I}-39\mathbf{A}+9.375\mathbf{A}^{2}-1\mathbf{A}^{3}\,.\label{eq:M1scheme-PolyQ}
\end{equation}
are returned in vectors $\texttt{pcoe}$, $\texttt{qcoe}$. Therefore,
the spectral radius is given by $\rho=p_{M}/q_{M}=-0.125$. The polynomial
$\mathbf{Q}$ in Eq.~(\ref{eq:M1scheme-PolyQ}) has one real root
and two complex-conjugate roots which may be solved to be:
\begin{equation}
r_{1}=3.7821,\ \ r_{2,3}=2.7964\pm3.1665\ \mathrm{i}\,,\label{eq:roots-example}
\end{equation}
stored in $\texttt{rs}$. The coefficients of the polynomials $P_{\mathrm{L}}(r_{i})$
in Eq.~(\ref{eq:PF_Pl(ri)}) are collected in $\texttt{plcoe}$ by
use of Eq.~(\ref{eq:Coefficient-pli}):
\begin{equation}
P_{\mathrm{L}}(r_{i})=75.9375+23.6250r_{i}+5.2969r_{i}^{2}\,.\label{eq:PL-example}
\end{equation}

Next, the function $\texttt{polyPartialFraction}$ determines parameters
$a_{i}$ for each root using Eq.~(\ref{eq:PF-expansionSln}), returning
for the example case:
\begin{equation}
a_{1}=0.909,\ \ a_{2}=a_{3}=-0.0455+0.0142\mathrm{i}\,.\label{eq:ai-example}
\end{equation}

Finally, by employing the function $\texttt{timeIntgCoeffForce}$
in Fig.\,\ref{fig:MATLAB-function-CoeffForce} the coefficient matrix
$\mathbf{c}$ is found. For the example case, Eqs.~(\ref{eq:Stepping_Ck})
and~(\ref{eq:Stepping_C0}), matrix polynomials $\mathbf{C}_{k}$
are obtained as 
\begin{equation}
C_{k}(\mathbf{A})=\left[\begin{array}{c}
67.5\mathbf{I}-5.25\mathbf{A}+1.125\mathbf{A}^{2}\\
-0.5625\mathbf{A}+0.4375\mathbf{A}^{2}\\
5.625\mathbf{I}-0.4375\mathbf{A}+0.2812\mathbf{A}^{2}\\
-0.8438\mathbf{A}+0.1094\mathbf{A}^{2}
\end{array}\right]\,.\label{eq:Mscheme-M3-Ck}
\end{equation}
and, therefore, the coefficient matrix stored as: 
\begin{equation}
\mathbf{c}=\left[\begin{array}{ccc}
67.5 & -5.25 & 1.125\\
0 & -0.5625 & 0.4375\\
5.625 & -0.4375 & 0.2812\\
0 & -0.8438 & 0.1094
\end{array}\right]\,.\label{eq:Mscheme-M3-c}
\end{equation}

The function $\texttt{InitSchemePadePF(3, 0.125)}$ in Fig.\,\ref{fig:MATLAB-code-initial_multiple}
then returns the parameter $\rho=-0.125$, as well as the coefficients
$p_{\mathrm{L}i}$ (Eq.\,(\ref{eq:PL-example})), the distinct roots
$r_{i}$ (Eq.\,(\ref{eq:roots-example})), the coefficients $a_{i}$
(Eq.\,(\ref{eq:ai-example})) and the matrix $\mathbf{c}$ (Eq.\,(\ref{eq:Mscheme-M3-c})),
as required.

\subparagraph{Single multiple root }\label{subsec:Appendix-Single-multiple-root}

Sample MATLAB code to determine the root $r$, the coefficients $p_{\mathrm{r}i}$,
and matrix $\mathbf{c}_{\mathrm{r}}$ required to initialize the time-stepping
scheme for the case of single multiple roots is provided in Fig.\,\ref{fig:MATLAB-code-initial_multiple}.
Note that the function $\texttt{shiftPolyCoefficients}$ was provided
earlier in Fig.\,\ref{fig:MATLAB-function-shiftingroots}.

\begin{figure}

\texttt{function {[}r, prcoe, cfr{]} = InitSchemeRhoPF(M, rhoInfty,
pf)}

\texttt{r = MschemeRoot(M, rhoInfty);}

\texttt{pcoe = pCoefficients(M, r);}

\texttt{prcoe = shiftPolyCoefficients(pcoe,r);}

\texttt{qrcoe = {[}zeros(1,M) 1{]};}

\texttt{qcoe = shiftPolyCoefficents(qrcoe,r);}

\texttt{cf = TimeIntgCoeffForce(pcoe,qcoe,pf);}

\texttt{cfr = shiftPolyCoefficients(cf,r);}

\caption{MATLAB code for initializing the time-stepping scheme for case of
single multiple root}\label{fig:MATLAB-code-initial_multiple}
\end{figure}
Consider the case of finding the root and coefficients for a third-order
scheme ($M=3$) with spectral radius $\rho_{\infty}=0.125$ required
to implement the final time-stepping equation for the case of single-multiple
roots as given in Eq.~(\ref{eq:SteppingEq-final}). This is equivalent
to executing the function $\texttt{InitSchemeRhoPF(3, 0.125)}$ using
Fig.\,\ref{fig:MATLAB-code-initial_multiple}. The first step uses
the function $\texttt{MschemeRoot}$ (provided in \citep{Song2024})
to solve for the single multiple root $r$. The expansion of polynomial
$p_{M}(r)$ for the third-order scheme ($M=3$) is:
\begin{equation}
p_{M}(r)=-1+3r-\frac{3}{2}r^{2}+\frac{1}{8}r^{3}=\pm\rho_{\infty}\,,\label{eq:Mscheme-PolyP}
\end{equation}
where, for $M=3$ the RHS is set equal to $-\rho_{\infty}$ (see \citep{Song2024})
\begin{equation}
-1+3r-\frac{3}{2}r^{2}+\frac{1}{6}r^{3}=-0.125\,,\label{eq:eq:Mscheme-M2-rEquation}
\end{equation}
which, by solving this cubic equation and testing the three roots
by considering the stability of the schemes, the root 
\begin{equation}
r=2.3917\label{eq:Mscheme-M2-r}
\end{equation}
is selected. Using the function $\texttt{pCoefficients}$ (provided
in \citep{Song2024}) with the newly found root $r=2.3917$, the polynomial
$P(\mathbf{A})$ can be given as follows 
\begin{equation}
P(\mathbf{A})=13.6802\mathbf{I}-3.4798\mathbf{A}-3.1449\mathbf{A}^{2}-0.125\mathbf{A}^{3}\,.\label{eq:Mscheme-M2-1-p}
\end{equation}
with coefficients stored in $\texttt{pcoe}$. Through use of Eq.~(\ref{eq:Ar-def}),
the algorithm $\texttt{shiftPolyCoefficients}$ transforms Eq.\,(\ref{eq:Mscheme-M2-1-p})
to a polynomial in terms of the matrix $\mathbf{A}_{r}$ 
\begin{align}
P_{\mathrm{r}}(\mathbf{A}_{\mathrm{r}}) & =13.6802\mathbf{I}-3.4798(2.3917\thinspace\mathbf{I}-\textbf{A}_{\mathrm{r}})-3.1449(2.3917\thinspace\mathbf{I}-\mathbf{A}_{\mathrm{r}})^{2}-0.125(2.3917\thinspace\mathbf{I}-\mathbf{A}_{\mathrm{r}})^{3}\nonumber \\
 & =-14.3410\mathbf{I}+20.6678\mathbf{A}_{\mathrm{r}}-4.0418\mathbf{A}_{\mathrm{r}}^{2}+0.125\mathbf{A}_{r}^{3}\,,\label{eq:Mscheme-M2-Pr}
\end{align}
and returns the coefficients in the vector $\texttt{prcoe}$. The
denominator of the expansion in Eq.\,(\ref{eq:R-def}) is equivalent
to $\mathbf{A}_{r}^{M}$, and therefore all coefficients of polynomial
$\mathbf{Q}_{r}$ are zero except the term of order $M$ as defined
for the variable $\texttt{qrcoe}$. Coefficients of $\mathbf{Q}$
are found using the shifting algorithm $\texttt{shiftPolyCoefficients}$. 

The MATLAB function $\texttt{TimeIntgCoeffForce}$ for determining
the coefficient matrix related to the expression in Eq.~(\ref{eq:Mscheme-M2-Ck})
is provided in Fig.~\ref{fig:MATLAB-function-CoeffForce}. 
\begin{figure}

\texttt{function {[}C{]} = TimeIntgCoeffForce(p, q, pf)}

\texttt{M = length(q) - 1;}

\texttt{tmp = p - q;}

\texttt{C = zeros(pf+1,M);}

\texttt{C(1,:) = tmp(2:end); \%ltx order reduced by one \{A\}\textasciicircum\{-1\}\textbackslash left(\{P\}-\{Q\}\textbackslash right) }

\texttt{for k = 1:pf}

\texttt{tmp = ((-1/2)\textasciicircum k){*}(p-((-1)\textasciicircum k){*}q);
\%ltx order of \$C\_\{k-1\}\$ is lower than \$q\$ by one}

\texttt{tmp(1:M) = tmp(1:M) + k{*}C(k,:); }

\texttt{C(k+1,:) = tmp(2:end);}

\texttt{end}

\caption{MATLAB function to determine the coefficients required to integrate
the non-homogeneous term $C_{k}(\mathbf{A})$. }\label{fig:MATLAB-function-CoeffForce}
\end{figure}
Following Eqs.~(\ref{eq:Stepping_Ck}) and Eqs.~(\ref{eq:Stepping_C0}),
the matrix polynomials $\mathbf{C}_{k}$ can be obtained as 
\begin{equation}
C_{k}(\mathbf{A})=\left[\begin{array}{c}
13.6802\mathbf{I}-10.3199\mathbf{A}+0.875\mathbf{A}^{2}\\
-1.14\mathbf{A}+0.5625\mathbf{A}^{2}\\
1.14\mathbf{I}-1.455\mathbf{A}+0.2187\mathbf{A}^{2}\\
-1.7849\mathbf{I}+0.1525\mathbf{A}+0.1406\mathbf{A}^{2}
\end{array}\right]\,,\label{eq:Mscheme-M2-Ck}
\end{equation}
which once more, by employing Eq.~(\ref{eq:Ar-def}), are re-written
in terms of $\mathbf{A}_{\mathrm{r}}$: 
\begin{align}
C_{\mathrm{r}k}(\mathbf{A}_{\mathrm{r}}) & =\left[\begin{array}{c}
13.6802\mathbf{I}-10.3199(2.3917\thinspace\mathbf{I}-\mathbf{A}_{\mathrm{r}})+0.875(2.3917\thinspace\mathbf{I}-\mathbf{A}_{\mathrm{r}})^{2}\\
-1.14(2.3917\thinspace\mathbf{I}-\mathbf{A}_{\mathrm{r}})+0.5625(2.3917\thinspace\mathbf{I}-\mathbf{A}_{\mathrm{r}})^{2}\\
1.14\mathbf{I}-1.455(2.3917\thinspace\mathbf{I}-\mathbf{A}_{\mathrm{r}})+0.2187(2.3917\thinspace\mathbf{I}-\mathbf{A}_{\mathrm{r}})^{2}\\
-1.7849\mathbf{I}+0.1525(2.3917\thinspace\mathbf{I}-\mathbf{A}_{\mathrm{r}})+0.1406(2.3917\thinspace\mathbf{I}-\mathbf{A}_{\mathrm{r}})^{2}
\end{array}\right]\nonumber \\
 & =\left[\begin{array}{c}
-5.9963\mathbf{I}+6.1345\mathbf{A}+0.875\mathbf{A}^{2}\\
0.491\mathbf{I}-1.5506\mathbf{A}+0.5625\mathbf{A}^{2}\\
-1.0885\mathbf{I}+0.4086\mathbf{A}+0.2187\mathbf{A}^{2}\\
-0.6158\mathbf{I}-0.8251\mathbf{A}+0.1406\mathbf{A}^{2}
\end{array}\right]\,,\label{eq:Mscheme-M2-Crk}
\end{align}
and, consequently, the coefficient matrix is stored as 
\begin{equation}
\mathbf{c}_{\mathrm{r}}=\left[\begin{array}{ccc}
-5.9963 & 6.1345 & 0.8750\\
0.4910 & -1.5506 & 0.5625\\
-1.0885 & 0.4086 & 0.2187\\
-0.6158 & -0.8251 & 0.1406
\end{array}\right]\,.\label{eq:Mscheme-M2-cr}
\end{equation}

The function $\texttt{InitSchemeRhoPF(3, 0.125)}$ in Fig.\,\ref{fig:MATLAB-code-initial_multiple}
then returns the root $r$ in Eq.\,(\ref{eq:Mscheme-M2-r}), the
coefficients $p_{\mathrm{r}i}$ of Eq.\,(\ref{eq:Mscheme-M2-Pr}),
and the coefficient matrix $\mathbf{c}_{\mathrm{r}}$ in Eq.\,(\ref{eq:Mscheme-M2-cr}),
as required.

\section{Sample MATLAB code for linear time-stepping solution algorithms}\label{sec:Sample-MATLAB-code}

Figure\,\ref{fig:MATLAB-function-TimeSolverPF} and Figure\,\ref{fig:MATLAB-function-TimeSolverTRhoPF}
provide sample MATLAB code for executing the algorithms presented
in this work for linear problems. These can be simply extended to
nonlinear problems by introducing an iterative solution procedure
as presented in Tables \ref{tab:Time-integration-solution_distinct}
and \ref{tab:Time-integration-solution_multiple}. Note that the functions
$\texttt{InitSchemePadePF}$ and $\texttt{InitSchemeRhoPF}$ used
to initialize the coefficients required for the time-stepping solution
for each procedure were provided in Fig.\,\ref{fig:MATLAB-code-initial_distinct}
and Fig.\,\ref{fig:MATLAB-code-initial_multiple}. The function \texttt{lu
}is used in Fig.\,\ref{fig:MATLAB-function-TimeSolverTRhoPF} for
the solution of complex equations since it is found to be faster than
the \texttt{decomposition} function in MATLAB 2023b. 

\begin{figure}

\texttt{function {[}RespHist{]} = TimeSolverPF(porder,signal,ns,dt,K,M,C,F,MLumped,u0,v0)}

\texttt{p = mod(porder(1),10);}

\texttt{rhoInfty = porder(2);}

\texttt{nCmplx = floor(p/2);}

\texttt{nReal = mod(p,2);}

\texttt{pf = min(p+1,size(cfr,1));}

\texttt{{[}rho, plr, rs, a, cfr{]} = InitSchemePadePF(p,rhoInfty,pf);}

\texttt{s = forceSamplingPoints(pf+0);}

\texttt{Tcfr = transMtxPointsToPoly(s, pf){*}cfr(1:pf,:);}

\texttt{if nReal > 0}

\texttt{r = rs(1);}

\texttt{dKd1 = decomposition(sparse((r{*}r){*}M + r{*}C + K));}

\texttt{end}

\texttt{if nCmplx > 0}

\texttt{L = cell(nCmplx,1); U = L; LUp = L; LUq = L;}

\texttt{for ic = 1:nCmplx}

\texttt{r = rs(ic+nReal);}

\texttt{{[}L\{ic\},U\{ic\},LUp\{ic\},LUq\{ic\}{]} = lu(sparse((r{*}r){*}M
+ r{*}C + K),'vector');}

\texttt{end}

\texttt{end}

\texttt{n = size(K,1);}

\texttt{tm = 0; z0 = {[}v0; u0{]}; \% Initial conditions}

\texttt{an = (F{*}signal(0)+ftmp)./MLumped;}

\texttt{for is= 1:ns-1}

\texttt{ts = tm + dt{*}s;}

\texttt{Fp = F.{*}reshape(signal(ts), 1, {[}{]});}

\texttt{tm = tm + dt;}

\texttt{z = rho{*}z0;}

\texttt{Az1 = rho{*}Az1;}

\texttt{if nReal > 0}

\texttt{r = real(rs(1));}

\texttt{g = real(plr(1)){*}z0;}

\texttt{fri = Fp{*}(real(Tcfr(:,1)));}

\texttt{ftmp = dKd1\textbackslash (r{*}(M{*}g(1:n)) - K{*}g(n+1:end)
+ r{*}fri);}

\texttt{y = {[}ftmp ; (ftmp+g(n+1:end))/r{]};}

\texttt{z = z + real(a(1)){*}y;}

\texttt{ftmp = r{*}y(1:n) - ( g(1:n) + fri./MLumped );}

\texttt{Az1 = Az1 + real(a(1)){*}ftmp;}

\texttt{end}

\caption{MATLAB function to solve the time-stepping equation for case of distinct
roots. }\label{fig:MATLAB-function-TimeSolverPF}
\end{figure}

\addtocounter{figure}{-1}
\begin{figure}

\texttt{if nCmplx > 0}

\texttt{for ic = 1:nCmplx}

\texttt{r = rs(ic+nReal);}

\texttt{cg = plr(ic+nReal){*}z0;}

\texttt{cfri = Fp{*}(Tcfr(:,ic+nReal));}

\texttt{ctmp = plr(ic+nReal){*}(r{*}(M{*}z0(1:n)) - K{*}z0(n+1:end))
+ r{*}cfri;}

\texttt{tmp(LUp\{ic\},:) = U\{ic\}\textbackslash (L\{ic\}\textbackslash ctmp(LUq\{ic\},:));}

\texttt{y = {[}tmp; (tmp + cg(n+1:end))/r {]};}

\texttt{z = z + 2{*}real(a(ic+nReal){*}y);}

\texttt{ctmp = r{*}y(1:n) - ( cg(1:n) + cfri./MLumped );}

\texttt{Az1 = Az1 + 2{*}real( a(ic+nReal){*}ctmp );}

\texttt{end}

\texttt{end}

\texttt{z0 = z;}

\texttt{an = Az1(1:n) + (F{*}signal(tm))./MLumped;}

\texttt{end}

\caption{MATLAB function to solve the time-stepping equation for case of distinct
roots (continued). }
\end{figure}
\begin{figure}

\texttt{function {[}RespHist{]} = TimeSolverTRhoPF(scheme,p,rho,signal,ns,dt,K,M,C,F,MLumped,u0,v0)}

\texttt{{[}r, prcoe, cfr1{]} = InitSchemeRhoPF(scheme, p, rho, pf);}

\texttt{pf = min(p+1,size(cfr1,1));}

\texttt{s = forceSamplingPoints(pf);}

\texttt{Tcfr1 = transMtxPointsToPoly(s, pf){*}cfr1(1:pf,:);}

\texttt{nz = size(K,1); \% number of DOFs}

\texttt{tm = 0; z = {[}v0; u0{]}; \% Initial conditions}

\texttt{dKd = decomposition(sparse((r{*}r){*}M + r{*}C + K));}

\texttt{ftmp = -(C{*}v0 + K{*}u0);}

\texttt{Az1 = ftmp./MLumped;}

\texttt{an = (F{*}signal(0)+ftmp)./MLumped;}

\texttt{for is = 1:ns-1 }

\texttt{ts = tm + dt{*}s;}

\texttt{Fp = F.{*}reshape(signal(ts), 1, {[}{]});}

\texttt{tm = tm + dt;}

\texttt{zi = zeros(2{*}nz,1);}

\texttt{for ip = 1:p}

\texttt{g = zi + prcoe(ip){*}z;}

\texttt{rfri = Fp{*}Tcfr1(:,ip);}

\texttt{zi(1:nz) = dKd\textbackslash (r{*}(M{*}g(1:nz)) - K{*}g(nz+1:end)
+ rfri);}

\texttt{zi(nz+1:end) = (zi(1:nz) + g(nz+1:end))/r;}

\texttt{end}

\texttt{z = prcoe(p+1){*}z + zi;}

\texttt{Az1 = prcoe(p+1){*}Az1 + r{*}zi(1:nz) - (g(1:nz) + (rfri/r)./MLumped);}

\texttt{an = Az1 + (F{*}signal(tm))./MLumped;}

\texttt{end}

\caption{MATLAB function to solve the time-stepping equation for case of single
multiple root. }\label{fig:MATLAB-function-TimeSolverTRhoPF}
\end{figure}

For each of the MATLAB codes in Figs.\,\ref{fig:MATLAB-function-TimeSolverPF}
and \ref{fig:MATLAB-function-TimeSolverTRhoPF}, the functions $\texttt{forceSamplingPoints}$
and $\texttt{transMtxPointsToPoly}$ not listed in this paper can
be found in the source code package at: <to be inserted once the manuscript
is accepted for publication>.

\end{document}